\documentclass[12pt,leqno]{article}
\usepackage{amssymb}
\usepackage[mathscr]{eucal}
\usepackage{amsmath,amssymb,latexsym,theorem,bbm}
\usepackage{amsmath,amssymb,latexsym,theorem}
\usepackage{color}
\usepackage{appendix}

\newcommand{\SC}{\scriptstyle}

\newcommand{\CC}{\mathbb{C}}
\newcommand{\DD}{\mathbb{D}}

\newcommand{\NN}{\mathbb{N}}

\newcommand{\RR}{\mathbb{R}}

\newcommand{\ZZ}{\mathbb{Z}}

\newcommand{\bA}{{\boldsymbol{A}}}

\newcommand{\bB}{{\boldsymbol{B}}}

\newcommand{\tb}{\widetilde{b}}
\newcommand{\htb}{\widehat{\tb}}

\newcommand{\tbB}{{\widetilde{\bB}}}
\newcommand{\tB}{\widetilde{B}}
\newcommand{\htbB}{\widehat{\tbB}}

\newcommand{\bc}{{\boldsymbol{c}}}
\newcommand{\bC}{{\boldsymbol{C}}}

\newcommand{\tbC}{\widetilde{\bC}}

\newcommand{\obC}{\overline{\bC}}

\newcommand{\bD}{{\boldsymbol{D}}}

\newcommand{\be}{{\boldsymbol{e}}}

\newcommand{\Bf}{{\boldsymbol{f}}}

\newcommand{\bF}{{\boldsymbol{F}}}
\newcommand{\bG}{{\boldsymbol{G}}}

\newcommand{\bh}{{\boldsymbol{h}}}
\newcommand{\bI}{{\boldsymbol{I}}}

\newcommand{\bK}{{\boldsymbol{K}}}
\newcommand{\bL}{{\boldsymbol{L}}}

\newcommand{\bM}{{\boldsymbol{M}}}

\newcommand{\bP}{{\boldsymbol{P}}}

\newcommand{\bcQ}{{\boldsymbol{\cQ}}}
\newcommand{\br}{{\boldsymbol{r}}}
\newcommand{\bR}{{\boldsymbol{R}}}
\newcommand{\tbR}{{\widetilde{\bR}}}

\newcommand{\bS}{{\boldsymbol{S}}}
\newcommand{\tbS}{\widetilde{\bS}}
\newcommand{\ttbS}{\widetilde{\tbS}}
\newcommand{\bT}{{\boldsymbol{T}}}
\newcommand{\tU}{\widetilde{U}}
\newcommand{\bu}{{\boldsymbol{u}}}
\newcommand{\tbu}{\widetilde{\bu}}
\newcommand{\bv}{{\boldsymbol{v}}}
\newcommand{\bV}{{\boldsymbol{V}}}

\newcommand{\bx}{{\boldsymbol{x}}}
\newcommand{\bX}{{\boldsymbol{X}}}
\newcommand{\by}{{\boldsymbol{y}}}
\newcommand{\bY}{{\boldsymbol{Y}}}

\newcommand{\tcY}{\widetilde{\cY}}

\newcommand{\bz}{{\boldsymbol{z}}}
\newcommand{\bZ}{{\boldsymbol{Z}}}
\newcommand{\bU}{{\boldsymbol{U}}}
\newcommand{\bw}{{\boldsymbol{w}}}
\newcommand{\bW}{{\boldsymbol{W}}}
\newcommand{\Bbeta}{{\boldsymbol{\beta}}}
\newcommand{\oBbeta}{\overline{\Bbeta}}
\newcommand{\obeta}{\overline{\beta}}
\newcommand{\hoBbeta}{\widehat{\oBbeta}}
\newcommand{\tbeta}{\widetilde{\beta}}
\newcommand{\htbeta}{\widehat{\tbeta}}
\newcommand{\tBbeta}{\widetilde{\Bbeta}}
\newcommand{\htBbeta}{\widehat{\tBbeta}}
\newcommand{\tdelta}{\widetilde{\delta}}
\newcommand{\bgamma}{{\boldsymbol{\gamma}}}

\newcommand{\blambda}{{\boldsymbol{\lambda}}}

\newcommand{\tvarrho}{\widetilde{\varrho}}

\newcommand{\bxi}{{\boldsymbol{\xi}}}
\newcommand{\bXi}{{\boldsymbol{\Xi}}}
\newcommand{\bmu}{{\boldsymbol{\mu}}}
\newcommand{\balpha}{{\boldsymbol{\alpha}}}
\newcommand{\bbeta}{{\boldsymbol{\beta}}}
\newcommand{\Beta}{{\boldsymbol{\eta}}}
\newcommand{\btheta}{{\boldsymbol{\theta}}}

\newcommand{\bPi}{{\boldsymbol{\Pi}}}
\newcommand{\bSigma}{{\boldsymbol{\Sigma}}}

\newcommand{\bzero}{{\boldsymbol{0}}}

\newcommand{\cA}{{\mathcal A}}
\newcommand{\cB}{{\mathcal B}}

\newcommand{\cD}{{\mathcal D}}

\newcommand{\cI}{\mathcal{I}}
\newcommand{\cJ}{\mathcal{J}}

\newcommand{\cF}{{\mathcal F}}
\newcommand{\cH}{{\mathcal H}}

\newcommand{\cM}{{\mathcal M}}
\newcommand{\bcM}{\boldsymbol{\cM}}

\newcommand{\cN}{{\mathcal N}}
\newcommand{\bcN}{\boldsymbol{\cN}}
\newcommand{\cP}{{\mathcal P}}
\newcommand{\bcP}{\boldsymbol{\cP}}
\newcommand{\cQ}{{\mathcal Q}}
\newcommand{\cR}{{\mathcal R}}
\newcommand{\cS}{{\mathcal S}}
\newcommand{\cU}{{\mathcal U}}

\newcommand{\bcU}{\boldsymbol{\cU}}
\newcommand{\cV}{{\mathcal V}}

\newcommand{\cX}{{\mathcal X}}

\newcommand{\bcX}{\boldsymbol{\cX}}
\newcommand{\tbcX}{\widetilde{\bcX}}

\newcommand{\cY}{{\mathcal Y}}
\newcommand{\cZ}{{\mathcal Z}}
\newcommand{\cW}{{\mathcal W}}
\newcommand{\bcW}{\boldsymbol{\cW}}
\newcommand{\tbcW}{\widetilde{\bcW}}
\newcommand{\bcY}{\boldsymbol{\cY}}

\newcommand{\bcZ}{\boldsymbol{\cZ}}
\newcommand{\tcW}{\widetilde{\cW}}
\newcommand{\ttcW}{\widetilde{\tcW}}

\newcommand{\cc}{\mathrm{c}}
\newcommand{\dd}{\mathrm{d}}
\newcommand{\ee}{\mathrm{e}}

\newcommand{\slu}{{\SC\mathrm{lu}}}

\newcommand{\INARp}{\textup{INAR($p$)}}

\newcommand{\EE}{\operatorname{\mathbb{E}}}

\newcommand{\PP}{\operatorname{\mathbb{P}}}
\newcommand{\OO}{\operatorname{O}}

\newcommand{\var}{\operatorname{Var}}

\renewcommand{\Re}{\operatorname{Re}}
\newcommand{\argmin}{\operatorname{arg\,min}}

\newcommand{\hgamma}{\widehat{\gamma}}
\newcommand{\hdelta}{\widehat{\delta}}

\newcommand{\hkappa}{\widehat{\kappa}}
\newcommand{\hvarrho}{\widehat{\varrho}}

\newcommand{\hs}{\widehat{s}}

\newcommand{\tH}{\widetilde{H}}

\newcommand{\tN}{\widetilde{N}}

\newcommand{\vare}{\varepsilon}

\renewcommand{\mid}{\,|\,}
\newcommand{\bmid}{\,\big|\,}

\renewcommand{\leq}{\leqslant}
\renewcommand{\geq}{\geqslant}

\newcommand{\stoch}{\stackrel{\PP}{\longrightarrow}}
\newcommand{\distr}{\stackrel{\cD}{\longrightarrow}}
\newcommand{\distre}{\stackrel{\cD}{=}}

\newcommand{\lu}{\stackrel{\slu}{\longrightarrow}}

\newcommand{\ase}{\stackrel{{\mathrm{a.s.}}}{=}}

\newcommand{\bbone}{\mathbbm{1}}
\newcommand{\ns}{{\lfloor ns\rfloor}}
\newcommand{\nt}{{\lfloor nt\rfloor}}
\newcommand{\nT}{{\lfloor nT\rfloor}}
\newcommand{\proofend}{\hfill\mbox{$\Box$}}

\numberwithin{equation}{section}

\theoremstyle{change} \theorembodyfont{\em}
\newtheorem{Lem}{Lemma.}[section]
\newtheorem{Thm}[Lem]{Theorem.}
\newtheorem{Pro}[Lem]{Proposition.}
\newtheorem{Cor}[Lem]{Corollary.}
\newtheorem{Def}[Lem]{Definition.}

\theorembodyfont{\rm}
\newtheorem{Rem}[Lem]{Remark.}

\begin{document}

\begin{center}
 {\bfseries\Large
   Statistical inference for 2-type doubly symmetric critical
    irreducible continuous state and continuous time \\[2mm]
    branching processes with immigration}

\vspace*{3mm}

 {\sc\large
  M\'aty\'as $\text{Barczy}^{*,\diamond}$,
  \ Krist\'of $\text{K\"ormendi}^{**}$,
  \ Gyula $\text{Pap}^{***}$}

\end{center}

\vskip0.2cm

\noindent
 * Faculty of Informatics, University of Debrecen,
   Pf.~12, H--4010 Debrecen, Hungary.

\noindent
 ** MTA-SZTE Analysis and Stochastics Research Group,
     Bolyai Institute, University of Szeged,
     Aradi v\'ertan\'uk tere 1, H--6720 Szeged, Hungary.

\noindent
 *** Bolyai Institute, University of Szeged,
     Aradi v\'ertan\'uk tere 1, H--6720 Szeged, Hungary.

\noindent e--mails: barczy.matyas@inf.unideb.hu (M. Barczy),
                    kormendi@math.u-szeged.hu (K. K\"ormendi),
                    papgy@math.u-szeged.hu (G. Pap).

\noindent $\diamond$ Corresponding author.

\vskip0.2cm


\renewcommand{\thefootnote}{}
\vspace*{0.2cm}
\footnote{The research of M. Barczy and G. Pap was realized in the frames of
 T\'AMOP 4.2.4.\ A/2-11-1-2012-0001 ,,National Excellence Program --
 Elaborating and operating an inland student and researcher personal support
 system''.
The project was subsidized by the European Union and co-financed by the
 European Social Fund.}

\vspace*{-10mm}

\begin{abstract}
We study asymptotic behavior of conditional least squares estimators for
 2-type doubly symmetric critical irreducible continuous state and continuous
 time branching processes with immigration based on discrete time (low
 frequency) observations.
\end{abstract}

\section{Introduction}
\label{section_intro}

Asymptotic behavior of conditional least squares (CLS) estimators for critical
 continuous state and continuous time branching processes with immigration
 (CBI processes) is available only for single-type processes.
Huang et al.\ \cite{HuaMaZhu} considered a single-type CBI process which can
 be represented as a pathwise unique strong solution of the stochastic
 differential equation (SDE)
 \begin{align}\label{SDE_atirasa_dim1}
  \begin{split}
   X_t
   &=X_0
     + \int_0^t (\beta + \tb X_s) \, \dd s
     + \int_0^t \sqrt{2 c \max \{0, X_s\}} \, \dd W_s \\
   &\quad
      + \int_0^t \int_0^\infty \int_0^\infty
         z \bbone_{\{u\leq X_{s-}\}} \, \tN(\dd s, \dd z, \dd u)
      + \int_0^t \int_0^\infty z \, M(\dd s, \dd z)
  \end{split}
 \end{align}
 for \ $t \in [0, \infty)$, \ where \ $\beta, c \in [0, \infty)$,
 \ $\tb \in \RR$, \ and \ $(W_t)_{t\geq0}$ \ is a standard Wiener process,
 \ $N$ \ and \ $M$ \ are independent Poisson random measures on
 \ $(0, \infty)^3$ \ and on \ $(0, \infty)^2$ \ with intensity measures
 \ $\dd s \, \mu(\dd z) \, \dd u$ \ and \ $\dd s \, \nu(\dd z)$,
 \ respectively,
 \ $\tN(\dd s, \dd z, \dd u)
    := N(\dd s, \dd z, \dd u) - \dd s \, \mu(\dd z) \, \dd u$
 \ is the compensated Poisson random measure
    corresponding to \ $N$,
 \ the measures \ $\mu$ \ and \ $\nu$ \ satisfy some moment conditions, and
 \ $(W_t)_{t\geq0}$, \ $N$ \ and \ $M$ \ are independent.
The model is called subcritical, critical or supercritical if \ $\tb < 0$,
 \ $\tb = 0$ \ or \ $\tb > 0$, \ see Huang et al.\ \cite[page 1105]{HuaMaZhu}
 or Definition \ref{Def_indecomposable_crit}.
Based on discrete time (low frequency) observations
 \ $(X_k)_{k \in \{0, 1, \ldots, n\}}$, \ $n \in \{1, 2, \ldots\}$,
 \ Huang et al.\ \cite{HuaMaZhu} derived weighted CLS estimator of
 \ $(\tb, \beta)$.
\ Under some second order moment assumptions, supposing that \ $c$, \ $\mu$ \ and
 \ $\nu$ \ are known, they showed the following results:
 in the subcritical case the estimator of \ $(\tb, \beta)$ \ is asymptotically normal;
 in the critical case the estimator of \ $\tb$ \ has a non-normal limit, but
 the asymptotic behavior of the estimator of \ $\beta$ \ remained open; in the
 supercritical case the estimator of \ $\tb$ \ is asymptotically normal with a
 random scaling, but the estimator of \ $\beta$ \ is not weakly consistent.

Based on the observations \ $(X_k)_{k \in \{0, 1, \ldots, n\}}$,
 \ $n \in \{1, 2, \ldots\}$, \ supposing that \ $c$, \ $\mu$ \ and
 \ $\nu$ \ are known, Barczy et al. \cite{BarKorPap} derived (non-weighted)
 CLS estimator \ $(\htb_n, \htbeta_n)$, \  of \ $(\tb, \tbeta)$, \ where
 \ $\tbeta := \beta + \int_0^\infty z \, \nu(\dd z)$.
\ In the critical case, under some moment assumptions, it has been shown that
 \ $\bigl(n (\htb_n - \tb), \htbeta_n - \tbeta\bigr)$ \ has a non-normal limit.
As a by-product, the estimator \ $\htbeta_n$ \ is not weakly consistent.

Overbeck and Ryd\'en \cite{OveRyd} considered CLS and weighted CLS estimators
 for the well-known Cox--Ingersoll--Ross model, which is, in fact, a
 single-type diffusion CBI process (without jump part), i.e., when \ $\mu = 0$
 \ and \ $\nu = 0$ \ in \eqref{SDE_atirasa_dim1}.
Based on discrete time observations
 \ $(X_k)_{k\in\{0,1,\ldots,n\}}$, \ $n \in \{1, 2, \ldots\}$, \ they derived CLS
 estimator of \ $(\tb, \beta, c)$ \ and proved its asymptotic normality in the
 subcritical case.
Note that Li and Ma \cite{LiMa} started to investigate the asymptotic
 behaviour of the CLS and weighted CLS estimators of the parameters
 \ $(\tb, \beta)$ \ in the subcritical case for a Cox--Ingersoll--Ross model
 driven by a stable noise, which is again a special single-type CBI process
 (with jump part).

In this paper we consider a 2-type CBI process which can be represented as a
 pathwise unique strong solution of the SDE
 \begin{align}\label{SDE_atirasa_dim2}
  \begin{split}
   \bX_t
   &=\bX_0
     + \int_0^t (\Bbeta + \tbB \bX_s) \, \dd s
     + \sum_{i=1}^2 \int_0^t \sqrt{2 c_i \max \{0, X_{s,i}\}} \, \dd W_{s,i} \, \be_i \\
   &\quad
      + \sum_{j=1}^2
         \int_0^t \int_{\cU_2} \int_0^\infty
          \bz \bbone_{\{u\leq X_{s-,j}\}} \, \tN_j(\dd s, \dd\bz,  \dd u)
      + \int_0^t \int_{\cU_2} \bz \, M(\dd s, \dd\bz)
  \end{split}
 \end{align}
 for \ $t \in [0, \infty)$.
 \ Here \ $X_{t,i}$, \ $i \in \{1, 2\}$, \ denotes the coordinates of \ $\bX_t$,
 \ $\Bbeta \in [0, \infty)^2$, \ $\tbB \in \RR^{2\times2}$ \ has
 non-negative off-diagonal entries, \ $c_1, c_2 \in [0, \infty)$,
 \ $\be_1, \be_2$ \ denotes the natural basis in \ $\RR^2$,
 \ $\cU_2 := [0, \infty)^2 \setminus \{(0, 0)\}$, \ $(W_{t,1})_{t\geq0}$ \ and
 \ $(W_{t,2})_{t\geq0}$ \ are independent standard Wiener processes, \ $N_j$,
 \ $j \in \{1, 2\}$, \ and \ $M$ \ are independent Poisson random measures on
 \ $(0, \infty) \times \cU_2 \times (0, \infty)$ \ and on
 \ $(0, \infty) \times \cU_2$ \ with intensity measures
 \ $\dd s \, \mu_j(\dd\bz) \, \dd u$, \ $j \in \{1, 2\}$, \ and
 \ $\dd s \, \nu(\dd\bz)$, \ respectively,
 \ $\tN_j(\dd s, \dd\bz, \dd u)
    := N_j(\dd s, \dd\bz, \dd u) - \dd s \, \mu_j(\dd\bz) \, \dd u$,
 \ $j \in \{1, 2\}$.
\ We suppose that the Borel measures \ $\mu_j$, \ $j \in \{1, 2\}$,
 \ and \ $\nu$ \ on \ $\cU_2$ \ satisfy some moment conditions, and
 \ $(W_{t,1})_{t\geq0}$, \ $(W_{t,2})_{t\geq0}$, \ $N_1$, \ $N_2$ \ and \ $M$
 \ are independent.
We will suppose that the process \ $(\bX_t)_{t\geq0}$ \ is doubly symmetric in
 the sense that
 \[
   \tbB = \begin{bmatrix}
           \gamma & \kappa \\
           \kappa & \gamma
          \end{bmatrix} ,
 \]
 where \ $\gamma \in \RR$ \ and \ $\kappa \in [0, \infty)$.
\ Note that the parameters \ $\gamma$ \ and \ $\kappa$ \ might be interpreted
 as the transformation rates of one type to the same type and one type to the
 other type, respectively, compare with Xu \cite{Xu}; that's why the model can
 be called doubly symmetric.

The model will be called subcritical, critical or supercritical if
 \ $s < 0$, \ $s = 0$ \ or \ $s > 0$, \ respectively, where
 \ $s := \gamma + \kappa$ \ denotes the criticality parameter, see Definition
 \ref{Def_indecomposable_crit}.

For the simplicity, we suppose \ $\bX_0 = (0, 0)^\top$.
\ We suppose that \ $c_1, c_2, \mu_1, \mu_2$ \ and \ $\nu$ \ are known,
 and we derive the CLS estimators of the parameters \ $s$, \ $\gamma$,
 \ $\kappa$ \ and \ $\Bbeta$ \ based on discrete time observations
 \ $(\bX_k)_{k \in \{1, \ldots, n\}}$, \ $n \in \{1, 2, \ldots\}$.
In the irreducible and critical case, i.e, when \ $\kappa > 0$ \ and
 \ $s = \gamma + \kappa = 0$, \ under some moment conditions, we describe the
 asymptotic behavior of these CLS estimators as \ $n \to \infty$, \ provided
 that \ $\Bbeta \ne (0, 0)^\top$ \ or \ $\nu \ne 0$, \ see Theorem \ref{main}.
We point out that the limit distributions are non-normal in general.
In the present paper we do not investigate the asymptotic behavior of CLS
 estimators of  \ $s$, \ $\gamma$, \ $\kappa$ \ and \ $\Bbeta$ \ in the
 subcritical and supercritical cases, it could be the topic of separate papers,
 needing different approaches.

Xu \cite{Xu} considered a 2-type diffusion CBI process (without jump part),
 i.e., when \ $\mu_j = 0$, \ $j \in \{1, 2\}$, \ and \ $\nu = 0$ \ in
 \eqref{SDE_atirasa_dim2}.
Based on discrete time (low frequency) observations
 \ $(\bX_k)_{k\in\{1,\ldots,n\}}$, \ $n \in \{1, 2, \ldots\}$, \ Xu \cite{Xu}
 derived CLS estimators and weighted CLS estimators of
 \ $(\Bbeta, \tbB, c_1, c_2)$.
\ Provided that \ $\Bbeta \in (0, \infty)^2$, \ the diagonal entries of
 \ $\tbB$ \ are negative, the off-diagonal entries of \ $\tbB$ \ are positive,
 \ the determinant of \ $\tbB$ \ is positive and \ $c_i > 0$,
 \ $i \in \{1, 2\}$ \ (which yields that the process \ $\bX$ \ is irreducible
 and subcritical, see Xu \cite[Theorem 2.2]{Xu} and Definitions
 \ref{Def_irreducible} and \ref{Def_indecomposable_crit}), it was shown that
 these CLS estimators are asymptotically normal, see Theorem 4.6 in Xu
 \cite{Xu}.

Finally, we give an overview of the paper.
In Section \ref{section_CBI}, for completeness and better readability,
 from Barczy et al.~\cite{BarLiPap2} and \cite{BarPap}, we recall some notions
 and statements for multi-type CBI processes such as the form of their
 infinitesimal generator, Laplace transform, a formula for their first moment,
 the definition of subcritical, critical and supercritical irreducible CBI
 processes, see Definitions \ref{Def_irreducible} and \ref{Def_indecomposable_crit}.
We recall a result due to Barczy and Pap \cite[Theorem 4.1]{BarPap} stating that,
 under some fourth order moment assumptions, a sequence of scaled random step
 functions \ $(n^{-1}\bX_{\lfloor nt\rfloor})_{t\geq 0}$, $n\geq 1$, \ formed from
 a critical, irreducible multi-type CBI process \ $\bX$ \ converges weakly towards
 a squared Bessel process supported by a ray determined by the Perron vector of a
 matrix related to the branching mechanism of \ $\bX$.

In Section \ref{section_CBI_2}, first we derive formulas of CLS estimators of the
 transformed parameters \ $\ee^{\gamma+\kappa}$, \ $\ee^{\gamma-\kappa}$ \ and
 \ $\int_0^1 \ee^{s\tbB} \tBbeta \, \dd s$, \ and then of the parameters
 \ $\gamma$, \ $\kappa$ \ and \ $\tBbeta$.
\ The reason for this parameter transformation is to reduce the minimization in
 the CLS method to a linear problem.
Then we formulate our main result about the asymptotic behavior of CLS estimators
 of \ $s$, \ $\gamma$, \ $\kappa$ \ and \ $\tBbeta$ \ in the irreducible and
 critical case, see Theorem \ref{main}.
These results will be derived from the corresponding statements for the
 transformed parameters, see Theorem \ref{main_rdb}.

In Section \ref{section_deco}, we give a decomposition of the process \ $\bX$ \
 and of the CLS estimators of the transformed parameters as
 well, related to the left eigenvectors of \ $\tbB$ \ belonging to the
 eigenvalues \ $\gamma + \kappa$ \ and \ $\gamma - \kappa$, \ see formulas
 \eqref{XUV}, \eqref{r-}, \eqref{d-} and \eqref{b-}.
By the help of these decompositions, Theorem \ref{main_rdb} will follow from
 Theorems \ref{main_Ad}, \ref{maint_Ad} and \ref{main1_Ad}.

Sections \ref{section_proof_main}, \ref{section_proof_maint} and
 \ref{section_proof_main1} are devoted to the proofs of Theorems \ref{main_Ad},
 \ref{maint_Ad} and \ref{main1_Ad}, respectively.
The proofs are heavily based on a careful analysis of the asymptotic behavior of
 some martingale differences related to the process \ $\bX$ \ and the
 decompositions given in Section \ref{section_deco}, and delicate moment
 estimations for the process \ $\bX$ \ and some auxiliary processes.

In Appendix \ref{section_SDE} we recall a representation of multi-type CBI
 processes as pathwise unique strong solutions of certain SDEs with jumps based
 on Barczy et al.\ \cite{BarLiPap2}.
In Appendix \ref{section_moments} we recall some results about the asymptotic
 behaviour of moments of irreducible and critical multi-type CBI processes based
 on Barczy, Li and Pap \cite{BarLiPap3}, and then, presenting new results as well,
 the asymptotic behaviour of the moments of some auxiliary processes is also
 investigated.
Appendix \ref{section_estimators} is devoted to study of the existence of the CLS
 estimator of the transformed parameters.
In Appendix \ref{CMT}, we present a version of the continuous mapping theorem.
In Appendix \ref{section_conv_step_processes}, we recall a useful result about
 convergence of random step processes towards a diffusion process due to Isp\'any
 and Pap \cite[Corollary 2.2]{IspPap2}.

\section{Multi-type CBI processes}
\label{section_CBI}

Let \ $\ZZ_+$, \ $\NN$, \ $\RR$, \ $\RR_+$  \ and \ $\RR_{++}$ \ denote the set
 of non-negative integers, positive integers, real numbers, non-negative real
 numbers and positive real numbers, respectively.
For \ $x , y \in \RR$, \ we will use the notations
 \ $x \land y := \min \{x, y\} $ \ and \ $x^+:= \max \{0, x\}$.
\ By \ $\|\bx\|$ \ and \ $\|\bA\|$, \ we denote the Euclidean norm of a vector
 \ $\bx \in \RR^d$ \ and the induced matrix norm of a matrix
 \ $\bA \in \RR^{d\times d}$, \ respectively.
The natural basis in \ $\RR^d$ \ will be denoted by \ $\be_1$, \ldots, $\be_d$.
\ The null vector and the null matrix will be denoted by \ $\bzero$.
\ By \ $C^2_\cc(\RR_+^d,\RR)$ \ we denote the set of twice continuously
 differentiable real-valued functions on \ $\RR_+^d$ \ with compact support.
Convergence in distribution and in probability will be denoted by \ $\distr$
 \ and \ $\stoch$, \ respectively.
Almost sure equality will be denoted by \ $\ase$.

\begin{Def}\label{Def_essentially_non-negative}
A matrix \ $\bA = (a_{i,j})_{i,j\in\{1,\ldots,d\}} \in \RR^{d\times d}$ \ is called
 essentially non-negative if \ $a_{i,j} \in \RR_+$ \ whenever
 \ $i, j \in \{1,\ldots,d\}$ \ with \ $i \ne j$, \ that is, if \ $\bA$ \ has
 non-negative off-diagonal entries.
The set of essentially non-negative \ $d \times d$ \ matrices will be denoted
 by \ $\RR^{d\times d}_{(+)}$.
\end{Def}

\begin{Def}\label{Def_admissible}
A tuple \ $(d, \bc, \Bbeta, \bB, \nu, \bmu)$ \ is called a set of admissible
 parameters if
 \renewcommand{\labelenumi}{{\rm(\roman{enumi})}}
 \begin{enumerate}
  \item
   $d \in \NN$,
  \item
   $\bc = (c_i)_{i\in\{1,\ldots,d\}} \in \RR_+^d$,
  \item
   $\Bbeta = (\beta_i)_{i\in\{1,\ldots,d\}} \in \RR_+^d$,
  \item
   $\bB = (b_{i,j})_{i,j\in\{1,\ldots,d\}} \in \RR^{d \times d}_{(+)}$,
  \item
   $\nu$ \ is a Borel measure on \ ${\cU_d := \RR_+^d \setminus \{\bzero\}}$
    \ satisfying \ $\int_{\cU_d} (1 \land \|\bz\|) \, \nu(\dd \bz) < \infty$,
  \item
   $\bmu = (\mu_1, \ldots, \mu_d)$, \ where, for each
    \ $i \in \{1, \ldots, d\}$, \ $\mu_i$ \ is a Borel measure on
    \ $\cU_d$ \ satisfying
    \[
      \int_{\cU_d}
       \left[ \|\bz\| \wedge \|\bz\|^2
              + \sum_{j \in \{1, \ldots, d\} \setminus \{i\}} (1 \land z_j) \right]
       \mu_i(\dd \bz)
      < \infty .
    \]
  \end{enumerate}
\end{Def}

\begin{Rem}
Our Definition \ref{Def_admissible} of the set of admissible parameters is a
 special case of Definition 2.6 in Duffie et al.~\cite{DufFilSch}, which is
 suitable for all affine processes, see
 Barczy et al.~\cite[Remark 2.3]{BarLiPap2}.
\proofend
\end{Rem}

\begin{Thm}\label{CBI_exists}
Let \ $(d, \bc, \Bbeta, \bB, \nu, \bmu)$ \ be a set of admissible parameters.
Then there exists a unique conservative transition semigroup
 \ $(P_t)_{t\in\RR_+}$ \ acting on the Banach space (endowed with the supremum norm)
 of real-valued bounded
 Borel-measurable functions on the state space \ $\RR_+^d$ \ such that its
 infinitesimal generator is
 \begin{equation}\label{CBI_inf_gen}
  \begin{aligned}
   (\cA f)(\bx)
   &= \sum_{i=1}^d c_i x_i f_{i,i}''(\bx)
      + \langle \Bbeta + \bB \bx, \Bf'(\bx) \rangle
      + \int_{\cU_d} \bigl( f(\bx + \bz) - f(\bx) \bigr) \, \nu(\dd \bz) \\
   &\phantom{\quad}
      + \sum_{i=1}^d
         x_i
         \int_{\cU_d}
          \bigl( f(\bx + \bz) - f(\bx) - f'_i(\bx) (1 \land z_i) \bigr)
          \, \mu_i(\dd \bz)
  \end{aligned}
 \end{equation}
 for \ $f \in C^2_\cc(\RR_+^d,\RR)$ \ and \ $\bx \in \RR_+^d$, \ where \ $f_i'$
\ and \ $f_{i,i}''$, \ $i \in \{1, \ldots, d\}$, \ denote the first and second
 order partial derivatives of \ $f$ \ with respect to its \ $i$-th variable,
 respectively, and \ $\Bf'(\bx) := (f_1'(\bx), \ldots, f_d'(\bx))^\top$.
\ Moreover, the Laplace transform of the transition semigroup
 \ $(P_t)_{t\in\RR_+}$ \ has a representation
 \begin{align*}
  \int_{\RR_+^d} \ee^{- \langle \blambda, \by \rangle} P_t(\bx, \dd \by)
  = \ee^{- \langle \bx, \bv(t, \blambda) \rangle - \int_0^t \psi(\bv(s, \blambda)) \, \dd s} ,
  \qquad \bx \in \RR_+^d, \quad \blambda \in \RR_+^d , \quad t \in \RR_+ ,
 \end{align*}
 where, for any \ $\blambda \in \RR_+^d$, \ the continuously differentiable
 function
 \ $\RR_+ \ni t \mapsto \bv(t, \blambda)
    = (v_1(t, \blambda), \ldots, v_d(t, \blambda))^\top \in \RR_+^d$
 \ is the unique locally bounded solution to the system of differential
 equations
 \begin{equation}\label{EES}
   \partial_t v_i(t, \blambda) = - \varphi_i(\bv(t, \blambda)) , \qquad
   v_i(0, \blambda) = \lambda_i , \qquad i \in \{1, \ldots, d\} ,
 \end{equation}
 with
 \[
   \varphi_i(\blambda)
   := c_i \lambda_i^2 -  \langle \bB \be_i, \blambda \rangle
      + \int_{\cU_d}
         \bigl( \ee^{- \langle \blambda, \bz \rangle} - 1
                + \lambda_i (1 \land z_i) \bigr)
         \, \mu_i(\dd \bz)
 \]
 for \ $\blambda \in \RR_+^d$, \ $i \in \{1, \ldots, d\}$, \ and
 \[
   \psi(\blambda)
   := \langle \bbeta, \blambda \rangle
      + \int_{\cU_d}
         \bigl( 1 - \ee^{- \langle \blambda, \bz \rangle} \bigr)
         \, \nu(\dd \bz) , \qquad
   \blambda \in \RR_+^d .
 \]
\end{Thm}

\begin{Rem}
This theorem is a special case of Theorem 2.7 of Duffie et
 al.~\cite{DufFilSch} with \ $m = d$, \ $n = 0$ \ and zero killing rate.
The unique existence of a locally bounded solution to the system of
 differential equations \eqref{EES} is proved by Li \cite[page 45]{Li}.
\proofend
\end{Rem}

\begin{Def}\label{Def_CBI}
A conservative Markov process with state space \ $\RR_+^d$ \ and with
 transition semigroup
 \ $(P_t)_{t\in\RR_+}$ \ given in Theorem \ref{CBI_exists} is called a multi-type
 CBI process with parameters \ $(d, \bc, \Bbeta, \bB, \nu, \bmu)$.
\ The function
 \ $\RR_+^d \ni \blambda
    \mapsto (\varphi_1(\blambda), \ldots, \varphi_d(\blambda))^\top \in \RR^d$
 \ is called its branching mechanism, and the function
 \ $\RR_+^d \ni \blambda \mapsto \psi(\blambda) \in \RR_+$ \ is called its
 immigration mechanism.
\end{Def}

Note that the branching mechanism depends only on the parameters \ $\bc$,
 $\bB$ \ and \ $\bmu$, \ while the immigration mechanism depends only on the
 parameters \ $\Bbeta$ \ and \ $\nu$.

Let \ $(\bX_t)_{t\in\RR_+}$ \ be a multi-type CBI process with parameters
 \ $(d, \bc, \Bbeta, \bB, \nu, \bmu)$ \ such that
 \ $\EE(\|\bX_0\|) < \infty$ \ and the moment condition
 \begin{equation}\label{moment_condition_m_new}
  \int_{\cU_d} \|\bz\| \bbone_{\{\|\bz\|\geq1\}} \, \nu(\dd \bz) < \infty
 \end{equation}
 holds.

 Then, by formula (3.4) in Barczy et al. \cite{BarLiPap2},
 \begin{equation}\label{EXcond}
  \EE(\bX_t \mid \bX_0 = \bx)
  = \ee^{t\tbB} \bx + \int_0^t \ee^{u\tbB} \tBbeta \, \dd u ,
  \qquad \bx \in \RR_+^d , \quad t \in \RR_+ ,
 \end{equation}
 where
 \begin{gather}
  \tbB := (\tb_{i,j})_{i,j\in\{1,\ldots,d\}} , \qquad
  \tb_{i,j} := b_{i,j} + \int_{\cU_d} (z_i - \delta_{i,j})^+ \, \mu_j(\dd \bz) ,
  \label{tbB} \\
  \tBbeta := \Bbeta + \int_{\cU_d} \bz \, \nu(\dd \bz) ,
  \label{tBbeta}
 \end{gather}
 with \ $\delta_{i,j}:=1$ \ if \ $i=j$, \ and \ $\delta_{i,j}:=0$ \ if
 \ $i \ne j$.
\ Note that \ $\tbB \in \RR^{d \times d}_{(+)}$ \ and \ $\tBbeta \in \RR_+^d$,
 \ since
 \begin{equation}\label{help}
  \int_{\cU_d} \|\bz\| \, \nu(\dd\bz) < \infty , \qquad
  \int_{\cU_d} (z_i - \delta_{i,j})^+ \, \mu_j(\dd \bz) < \infty , \quad
  i, j \in \{1, \ldots, d\} ,
 \end{equation}
 see Barczy et al. \cite[Section 2]{BarLiPap2}.
One can give probabilistic interpretations of the modified parameters
 \ $\tbB$ \ and \ $\tBbeta$, \ namely,
 \ $\ee^\tbB \be_j = \EE(\bY_1 \mid \bY_0 = \be_j)$, \ $j \in \{1, \ldots, d\}$,
 \ and \ $\tBbeta = \EE(\bZ_1 \mid \bZ_0 = \bzero)$, \ where
 \ $(\bY_t)_{t\in\RR_+}$ \ and \ $(\bZ_t)_{t\in\RR_+}$ \ are multi-type CBI
 processes with parameters \ $(d, \bc, \bzero, \bB, 0, \bmu)$ \ and
 \ $(d, \bzero, \Bbeta, \bzero, \nu, \bzero)$, \ respectively, see formula
 \eqref{EXcond}.
The processes \ $(\bY_t)_{t\in\RR_+}$ \ and \ $(\bZ_t)_{t\in\RR_+}$ \ can be
 considered as pure branching (without immigration) and pure immigration
 (without branching) processes, respectively.
Consequently, \ $\ee^\tbB$ \ and \ $\tBbeta$ \ may be called the branching mean
 matrix and the immigration mean vector, respectively.

Next we recall a classification of multi-type CBI processes.
For a matrix \ $\bA \in \RR^{d \times d}$, \ $\sigma(\bA)$ \ will denote the
 spectrum of \ $\bA$, \ that is, the set of the eigenvalues of \ $\bA$.
\ Then \ $r(\bA) := \max_{\lambda \in \sigma(\bA)} |\lambda|$ \ is the spectral
 radius of \ $\bA$.
\ Moreover, we will use the notation
 \[
   s(\bA) := \max_{\lambda \in \sigma(\bA)} \Re(\lambda) .
 \]
A matrix \ $\bA \in \RR^{d\times d}$ \ is called reducible if there exist a
 permutation matrix \ $\bP \in \RR^{ d \times d}$ \ and an integer \ $r$ \ with
 \ $1 \leq r \leq d-1$ \ such that
 \[
  \bP^\top \bA \bP
   = \begin{bmatrix} \bA_1 & \bA_2 \\ \bzero & \bA_3 \end{bmatrix},
 \]
 where \ $\bA_1 \in \RR^{r\times r}$, \ $\bA_3 \in \RR^{ (d-r) \times (d-r) }$,
 \ $\bA_2 \in \RR^{ r \times (d-r) }$, \ and \ $\bzero \in \RR^{ (d-r) \times r}$
 \ is a null matrix.
A matrix \ $\bA \in \RR^{d\times d}$ \ is called irreducible if it is not
 reducible, see, e.g.,
 Horn and Johnson \cite[Definitions 6.2.21 and 6.2.22]{HorJoh}.
We do emphasize that no 1-by-1 matrix is reducible.

\begin{Def}\label{Def_irreducible}
Let \ $(\bX_t)_{t\in\RR_+}$ \ be a multi-type CBI process with parameters
 \ $(d, \bc, \Bbeta, \bB, \nu, \bmu)$ \ such that the moment condition
 \eqref{moment_condition_m_new} holds.
Then \ $(\bX_t)_{t\in\RR_+}$ \ is called irreducible if \ $\tbB$ \ is
 irreducible.
\end{Def}

\begin{Def}\label{Def_indecomposable_crit}
Let \ $(\bX_t)_{t\in\RR_+}$ \ be a multi-type CBI process with parameters
 \ $(d, \bc, \Bbeta, \bB, \nu, \bmu)$ \ such that \ $\EE(\|\bX_0\|) < \infty$
 \ and the moment condition
 \eqref{moment_condition_m_new} holds.
Suppose that \ $(\bX_t)_{t\in\RR_+}$ \ is irreducible.
Then \ $(\bX_t)_{t\in\RR_+}$ \ is called
 \[
   \begin{cases}
    subcritical & \text{if \ $s(\tbB) < 0$,} \\
    critical & \text{if \ $s(\tbB) = 0$,} \\
    supercritical & \text{if \ $s(\tbB) > 0$.}
   \end{cases}
 \]
\end{Def}
For motivations of Definitions \ref{Def_irreducible} and
 \ref{Def_indecomposable_crit}, see Barczy and Pap \cite[Section 3]{BarPap}.

Next we will recall a convergence result for irreducible and critical
 multi-type CBI processes.

A function \ $f : \RR_+ \to \RR^d$ \ is called \emph{c\`adl\`ag} if it is right
 continuous with left limits.
\ Let \ $\DD(\RR_+, \RR^d)$ \ and \ $\CC(\RR_+, \RR^d)$ \ denote the space of
 all \ $\RR^d$-valued c\`adl\`ag and continuous functions on \ $\RR_+$,
 \ respectively.
Let \ $\cD_\infty(\RR_+, \RR^d)$ \ denote the Borel $\sigma$-field in
 \ $\DD(\RR_+, \RR^d)$ \ for the metric characterized by Jacod and Shiryaev
 \cite[VI.1.15]{JacShi} (with this metric \ $\DD(\RR_+, \RR^d)$ \ is a complete
 and separable metric space).
For \ $\RR^d$-valued stochastic processes \ $(\bcY_t)_{t\in\RR_+}$ \ and
 \ $(\bcY^{(n)}_t)_{t\in\RR_+}$, \ $n \in \NN$, \ with c\`adl\`ag paths we write
 \ $\bcY^{(n)} \distr \bcY$ \ as \ $n \to \infty$ \ if the distribution of
 \ $\bcY^{(n)}$ \ on the space \ $(\DD(\RR_+, \RR^d), \cD_\infty(\RR_+, \RR^d))$
 \ converges weakly to the distribution of \ $\bcY$ \ on the space
 \ $(\DD(\RR_+, \RR^d), \cD_\infty(\RR_+, \RR^d))$ \ as \ $n \to \infty$.
\ Concerning the notation \ $\distr$ \ we note that if \ $\xi$ \ and \ $\xi_n$,
 \ $n \in \NN$, \ are random elements with values in a metric space
 \ $(E, \rho)$, \ then we also denote by \ $\xi_n \distr \xi$ \ the weak
 convergence of the distributions of \ $\xi_n$ \ on the space \ $(E, \cB(E))$
 \ towards the distribution of \ $\xi$ \ on the space \ $(E, \cB(E))$ \ as
 \ $n \to \infty$, \ where \ $\cB(E)$ \ denotes the Borel $\sigma$-algebra on
 \ $E$ \ induced by the given metric \ $\rho$.

The proof of the following convergence theorem can be found in Barczy and Pap
 \cite[Theorem 4.1 and Lemma A.3]{BarPap}.

\begin{Thm}\label{conv}
Let \ $(\bX_t)_{t\in\RR_+}$ \ be a multi-type CBI process with
 parameters \ $(d, \bc, \Bbeta, \bB, \nu, \bmu)$ \ such that
 \ $\EE(\|\bX_0\|^4) < \infty$ \ and the moment conditions
 \begin{equation}\label{moment_condition_m}
  \int_{\cU_d} \|\bz\|^q \bbone_{\{\|\bz\|\geq1\}} \, \nu(\dd \bz) < \infty ,
  \qquad
  \int_{\cU_d} \|\bz\|^q \bbone_{\{\|\bz\|\geq1\}} \, \mu_i(\dd \bz) < \infty ,
  \quad
  i \in \{1, \ldots, d\}
 \end{equation}
 hold with \ $q = 4$.
\ Suppose that \ $(\bX_t)_{t\in\RR_+}$ \ is irreducible and critical.
Then
 \begin{gather}\label{Conv_X}
  (\bcX_t^{(n)})_{t\in\RR_+} := (n^{-1} \bX_{\nt})_{t\in\RR_+}
  \distr (\bcX_t)_{t\in\RR_+} := (\cZ_t \tbu)_{t\in\RR_+} \qquad
  \text{as \ $n \to \infty$}
 \end{gather}
 in \ $\DD(\RR_+, \RR^d)$, \ where \ $\tbu \in \RR_{++}^d$ \ is the right
 Perron vector of \ $\ee^{\tbB}$ \ (corresponding to the eigenvalue 1 with
 \ $\sum_{i=1}^d \be_i^\top \tbu = 1$),
 \ $(\cZ_t)_{t \in \RR_+}$ \ is the pathwise unique strong solution of the SDE
 \begin{equation}\label{SDE_Y}
  \dd \cZ_t
  = \langle \bu, \tBbeta \rangle \, \dd t
    + \sqrt{ \langle \obC \bu, \bu \rangle \cZ_t^+ }
      \, \dd \cW_t ,
  \qquad t \in \RR_+ , \qquad \cZ_0 = 0 ,
 \end{equation}
 where \ $\bu \in \RR_{++}^d$ \ is the left Perron vector of \ $\ee^{\tbB}$
 \ (corresponding to the eigenvalue 1 with \ $\bu^\top \tbu = 1$),
 \ $(\cW_t)_{t \in \RR_+}$ \ is a standard Brownian motion and
 \begin{equation}\label{obC=}
  \obC := \sum_{k=1}^d \langle \be_k, \tbu \rangle \, \bC_k
       \in \RR_+^{d \times d}
 \end{equation}
 with
 \begin{equation}\label{bCk}
  \bC_k := 2 c_k \be_k \be_k^\top + \int_{\cU_d} \bz \bz^\top \mu_k(\dd \bz)
        \in \RR_+^{d \times d} , \qquad
  k \in \{1, \ldots, d\} .
 \end{equation}
\end{Thm}

The moment conditions \eqref{moment_condition_m} with \ $q = 4$ \ in Theorem
 \ref{conv} are used only for checking the conditional Lindeberg condition,
 namely, condition (ii) of Theorem \ref{Conv2DiffThm}.
For a more detailed discussion, see Barczy and Pap \cite[Remark 4.2]{BarPap}.
Note also that Theorem \ref{conv} is in accordance with Theorem 3.1 in
 Isp\'any and Pap \cite{IspPap2}.

\begin{Rem}\label{REMARK_SDE}
The SDE \eqref{SDE_Y} has a pathwise unique strong solution
 \ $(\cZ_t^{(z)})_{t\in\RR_+}$ \ for all initial values
 \ $\cZ_0^{(z)} = z \in \RR$, \ and if the initial value \ $z$ \ is
 nonnegative, then \ $\cZ_t^{(z)}$ \ is nonnegative for all \ $t \in \RR_+$
 \ with probability 1, since
 \ $\langle \bu, \tBbeta \rangle \in \RR_+$, \ see, e.g., Ikeda and
 Watanabe \cite[Chapter IV, Example 8.2]{IkeWat}.
\proofend
\end{Rem}

\begin{Rem}\label{REMARK_par}
Note that for the definition of \ $\bC_k$, \ $k \in \{1, \ldots, d\}$ \ and
 \ $\obC$, \ the moment conditions \eqref{moment_condition_m} are
 needed only with \ $q = 2$.
\ Moreover, \ $\langle \obC \bu, \bu \rangle = 0$ \ if and only
 if \ $\bc = \bzero$ \ and \ $\bmu = \bzero$, \ when the pathwise unique
 strong solution of \eqref{SDE_Y} is the deterministic function
 \ $\cZ_t = \langle \bu, \tBbeta \rangle \, t$, \ $t \in \RR_+$.
\ Indeed,
 \[
   \langle \obC \bu, \bu \rangle
   = \sum_{k=1}^d
      \langle \be_k, \tbu \rangle
      \biggl( 2 c_k \langle \be_k, \bu \rangle^2
              + \int_{\cU_d}
                 \langle \bz, \bu \rangle^2 \, \mu_k(\dd \bz) \biggr) .
 \]
Further, \ $\obC$ \ in \eqref{SDE_Y} can be replaced by
 \begin{equation}\label{tbC=}
  \tbC := \sum_{i=1}^d \langle \be_i, \tbu \rangle \bV_i
        = \var(\bY_1 \mid \bY_0 = \tbu) ,
 \end{equation}
 where the matrices \ $\bV_i$, \ $i \in \{1, \ldots, d\}$, \ are defined in
 Proposition \ref{moment_formula_2}, \ and \ $(\bY_t)_{t\in\RR_+}$ \ is a multi-type
 CBI process with parameters \ $(d, \bc, \bzero, \bB, 0, \bmu)$ \ such that the
 moment conditions \eqref{moment_condition_m} hold with \ $q = 2$.
\ In fact, \ $(\bY_t)_{t\in\RR_+}$ \ is a multi-type CBI process without
 immigration such that its branching mechanism is the same as that of
 \ $(\bX_t)_{t\in\RR_+}$.
\ Indeed, by the spectral mapping theorem, \ $\bu$ \ and \ $\tbu$ \ are
 left and right eigenvectors of \ $\ee^{s\tbB}$, \ $s \in \RR_+$, \ belonging to the
 eigenvalue \ $1$, \ respectively, hence
 \begin{align*}
  &\langle \tbC \bu, \bu \rangle
   = \sum_{i=1}^d \langle \be_i, \tbu \rangle
      \bu^\top \bV_i \, \bu
   = \sum_{i=1}^d \langle \be_i, \tbu \rangle
      \sum_{\ell=1}^d
       \int_0^1
        \langle \ee^{(1-u)\tbB} \be_i, \be_\ell \rangle
        \bu^\top \ee^{u\tbB} \bC_\ell \, \ee^{u\tbB^\top} \!\! \bu
        \, \dd u \\
  &\qquad
   = \sum_{i=1}^d \langle \be_i, \tbu \rangle
      \sum_{\ell=1}^d
       \int_0^1
        \langle \ee^{(1-u)\tbB} \be_i, \be_\ell \rangle
        \bu^\top \bC_\ell \, \bu \, \dd u
   = \sum_{\ell=1}^d
      \int_0^1
       \be_\ell^\top \ee^{(1-u)\tbB} \sum_{i=1}^d \be_i \be_i^\top \tbu
       \langle \bC_\ell \bu, \bu \rangle \, \dd u \\
  &\qquad
   = \sum_{\ell=1}^d
      \int_0^1
       \be_\ell^\top \ee^{(1-u)\tbB} \tbu
       \langle \bC_\ell \bu, \bu \rangle \, \dd u
   = \sum_{\ell=1}^d
      \be_\ell^\top \tbu \langle \bC_\ell \bu, \bu \rangle
   = \langle \obC \bu, \bu \rangle .
 \end{align*}

Clearly, \ $\obC$ \ and \ $\tbC$ \ depend only on the branching mechanism.
Note that for each \ $i \in \{1, \ldots, d\}$,
 \ $\bV_i = \sum_{j=1}^d (\be_j^\top \be_i) \bV_j
    = \var(\bY_1 \mid \bY_0 = \be_i)$, \ see Proposition
 \ref{moment_formula_2}.
\proofend
\end{Rem}

\section{Main results}
\label{section_CBI_2}

Let \ $(\bX_t)_{t\in\RR_+}$ \ be a 2-type CBI process with parameters
 \ $(2, \bc, \Bbeta, \bB, \nu, \bmu)$ \ such that the moment condition
 \eqref{moment_condition_m_new} holds.
We call the process \ $(\bX_t)_{t\in\RR_+}$ \ doubly symmetric if
 \ $\tb_{1,1} = \tb_{2,2} =: \gamma \in \RR$ \ and
 \ $\tb_{1,2} = \tb_{2,1} =: \kappa \in \RR_+$, \ where
 \ $\tbB = (\tb_{i,j})_{i,j\in\{1,2\}}$ \ is defined in \eqref{tbB} with
 \ $d = 2$, that is, if \ $\tbB$ \ takes the form
 \begin{equation}\label{tbB2}
  \tbB = \begin{bmatrix} \gamma & \kappa \\ \kappa & \gamma \end{bmatrix}
 \end{equation}
 with some \ $\gamma \in \RR$ \ and \ $\kappa \in \RR_+$.
\ For the sake of simplicity, we suppose \ $\bX_0 = \bzero$.
\ In the sequel we also assume that \ $\Bbeta \ne \bzero$ \ or \ $\nu \ne 0$
 \ (i.e., the immigration mechanism is non-zero), equivalently,
 \ $\tBbeta \ne \bzero$ \ (where \ $\tBbeta$ \ is defined in \eqref{tBbeta}),
 otherwise \ $\bX_t = \bzero$ \ for all \ $t \in \RR_+$, \ following from
 \eqref{EXcond}.
Clearly \ $\tbB$ \ is irreducible if and only if \ $\kappa \in \RR_{++}$,
 \ since \ $\bP^\top \tbB \bP = \tbB$ \ for all permutation matrices
 \ $\bP \in \RR^{2\times2}$.
\ Hence \ $(\bX_t)_{t\in\RR_+}$ \ is irreducible if and only if
 \ $\kappa \in \RR_{++}$, \ see Definition \ref{Def_irreducible}.
The eigenvalues of \ $\tbB$ \ are \ $\gamma - \kappa$ \ and
 \ $\gamma + \kappa$, \ thus \ $s:= s(\tbB) = \gamma + \kappa$, \ which is
 called the \emph{criticality parameter}, and \ $(\bX_t)_{t\in\RR_+}$ \ is critical
 if and only if \ $s = 0$, \ see Definition \ref{Def_indecomposable_crit}.

For \ $k \in \ZZ_+$, \ let \ $\cF_k := \sigma(\bX_0, \bX_1 , \dots, \bX_k)$.
\ Since \ $(\bX_k)_{k\in\ZZ_+}$ \ is a time-homogeneous Markov process, by
 \eqref{EXcond},
 \begin{equation}\label{mart}
  \EE(\bX_k \mid \cF_{k-1}) = \EE(\bX_k \mid \bX_{k-1})
  = \ee^\tbB \bX_{k-1} + \oBbeta ,
  \qquad k \in \NN ,
 \end{equation}
 where
 \begin{equation}\label{ttBbeta}
  \oBbeta := \int_0^1 \ee^{s\tbB} \tBbeta \, \dd s
           \in \RR_+^2 .
  \end{equation}
Note that \ $\oBbeta = \EE(\bX_1 \mid \bX_0 = \bzero)$, \ see \eqref{EXcond}.
Note also that \ $\oBbeta$ \ depends both on the branching and immigration
 mechanisms, although \ $\tBbeta$ \ depends only on the immigration mechanism.
Let us introduce the sequence
 \begin{equation}\label{Mk}
  \bM_k
  := \bX_k - \EE(\bX_k \mid \cF_{k-1})
   = \bX_k - \ee^\tbB \bX_{k-1} - \oBbeta ,
  \qquad k \in \NN ,
 \end{equation}
 of martingale differences with respect to the filtration
 \ $(\cF_k)_{k \in \ZZ_+}$.
By \eqref{Mk}, the process \ $(\bX_k)_{k \in \ZZ_+}$ \ satisfies the recursion
 \begin{equation}\label{regr}
  \bX_k = \ee^\tbB \bX_{k-1} + \oBbeta + \bM_k ,
  \qquad k \in \NN .
 \end{equation}
By the so-called Putzer's spectral formula, see, e.g., Putzer \cite{Put}, we
 have
 \begin{equation}\label{Putzer}
   \ee^{t\tbB}
   = \frac{\ee^{(\gamma+\kappa)t}}{2}
       \begin{bmatrix}
       1 & 1 \\
       1 & 1
      \end{bmatrix}
     + \frac{\ee^{(\gamma-\kappa)t}}{2}
       \begin{bmatrix}
        1 & -1 \\
        -1 & 1
       \end{bmatrix}
   = \ee^{\gamma t}
     \begin{bmatrix}
      \cosh(\kappa t) & \sinh(\kappa t) \\
      \sinh(\kappa t) & \cosh(\kappa t)
     \end{bmatrix} ,
   \qquad t \in \RR_+ .
 \end{equation}
Consequently,
 \[
   \ee^\tbB = \begin{bmatrix} \alpha & \beta \\ \beta & \alpha \end{bmatrix}
   \qquad \text{with} \qquad
   \alpha := \ee^\gamma \cosh(\kappa) , \qquad
   \beta := \ee^\gamma \sinh(\kappa) .
 \]
Considering the eigenvalues \ $\varrho := \alpha + \beta$ \ and
 \ $\delta := \alpha - \beta$ \ of \ $\ee^\tbB$, \ we have
 \ $\alpha = (\varrho + \delta)/2$ \ and \ $\beta = (\varrho - \delta)/2$,
 \ thus the recursion \eqref{regr} can be written in the form
 \begin{equation}\label{Xregr}
   \bX_k = \frac{1}{2}
           \begin{bmatrix}
            \varrho + \delta & \varrho - \delta \\
            \varrho - \delta & \varrho + \delta
           \end{bmatrix}
           \bX_{k-1} + \bM_k + \oBbeta  ,
   \qquad k \in \NN .
 \end{equation}
For each \ $n \in \NN$, \ a CLS estimator
 \[
   (\hvarrho_n, \hdelta_n, \hoBbeta_n)
   := (\hvarrho_n(\bX_1, \ldots, \bX_n), \hdelta_n(\bX_1, \ldots, \bX_n),
       \hoBbeta_n(\bX_1, \ldots, \bX_n))
 \]
 of \ $(\varrho, \delta, \oBbeta)$ \ based on a sample \ $\bX_1, \ldots, \bX_n$
 \ can be obtained by minimizing the sum of squares
 \[
   \sum_{k=1}^n
    \left\| \bX_k
            - \frac{1}{2}
              \begin{bmatrix}
               \varrho + \delta & \varrho - \delta \\
               \varrho - \delta & \varrho + \delta
              \end{bmatrix}
              \bX_{k-1}
            - \oBbeta \right\|^2
 \]
 with respect to \ $(\varrho, \delta, \oBbeta)$ \ over \ $\RR^4$,
 \ and it has the form
 \begin{gather}
  \hvarrho_n
  = \frac{n \sum_{k=1}^n
            \langle \bu, \bX_k \rangle
            \langle \bu, \bX_{k-1} \rangle
           - \sum_{k=1}^n \langle \bu, \bX_k \rangle
             \sum_{k=1}^n \langle \bu, \bX_{k-1} \rangle}
          {n \sum_{k=1}^n \langle \bu, \bX_{k-1} \rangle^2
           - \bigl(\sum_{k=1}^n \langle \bu, \bX_{k-1} \rangle\bigr)^2} ,
  \label{CLSEr} \\
  \hdelta_n
  = \frac{n \sum_{k=1}^n
            \langle \bv, \bX_k \rangle
            \langle \bv, \bX_{k-1} \rangle
           - \sum_{k=1}^n \langle \bv, \bX_k \rangle
             \sum_{k=1}^n \langle \bv, \bX_{k-1} \rangle}
          {n \sum_{k=1}^n \langle \bv, \bX_{k-1} \rangle^2
           - \bigl(\sum_{k=1}^n \langle \bv, \bX_{k-1} \rangle\bigr)^2} ,
  \label{CLSEd} \\
  \hoBbeta_n
  = \frac{1}{n} \sum_{k=1}^n \bX_k
     - \frac{1}{2n}
       \sum_{k=1}^n
        \begin{bmatrix}
         \langle \bu, \bX_{k-1} \rangle & \langle \bv, \bX_{k-1} \rangle \\
         \langle \bu, \bX_{k-1} \rangle & - \langle \bv, \bX_{k-1} \rangle
        \end{bmatrix}
        \begin{bmatrix}
         \hvarrho_n \\
         \hdelta_n
        \end{bmatrix}
  \label{CLSEb}
 \end{gather}
 on the set \ $H_n \cap \tH_n$, \ where
 \begin{gather*}
  \bu = \begin{bmatrix} 1 \\ 1 \end{bmatrix} \in \RR_{++}^2 , \qquad
  \bv := \begin{bmatrix} 1 \\ -1 \end{bmatrix} \in \RR^2 , \\
  H_n:=\biggl\{\omega \in \Omega
               : n \sum_{k=1}^n
                  \langle \bu, \bX_{k-1}(\omega) \rangle^2
                 - \biggl(\sum_{k=1}^n \langle \bu, \bX_{k-1}(\omega) \rangle\biggr)^2
                 > 0\biggr\} , \\
  \tH_n:=\biggl\{\omega \in \Omega
                 : n \sum_{k=1}^n
                    \langle \bv, \bX_{k-1}(\omega) \rangle^2
                   - \biggl(\sum_{k=1}^n \langle \bv, \bX_{k-1}(\omega) \rangle\biggr)^2
                   > 0\biggr\} ,
 \end{gather*}
 see Lemma \ref{CLSE}.
Here \ $\bu$ \ and \ $\bv$ \ are left eigenvectors of \ $\tbB$
 \ belonging to the eigenvalues \ $\gamma + \kappa$ \ and \ $\gamma - \kappa$,
 \ respectively, hence they are left eigenvectors of \ $\ee^\tbB$ \ belonging
 to the eigenvalues \ $\varrho = \ee^{\gamma + \kappa}$ \ and
 \ $\delta = \ee^{\gamma - \kappa}$, \ respectively.
In a natural way, one can extend the CLS estimators \ $\hvarrho_n$ \ and
 \ $\hdelta_n$ \ to the set \ $H_n$ \ and \ $\tH_n$, \ respectively.

In the sequel we investigate the critical case.
By Lemma \ref{LEMMA_CLSE_exist_discrete}, \ $\PP(H_n) \to 1$ \ and
 \ $\PP(\tH_n) \to 1$ \ as \ $n \to \infty$ \ under appropriate assumptions.
Let us introduce the function \ $h : \RR^4 \to \RR_{++}^2 \times \RR^2$ \ by
 \[
   h(\gamma, \kappa, \tBbeta)
   := \biggl( \ee^{\gamma+\kappa}, \ee^{\gamma-\kappa},
              \Bigl(\int_0^1 \ee^{s\tbB} \, \dd s\Bigr) \tBbeta \biggr)
   = (\varrho, \delta, \oBbeta) , \qquad (\gamma, \kappa, \tBbeta) \in \RR^4 ,
 \]
 where, by formula \eqref{Putzer},
 \[
   \int_0^1 \ee^{s\tbB} \, \dd s
   = \frac{1}{2}
     \begin{bmatrix}
      \int_0^1 \ee^{(\gamma+\kappa)s} \, \dd s
      + \int_0^1 \ee^{(\gamma-\kappa)s} \, \dd s
       & \int_0^1 \ee^{(\gamma+\kappa)s} \, \dd s
         - \int_0^1 \ee^{(\gamma-\kappa)s} \, \dd s \\
      \int_0^1 \ee^{(\gamma+\kappa)s} \, \dd s
      - \int_0^1 \ee^{(\gamma-\kappa)s} \, \dd s
       & \int_0^1 \ee^{(\gamma+\kappa)s} \, \dd s
         + \int_0^1 \ee^{(\gamma-\kappa)s} \, \dd s
     \end{bmatrix} .
 \]
Note that \ $h$ \ is bijective having inverse
 \[
   h^{-1}(\varrho, \delta, \oBbeta)
   = \left( \frac{1}{2} \log(\varrho \delta),
            \frac{1}{2} \log\left(\frac{\varrho}{\delta}\right),
            \Bigl(\int_0^1 \ee^{s\tbB} \, \dd s\Bigr)^{-1} \oBbeta \right)
   = (\gamma, \kappa, \tBbeta) , \qquad
   (\varrho, \delta, \oBbeta) \in \RR_{++}^2 \times \RR^2 ,
 \]
 with
 \begin{equation}\label{int_ee^stbB}
   \biggl(\int_0^1 \ee^{s\tbB} \, \dd s\biggr)^{-1}
   = \frac{1}{2\int_0^1 \varrho^s \, \dd s}
     \begin{bmatrix} 1 & 1 \\ 1 & 1 \end{bmatrix}
     + \frac{1}{2\int_0^1 \delta^s \, \dd s}
       \begin{bmatrix} 1 & -1 \\ -1 & 1 \end{bmatrix} .
 \end{equation}
Theorem \ref{main_rdb} will imply that, under appropriate assumptions, the
 CLS estimator \ $(\hvarrho_n, \hdelta_n)$ \ of \ $(\varrho, \delta)$ \ is
 weakly consistent, hence, \ $\PP((\hvarrho_n, \hdelta_n) \in \RR_{++}^2) \to 1$
 \ as \ $n \to \infty$, \ and
 \[
   (\hvarrho_n, \hdelta_n, \hoBbeta_n)
   = \argmin_{(\varrho,\delta,\oBbeta)\in\RR_{++}^2\times\RR^2}
      \sum_{k=1}^n
       \left\| \bX_k
            - \frac{1}{2}
              \begin{bmatrix}
               \varrho + \delta & \varrho - \delta \\
               \varrho - \delta & \varrho + \delta
              \end{bmatrix}
              \bX_{k-1}
            - \oBbeta \right\|^2
 \]
 on the set
 \ $\{\omega \in \Omega : (\hvarrho_n(\omega), \hdelta_n(\omega)) \in \RR_{++}^2\}$.
\ Thus one can introduce a natural estimator of \ $(\gamma, \kappa, \tBbeta)$
 \ by applying the inverse of \ $h$ \ to the CLS estimator of
 \ $(\varrho, \delta, \oBbeta)$, \ that is,
 \[
   (\hgamma_n, \hkappa_n, \htBbeta_n)
   := h^{-1}(\hvarrho_n, \hdelta_n, \hoBbeta_n)
   = \left( \frac{1}{2} \log(\hvarrho_n \hdelta_n),
            \frac{1}{2} \log\left(\frac{\hvarrho_n}{\hdelta_n}\right),
            \Bigl(\int_0^1 \ee^{s\htbB_n} \, \dd s\Bigr)^{-1} \hoBbeta_n \right) ,
   \qquad \ n \in \NN ,
 \]
 on the set
 \ $\{\omega \in \Omega
      : (\hvarrho_n(\omega), \hdelta_n(\omega)) \in \RR_{++}^2\}$,
 \ where
 \[
   \htbB_n
   := \begin{bmatrix} \hgamma_n & \hkappa_n \\ \hkappa_n & \hgamma_n \end{bmatrix} ,
   \qquad n \in \NN .
 \]
We also obtain
 \[
   \bigl(\hgamma_n, \hkappa_n, \htBbeta_n\bigr)
   = \argmin_{(\gamma,\kappa,\tBbeta)\in\RR^4}
      \sum_{k=1}^n
       \left\| \bX_k
               - \ee^\gamma
                 \begin{bmatrix}
                  \cosh(\kappa) & \sinh(\kappa) \\
                  \sinh(\kappa) & \cosh(\kappa)
                 \end{bmatrix}
                 \bX_{k-1}
                - \Bigl(\int_0^1 \ee^{s\tB} \, \dd s\Bigr) \tBbeta \right\|^2
 \]
 on the set
 \ $\{\omega \in \Omega : (\hvarrho_n(\omega), \hdelta_n(\omega)) \in \RR_{++}^2\}$,
 \ hence the probability that \ $\bigl(\hgamma_n, \hkappa_n, \htBbeta_n\bigr)$
 \ is the CLS estimator of \ $(\gamma, \kappa, \tBbeta)$ \ converges to 1 as
 \ $n \to \infty$.
\ In a similar way, the probability that
 \[
   \hs_n := \log \hvarrho_n
 \]
 is the CLS estimator of the criticality parameter \ $s = \gamma + \kappa$
 \ converges to 1 as \ $n \to \infty$.

\begin{Thm}\label{main}
Let \ $(\bX_t)_{t\in\RR_+}$ \ be a 2-type CBI process with parameters
 \ $(2, \bc, \Bbeta, \bB, \nu, \bmu)$ \ such that \ $\bX_0 = \bzero$, \ the
 moment conditions \eqref{moment_condition_m} hold with \ $q = 8$,
 \ $\Bbeta \ne \bzero$ \ or \ $\nu \ne 0$, \ and \eqref{tbB2} holds with some
 \ $\gamma \in \RR$ \ and \ $\kappa \in \RR_{++}$ \ such that
 \ $s = \gamma + \kappa = 0$ \ (hence the process is irreducible and critical).
Then the probability of the existence of the estimator \ $\hs_n$ \ converges
 to 1 as \ $n \to \infty$ \ and
 \begin{equation}\label{h}
  n (\hs_n - s)
  \distr
  \frac{\int_0^1 \cY_t \, \dd (\cM_{t,1} + \cM_{t,2})
        - (\cM_{1,1} + \cM_{1,2}) \int_0^1 \cY_t \, \dd t}
       {\int_0^1 \cY_t^2 \, \dd t - \bigl(\int_0^1 \cY_t \, \dd t\bigr)^2}
  =: \cI
  \qquad \text{as \ $n \to \infty$,}
 \end{equation}
 where \ $(\bcM_t)_{t \in \RR_+} = (\cM_{t,1}, \cM_{t,2})_{t \in \RR_+}$ \ is
 the pathwise unique strong solution of the SDE
 \begin{equation}\label{SDE_M}
   \dd \bcM_t \\
   = ((\cM_{t,1} + \cM_{t,2} + (\tbeta_1 + \tbeta_2) t)^+)^{1/2} \, \tbC^{1/2} \,
     \dd \bcW_t , \qquad t \in \RR_+ , \qquad \bcM_0 = \bzero ,
 \end{equation}
 where \ $\tbC^{1/2}$ \ denotes the unique symmetric and positive semidefinite
 square root of \ $\tbC$, \ $(\bcW_t)_{t \in \RR_+}$ \ is a 2-dimensional
 standard Wiener process, and
 \[
   \cY_t := \cM_{t,1} + \cM_{t,2} + (\tbeta_1 + \tbeta_2) t , \qquad t \in \RR_+ .
 \]

If \ $\bc = \bzero$ \ and \ $\bmu = \bzero$, \ then
 \begin{equation}\label{h1}
  n^{3/2} (\hs_n - s)
  \distr
  \cN\biggl(0, \frac{12}{(\tbeta_1 + \tbeta_2)^2}
              \int_{\cU_2} (z_1 + z_2)^2 \, \nu(\dd\bz)\biggr)
  \qquad \text{as \ $n \to \infty$.}
 \end{equation}

If
 \ $\|\bc\|^2 + \sum_{i=1}^2 \int_{\cU_2} (z_1 - z_2)^2 \, \mu_i(\dd\bz) > 0$,
 \ then the probability of the existence of the estimator
 \ $(\hgamma_n, \hkappa_n, \htBbeta_n)$ \ converges to 1 as \ $n \to \infty$,
 \ and
 \begin{equation}\label{gamma,kappa,tBbeta}
  \begin{bmatrix}
   n^{1/2} (\hgamma_n - \gamma) \\
   n^{1/2} (\hkappa_n - \kappa) \\
   \htBbeta_n - \tBbeta
  \end{bmatrix}
  \distr
  \begin{bmatrix}
   \frac{1}{2} \sqrt{\ee^{2(\kappa-\gamma)} - 1} \,
   \frac{\int_0^1 \cY_t \, \dd\tcW_t}
        {\int_0^1 \cY_t \, \dd t}
   \begin{bmatrix} 1 \\ -1 \end{bmatrix} \\
   \frac{1}{2}
   \Biggl( \begin{bmatrix} 1 & 1 \\ 1 & 1 \end{bmatrix}
           + \frac{\kappa - \gamma}{1 - \ee^{\gamma-\kappa}}
             \begin{bmatrix} 1 & -1 \\ -1 & 1 \end{bmatrix} \Biggr)
   \bcM_1
   - \frac{1}{2} \cI \int_0^1 \cY_t \, \dd t
     \begin{bmatrix} 1 \\ 1 \end{bmatrix}
  \end{bmatrix}
 \end{equation}
 as \ $n \to \infty$, \ where \ $(\tcW_t)_{t \in \RR_+}$ \ is a standard Wiener
 process, independent from \ $(\bcW_t)_{t \in \RR_+}$.

If
 \ $\|\bc\|^2 + \sum_{i=1}^2 \int_{\cU_2} (z_1 - z_2)^2 \, \mu_i(\dd\bz) = 0$
 \ and \ $\int_{\cU_2} (z_1 - z_2)^2 \, \nu(\dd\bz) > 0$, \ then the probability
 of the existence of the estimator \ $(\hgamma_n, \hkappa_n, \htBbeta_n)$
 \ converges to 1 as \ $n \to \infty$, \ and
 \begin{equation}\label{gamma,kappa,tBbetat}
  \begin{bmatrix}
   n^{1/2} (\hgamma_n - \gamma) \\
   n^{1/2} (\hkappa_n - \kappa) \\
   \htBbeta_n - \tBbeta
  \end{bmatrix}
  \distr
  \begin{bmatrix}
   \frac{1}{2} \sqrt{\ee^{2(\kappa-\gamma)} - 1} \,
   \tcW_1
   \begin{bmatrix} 1 \\ -1 \end{bmatrix} \\
   \frac{1}{2}
   \bigl(\cM_{1,1} + \cM_{1,2} - \cI \int_0^1 \cY_t \, \dd t\bigr)
   \begin{bmatrix} 1 \\ 1 \end{bmatrix}
  \end{bmatrix}
 \end{equation}
 as \ $n \to \infty$, \ where \ $\tcW_1$ \ is a random variable with standard
 normal distribution, independent from \ $(\bcW_t)_{t \in \RR_+}$.

Furthermore, if \ $\bc = \bzero$, \ $\bmu = \bzero$ \ and
 \ $\int_{\cU_2} (z_1 - z_2)^2 \, \nu(\dd\bz) > 0$, \ then
 \begin{equation}\label{gamma,kappa,tBbeta1}
  \begin{bmatrix}
   n^{1/2} (\hgamma_n - \gamma) \\
   n^{1/2} (\hkappa_n - \kappa) \\
   n^{1/2} (\htBbeta_n - \tBbeta)
  \end{bmatrix}
  \distr
  \cN_4\left(\bzero, \begin{bmatrix}
                      \bR_{1,1} & \bR_{1,2} \\
                      \bR_{2,1} & \bR_{2,2}
                     \end{bmatrix}\right)
  \qquad \text{as \ $n \to \infty$,}
 \end{equation}
 with
 \begin{gather*}
  \bR_{1,1} := \frac{\ee^{2(\kappa-\gamma)}-1}{4}
         \begin{bmatrix} 1 & -1 \\ -1 & 1 \end{bmatrix} , \qquad
  \bR_{2,1} := \bR_{1,2}^\top
  := - \frac{(\tbeta_1-\tbeta_2)(\ee^{2(\kappa-\gamma)}-1)}{4(\kappa-\gamma)}
       \begin{bmatrix} 1 & -1 \\ -1 & 1 \end{bmatrix} , \\
  \begin{aligned}
   \bR_{2,2}
   &:= \int_{\cU_2} (z_1 + z_2)^2 \, \nu(\dd\bz)
       \begin{bmatrix} 1 & 1 \\ 1 & 1 \end{bmatrix}
       + \frac{1}{2} \int_{\cU_2} (z_1^2 - z_2^2) \, \nu(\dd\bz)
         \begin{bmatrix} 1 & 0 \\ 0 & -1 \end{bmatrix} \\
   &\:\quad
       + \frac{(1 - \ee^{2(\gamma-\kappa)})}{4}
         \left\{\frac{\kappa-\gamma}{2(1-\ee^{\gamma-\kappa})^2}
                \int_{\cU_2} (z_1 - z_2)^2 \, \nu(\dd\bz)
                + \frac{(\tbeta_1-\tbeta_2)^2}
                       {(\kappa-\gamma)^2\ee^{2(\gamma-\kappa)}}\right\}
         \begin{bmatrix} 1 & -1 \\ -1 & 1 \end{bmatrix} .
  \end{aligned}
 \end{gather*}
\end{Thm}

Under the assumptions of Theorem \ref{main}, we have the following remarks.

\begin{Rem}\label{REMARK0.5}
If
 \ $\|\bc\|^2 + \sum_{i=1}^2 \int_{\cU_2} (z_1 - z_2)^2 \, \mu_i(\dd\bz)
    + \int_{\cU_2} (z_1 - z_2)^2 \, \nu(\dd\bz) + (\tbeta_1 - \tbeta_2)^2 = 0$,
 \ then, by Lemma \ref{main_VVt} and \eqref{tbCbv_LEFT,bv_LEFT=0},
 \ $X_{k,1}\ase X_{k,2}$ \ for all \ $k \in \NN$, \ hence \ $\hdelta_n$ \ and
 \ $\hoBbeta_n$, \ $n \in \NN$, \ are not defined (see Lemma \ref{CLSE}).
Note that in Theorem \ref{main} we have covered all the possible cases when in
 Lemma \ref{CLSE} we showed that the probability of the existence of the
 estimators converges to 1 as the sample size converges to infinity.
\proofend
\end{Rem}

\begin{Rem}\label{REMARK0.7}
If
 \ $\|\bc\|^2 + \sum_{i=1}^2 \int_{\cU_2} (z_1 - z_2)^2 \, \mu_i(\dd\bz) = 0$
 \ and \ $\int_{\cU_2} (z_1 - z_2)^2 \, \nu(\dd\bz) > 0$, \ then the last two
 coordinates of the limit in \eqref{gamma,kappa,tBbeta} and in
 \eqref{gamma,kappa,tBbetat} are the same.
Indeed, by \eqref{tbCbv_LEFT,bv_LEFT=0},
 \ $\|\bc\|^2 + \sum_{i=1}^2 \int_{\cU_2} (z_1 - z_2)^2 \, \mu_i(\dd\bz) = 0$
 \ is equivalent to \ $\langle \tbC \bv, \bv \rangle = 0$.
Moreover,
 \ $\langle \tbC \bv, \bv \rangle
    = \|\bv^\top \tbC^{1/2}\|^2 = 0$
 \ implies \ $\bv^\top \tbC^{1/2} = \bzero^\top$, \ hence, by It\^o's
 formula, \ $\dd(\bv^\top  \bcM_t) = 0$, \ $t \in \RR_+$, \ implying that
 \ $\cM_{t,1} - \cM_{t,2} = \bv^\top  \bcM_t = 0$, \ $t \in \RR_+$.
\proofend
\end{Rem}

\begin{Rem}\label{REMARK1}
By It\^o's formula, \eqref{SDE_M} yields that \ $(\cY_t)_{t \in \RR_+}$ \ satisfies
 the SDE \eqref{SDE_Y} with initial value \ $\cY_0 = 0$.
\ Indeed, by It\^o's formula and the SDE \eqref{SDE_M} we obtain
 \[
   \dd \cY_t
   = \langle \bu, \tBbeta \rangle \, \dd t
     + (\cY_t^+)^{1/2} \bu^\top \tbC^{1/2} \, \dd \bcW_t , \qquad
   t \in \RR_+ .
 \]
If
 \ $\langle \tbC \bu, \bu \rangle
    = \|\bu^\top \tbC^{1/2}\|^2 = 0$
 \ then \ $\bu^\top \tbC^{1/2} = \bzero$, \ hence
 \ $\dd \cY_t = \langle \bu, \tBbeta \rangle \dd t$, \ $t \in \RR_+$,
 \ implying that the process \ $(\cY_t)_{t\in\RR_+}$ \ satisfies the SDE
 \eqref{SDE_Y}.
If \ $\langle \tbC \bu, \bu \rangle \ne 0$ \ then the process
 \[
   \tcW_t := \frac{\langle \tbC^{1/2} \bu, \bcW_t \rangle}
                  {\langle \tbC \bu, \bu \rangle^{1/2}} ,
   \qquad t \in \RR_+ ,
 \]
 is a (one-dimensional) standard Wiener process.
Consequently, the process \ $(\cY_t)_{t\in\RR_+}$ \ satisfies the SDE
 \eqref{SDE_Y}, since \ $\obC$ \ can be replaced by \ $\tbC$ \ in
 \eqref{SDE_Y}.

Consequently, \ $(\cY_t)_{t\in\RR_+} \distre (\cZ_t)_{t\in\RR_+}$, \ where
 \ $(\cZ_t)_{t \in \RR_+}$ \ is the unique strong solution of the SDE
 \eqref{SDE_Y} with initial value \ $\cZ_0 = 0$, \ hence, by Theorem
 \ref{conv},
 \begin{equation}\label{convX}
  (\bcX^{(n)}_t)_{t\in\RR_+}
  \distr (\bcX_t)_{t\in\RR_+}
  \distre (\cY_t \tbu)_{t\in\RR_+}
  \qquad \text{as \ $n \to \infty$.}
 \end{equation}

If \ $\langle \tbC \, \bu, \bu \rangle = 0$, \ which is
 equivalent to \ $\bc = \bzero$ \ and \ $\bmu = \bzero$
 \ (see Remark \ref{REMARK_par}), then the unique strong solution of
 \eqref{SDE_Y} is the deterministic function
 \ $\cZ_t = \langle \bu, \tBbeta \rangle \, t$, \ $t \in \RR_+$,
 \ hence \ $\cY_t = \langle \bu, \tBbeta \rangle \, t$,
 \ $t \in \RR_+$, \ and \ $\cM_{t,1} + \cM_{t,2} = 0$, \ $t \in \RR_+$.
\ Thus, by \eqref{h}, \ $n \hs_n \distr 0$, \ i.e., the scaling \ $n$ \ is
 not suitable.

If \ $\bc = \bzero$, \ $\bmu = \bzero$ \ and
 \ $\int_{\cU_2} (z_1 - z_2)^2 \, \nu(\dd\bz) > 0$, \ then, by
 \eqref{gamma,kappa,tBbetat},
 \[
   \begin{bmatrix}
    n^{1/2} (\hgamma_n - \gamma) \\
    n^{1/2} (\hkappa_n - \kappa)
   \end{bmatrix}
   \distr
   \cN_2(\bzero, \bR_{1,1})
  \qquad \text{as \ $n \to \infty$,}
 \]
 and \ $\htBbeta_n - \tBbeta \stoch \bzero$ \ as \ $n \to \infty$,
 \ hence we obtain the convergence of the first two coordinates in
 \eqref{gamma,kappa,tBbeta1}, but a suitable scaling of
 \ $\htBbeta_n - \tBbeta$ \ is needed.
\proofend
\end{Rem}

\begin{Rem}\label{REMARK1.5}
In \eqref{h1}, the limit distribution depends on the unknown parameter
 \ $\tBbeta$, \ but one can get rid of this dependence by random scaling
 in the following way.
If \ $\int_{\cU_2} (z_1 - z_2)^2 \, \nu(\dd\bz) > 0$, \ then
 \ $\nu \ne 0$, \ hence
 \ $\int_{\cU_2} (z_1 + z_2)^2 \, \nu(\dd\bz) \ne 0$.
\ If, in addition, \ $\bc = \bzero$ \ and \ $\bmu = \bzero$, \ then, by
 \eqref{gamma,kappa,tBbeta1}, \ $\htBbeta_n \stoch \tBbeta$ \ as
 \ $n \to \infty$, \ hence by \eqref{h1},
 \[
   n^{3/2} \frac{\langle \bu, \htBbeta_n \rangle}
                {\sqrt{12\int_{\cU_2} (z_1 + z_2)^2 \, \nu(\dd\bz)}}
           (\hs_n - s)
   \distr
   \cN(0, 1)
   \qquad \text{as \ $n \to \infty$.}
 \]
A similar random scaling may be applied in case of \eqref{gamma,kappa,tBbeta1},
 however, the variance matrix of the limiting normal distribution is singular
 (since the sum of the first two columns is $\bzero$),
 hence one may use Moore-Penrose pseudoinverse of the variance matrix.
Unfortunately, we can not see how one could get rid of the dependence on the
 unknown parameters in case of convergences \eqref{h},
 \eqref{gamma,kappa,tBbeta} and \eqref{gamma,kappa,tBbetat}, since we have not
 found a way how to eliminate the dependence of the process \ $\cY$ \ on the
 unknown parameters.
In order to perform hypothesis testing, one should investigate the subcritical
 and supercritical cases as well.
\proofend
\end{Rem}

Theorem \ref{main} will follow from the following statement.

\begin{Thm}\label{main_rdb}
Under the assumptions of Theorem \ref{main}, the probability of the existence
 of a unique CLS estimator \ $\hvarrho_n$ \ converges to 1 as \ $n \to \infty$
 \ and
 \begin{equation}\label{rho}
  n (\hvarrho_n - \varrho) \distr \cI \qquad \text{as \ $n \to \infty$.}
 \end{equation}

If \ $\bc = \bzero$ \ and \ $\bmu = \bzero$, \ then
 \begin{equation}\label{rho1}
  n^{3/2} (\hvarrho_n - \varrho)
  \distr
  \cN\left(0, \frac{12}{(\tbeta_1 + \tbeta_2)^2}
              \int_{\cU_2} (z_1 + z_2)^2 \, \nu(\dd\bz)\right)
  \qquad \text{as \ $n \to \infty$.}
 \end{equation}

If
 \ $\|\bc\|^2 + \sum_{i=1}^2 \int_{\cU_2} (z_1 - z_2)^2 \, \mu_i(\dd\bz) > 0$,
 \ then the probability of the existence of a unique CLS estimator
 \ $(\hvarrho_n, \hdelta_n, \hoBbeta_n)$ \ converges to 1 as \ $n \to \infty$,
 \ and
 \begin{equation}\label{rhodeltab}
  \begin{bmatrix}
   n (\hvarrho_n - \varrho) \\
   n^{1/2} (\hdelta_n - \delta) \\
   \hoBbeta_n - \oBbeta
  \end{bmatrix}
  \distr
  \begin{bmatrix}
   \cI \\
   \sqrt{1 - \delta^2} \,
   \frac{\int_0^1 \cY_t \, \dd\tcW_t}{\int_0^1 \cY_t \, \dd t} \\[-1mm]
   \bcM_1
   - \frac{1}{2} \cI \int_0^1 \cY_t \, \dd t
     \begin{bmatrix} 1 \\ 1 \end{bmatrix}
  \end{bmatrix}
 \end{equation}
 as \ $n \to \infty$, \ where \ $(\tcW_t)_{t \in \RR_+}$ \ is a standard Wiener
 process, independent from \ $(\bcW_t)_{t \in \RR_+}$.

If
 \ $\|\bc\|^2 + \sum_{i=1}^2 \int_{\cU_2} (z_1 - z_2)^2 \, \mu_i(\dd\bz) = 0$
 \ and \ $\int_{\cU_2} (z_1 - z_2)^2 \, \nu(\dd\bz) > 0$, \ then the
 probability of the existence of a unique CLS estimator
 \ $(\hvarrho_n, \hdelta_n, \hoBbeta_n)$ \ converges to 1 as \ $n \to \infty$,
 \ and
 \begin{equation}\label{rhodeltabt}
  \begin{bmatrix}
   n (\hvarrho_n - \varrho) \\
   n^{1/2} (\hdelta_n - \delta) \\
   \hoBbeta_n - \oBbeta
  \end{bmatrix}
  \distr
  \begin{bmatrix}
   \cI \\
   \sqrt{1 - \delta^2} \, \tcW_1 \\
   \frac{1}{2} \bigl(\cM_{1,1} + \cM_{1,2} - \cI \int_0^1 \cY_t \, \dd t\bigr)
   \begin{bmatrix} 1 \\ 1 \end{bmatrix}
  \end{bmatrix}
 \end{equation}
 as \ $n \to \infty$.

If \ $\bc = \bzero$, \ $\bmu = \bzero$ \ and
 \ $\int_{\cU_2} (z_1 - z_2)^2 \, \nu(\dd\bz) > 0$, \ then
 \begin{equation}\label{rhodeltab1}
  \begin{bmatrix}
   n^{3/2} (\hvarrho_n - \varrho) \\
   n^{1/2} (\hdelta_n - \delta) \\
   n^{1/2} (\hoBbeta_n - \oBbeta)
  \end{bmatrix}
  \distr
  \cN_4\left(\bzero, \bS\right)
  \qquad \text{as \ $n \to \infty$,}
 \end{equation}
 with
 \begin{align*}
  \bS
  &:=\begin{bmatrix}
      \bzero & \bzero \\
      \bzero & \int_0^1 \int_{\cU_2}
                (\ee^{t\tbB} \bz) (\ee^{t\tbB} \bz)^\top \nu(\dd\bz) \, \dd t
     \end{bmatrix}
     + \frac{3}{4}
       \int_{\cU_2} (z_1 + z_2)^2 \, \nu(\dd\bz)
       \begin{bmatrix}
        \frac{4}{\tbeta_1+\tbeta_2} \be_1 \\
        -\bu
       \end{bmatrix}
       \begin{bmatrix}
        \frac{4}{\tbeta_1+\tbeta_2} \be_1 \\
        -\bu
       \end{bmatrix}^\top \\
  &\:\quad
     + (1 - \delta^2)
       \begin{bmatrix}
        \be_2 \\
        -\frac{\tbeta_1-\tbeta_2}{2\log(\delta^{-1})} \bv
       \end{bmatrix}
       \begin{bmatrix}
        \be_2 \\
        -\frac{\tbeta_1-\tbeta_2}{2\log(\delta^{-1})} \bv
       \end{bmatrix}^\top ,
 \end{align*}
 where \ $\be_1 = \begin{bmatrix} 1 \\ 0 \end{bmatrix}$ \ and
 \ $\be_2 = \begin{bmatrix} 0 \\ 1 \end{bmatrix}$.
\end{Thm}

\noindent
\textbf{Proof of Theorem \ref{main}.}
Before Theorem \ref{main} we have already investigated the existence of the
 estimators \ $\hs_n$ \ and \ $(\hgamma_n, \hkappa_n, \htBbeta_n)$.

In order to prove \eqref{h}, we apply Lemma \ref{Lem_Kallenberg} with
 \ $S = T = \RR$, \ $C = \RR$, \ $\xi = \cI$,
 \ $\xi_n = n(\hvarrho_n - \varrho) = n(\hvarrho_n - 1)$, \ $n \in \NN$,
 \ and with functions \ $f : \RR \to \RR$ \ and \ $f_n : \RR \to \RR$,
 \ $n \in \NN$, \ given by
 \[
   f(x) := x , \quad x \in \RR , \qquad
   f_n(x)
   := \begin{cases}
       n \log\bigl(1 + \frac{x}{n}\bigr) & \text{if \ $x > -n$,} \\
       0 & \text{if \ $x \leq -n$.}
      \end{cases}
 \]
We have \ $f_n(n (\hvarrho_n - 1)) = n \hs_n = n(\hs_n - s)$ \ on the set
 \ $\{\omega \in \Omega : \hvarrho_n(\omega) \in \RR_{++}\}$, \ and
 \ $f_n(x_n) \to f(x)$ \ as \ $n \to \infty$ \ if \ $x_n \to x$ \ as
 \ $n \to \infty$, \ since
 \[
   \lim_{n\to\infty} f_n(x_n)
   = \lim_{n\to\infty} \log\Big[\Bigl( 1 + \frac{x_n}{n} \Bigr)^n\Bigr]
   = \log(\ee^x) = x , \qquad x \in \RR .
 \]
Consequently, \eqref{rho} implies \eqref{h}.

By the method of the proof of \eqref{h} as above, \eqref{rho1} implies \eqref{h1}
 with functions \ $f : \RR \to \RR$ \ and \ $f_n : \RR \to \RR$,
 \ $n \in \NN$, \ given by
 \[
   f(x) := x , \quad x \in \RR , \qquad
   f_n(x)
   := \begin{cases}
       n^{3/2} \log\bigl(1 + \frac{x}{n^{3/2}}\bigr) & \text{if \ $x > -n^{3/2}$,} \\
       0 & \text{if \ $x \leq -n^{3/2}$.}
      \end{cases}
 \]

In order to prove \eqref{gamma,kappa,tBbeta}, first note that
 \begin{equation}\label{gammakappatBbeta}
   \begin{bmatrix}
    n^{1/2} (\hgamma_n - \gamma) \\
    n^{1/2} (\hkappa_n - \kappa) \\
    \htBbeta_n - \tBbeta
   \end{bmatrix}
   = \begin{bmatrix}
      \frac{1}{2}
      \begin{bmatrix} 1 & 1 \\ 1 & -1 \end{bmatrix}
      \begin{bmatrix}
       n^{1/2} (\log(\hvarrho_n) - \log(\varrho)) \\
       n^{1/2} (\log(\hdelta_n) - \log(\delta))
      \end{bmatrix} \\[4mm]
      \bigl(\int_0^1 \ee^{s\htbB_n} \, \dd s\bigr)^{-1}
      \Bigl(\hoBbeta_n - \oBbeta\Bigr)
      + \Bigl\{\bigl(\int_0^1 \ee^{s\htbB_n} \, \dd s\bigr)^{-1}
               - \bigl(\int_0^1 \ee^{s\tbB} \, \dd s\bigr)^{-1}\Bigr\}
        \oBbeta
     \end{bmatrix}
 \end{equation}
 on the set
 \ $\{\omega \in \Omega
      : (\hvarrho_n(\omega), \hdelta_n(\omega)) \in \RR_{++}^2\}$.
\ Clearly, \eqref{h} implies
 \ $n^{1/2} (\log(\hvarrho_n) - \log(\varrho)) = n^{-1/2} n \hs_n \stoch 0$
 \ as \ $n\to \infty$.
\ Under the assumption
 \ $\|\bc\|^2 + \sum_{i=1}^2 \int_{\cU_2} (z_1 - z_2)^2 \, \mu_i(\dd\bz) > 0$,
 \ \eqref{rhodeltab} implies \ $\hvarrho_n \stoch \varrho$ \ and
 \ $\hdelta_n \stoch \delta$ \ as \ $n \to \infty$,
 \ hence
 \begin{equation}\label{htbB}
   \biggl(\int_0^1 \ee^{s\htbB_n} \, \dd s\biggr)^{-1}
   \stoch
   \biggl(\int_0^1 \ee^{s\tbB} \, \dd s\biggr)^{-1}
   = \frac{1}{2}
     \begin{bmatrix} 1 & 1 \\ 1 & 1 \end{bmatrix}
     + \frac{1}{2\int_0^1 \ee^{(\gamma-\kappa)s} \, \dd s}
       \begin{bmatrix} 1 & -1 \\ -1 & 1 \end{bmatrix}
 \end{equation}
 as \ $n \to \infty$, \ since the function
 \ $\RR_{++}^2 \ni (\varrho, \delta)
    \mapsto \bigl(\int_0^1 \ee^{s\tbB} \, \dd s\bigr)^{-1}$
 \ is continuous, see \eqref{int_ee^stbB}.
Moreover, by the method of the proof of \eqref{h} as above,
 \eqref{rhodeltab} implies
 \begin{align*}
  \begin{bmatrix}
   n^{1/2} (\log(\hdelta_n) - \log(\delta)) \\
   \hoBbeta_n - \oBbeta
  \end{bmatrix}
  = g_n\left(\begin{bmatrix}
              n^{1/2}(\hdelta_n - \delta) \\
              \hoBbeta_n - \oBbeta
             \end{bmatrix}\right)
  \distr
  \begin{bmatrix}
   \frac{\sqrt{1-\delta^2}}{\delta} \,
   \frac{\int_0^1 \cY_t \, \dd\tcW_t}{\int_0^1 \cY_t \, \dd t} \\
   \bcM_1
   - \frac{1}{2} \cI \int_0^1 \cY_t \, \dd t
     \begin{bmatrix} 1 \\ 1 \end{bmatrix}
  \end{bmatrix}
 \end{align*}
 as \ $n \to \infty$, \ with functions
 \ $g : \RR^3 \to \RR^3$ \ and \ $g_n : \RR^3 \to \RR^3$, \ $n \in \NN$,
 \ given by
 \[
   g(\bx)
   := \begin{bmatrix}
       \frac{x_1}{\delta} \\
       x_2 \\
       x_3
      \end{bmatrix} , \quad
   \bx = \begin{bmatrix} x_1 \\ x_2 \\ x_3 \end{bmatrix} \in \RR^3 , \qquad
   g_n(\bx)
   := \begin{bmatrix}
       n^{1/2} \log\bigl( 1 + \frac{x_1}{\delta n^{1/2}} \bigr) \\
       x_2 \\
       x_3
      \end{bmatrix}
 \]
 for \ $\bx = (x_1, x_2, x_3)^\top \in \RR^3$ \ with
 \ $x_1 > - \delta n^{1/2}$, \ and \ $g_n(\bx) := \bzero$ \ otherwise.
Hence, by the continuous mapping theorem, Slutsky's lemma,
 \eqref{gammakappatBbeta} and \eqref{rhodeltab} imply
 \[
   \begin{bmatrix}
    n^{1/2} (\hgamma_n - \gamma) \\
    n^{1/2} (\hkappa_n - \kappa) \\
    \htBbeta_n - \tBbeta
   \end{bmatrix}
   \distr
   \begin{bmatrix}
    \frac{1}{2}
    \begin{bmatrix} 1 & 1 \\ 1 & -1 \end{bmatrix}
    \begin{bmatrix}
     0 \\
     \frac{\sqrt{1-\delta^2}}{\delta} \,
     \frac{\int_0^1 \cY_t \, \dd \tcW_t}{\int_0^1 \cY_t \, \dd t}
    \end{bmatrix} \\
    \bigl(\int_0^1 \ee^{s\tbB} \, \dd s\bigr)^{-1}
    \Biggl(\bcM_1
           - \frac{1}{2} \cI \int_0^1 \cY_t \, \dd t
             \begin{bmatrix} 1 \\ 1 \end{bmatrix}\Biggr)
   \end{bmatrix}
 \]
 as \ $n \to \infty$, \ thus, by \eqref{int_ee^stbB} we obtain
 \eqref{gamma,kappa,tBbeta}.

Under the assumptions
 \ $\|\bc\|^2 + \sum_{i=1}^2 \int_{\cU_2} (z_1 - z_2)^2 \, \mu_i(\dd\bz) = 0$
 \ and \ $\int_{\cU_2} (z_1 - z_2)^2 \, \nu(\dd\bz) > 0$, \ by the method of
 the proof of \eqref{h} as above, \eqref{rhodeltabt} implies
 \begin{align*}
  \begin{bmatrix}
   n^{1/2} (\log(\hdelta_n) - \log(\delta)) \\
   \hoBbeta_n - \oBbeta
  \end{bmatrix}
  = g_n\left(\begin{bmatrix}
              n^{1/2}(\hdelta_n - \delta) \\
              \hoBbeta_n - \oBbeta
             \end{bmatrix}\right)
  \distr
  \begin{bmatrix}
   \frac{\sqrt{1-\delta^2}}{\delta} \,
   \tcW_1 \\
   \frac{1}{2}
   \bigl(\cM_{1,1} + \cM_{1,2} - \cI \int_0^1 \cY_t \, \dd t\bigr)
   \begin{bmatrix} 1 \\ 1 \end{bmatrix}
  \end{bmatrix}
 \end{align*}
 as \ $n \to \infty$.
\ Recall that \ $n^{1/2} (\log(\hvarrho_n) - \log(\varrho)) \stoch 0$ \ as
 \ $n\to \infty$.
\ Hence, by the continuous mapping theorem, Slutsky's lemma,
 \eqref{gammakappatBbeta} and \eqref{rhodeltabt} imply
 \[
   \begin{bmatrix}
    n^{1/2} (\hgamma_n - \gamma) \\
    n^{1/2} (\hkappa_n - \kappa) \\
    \htBbeta_n - \tBbeta
   \end{bmatrix}
   \distr
   \begin{bmatrix}
    \frac{1}{2}
    \begin{bmatrix} 1 & 1 \\ 1 & -1 \end{bmatrix}
    \begin{bmatrix}
     0 \\
     \frac{\sqrt{1-\delta^2}}{\delta} \,
     \tcW_1
    \end{bmatrix} \\
    \bigl(\int_0^1 \ee^{s\tbB} \, \dd s\bigr)^{-1}
    \Biggl(\frac{1}{2}
           \bigl(\cM_{1,1} + \cM_{1,2} - \cI \int_0^1 \cY_t \, \dd t\bigr)
           \begin{bmatrix} 1 \\ 1 \end{bmatrix}\Biggr)
   \end{bmatrix}
 \]
 as \ $n \to \infty$, \ thus, by \eqref{int_ee^stbB} we obtain
 \eqref{gamma,kappa,tBbetat}.

In order to prove \eqref{gamma,kappa,tBbeta1}, first note that
 \begin{equation}\label{gammakappatBbeta1}
   \begin{bmatrix}
    n^{1/2} (\hgamma_n - \gamma) \\
    n^{1/2} (\hkappa_n - \kappa) \\
    n^{1/2} (\htBbeta_n - \tBbeta)
   \end{bmatrix}
   = \begin{bmatrix}
      \frac{1}{2}
      \begin{bmatrix} 1 & 1 \\ 1 & -1 \end{bmatrix}
      \begin{bmatrix}
       n^{1/2} (\log(\hvarrho_n) - \log(\varrho)) \\
       n^{1/2} (\log(\hdelta_n) - \log(\delta))
      \end{bmatrix} \\[4mm]
      \bigl(\int_0^1 \ee^{s\htbB_n} \, \dd s\bigr)^{-1}
      n^{1/2} \Bigl(\hoBbeta_n - \oBbeta\Bigr)
      + \bXi_n \oBbeta
     \end{bmatrix} ,
 \end{equation}
 with
 \[
   \bXi_n
   := n^{1/2}
      \biggl\{\biggl(\int_0^1 \ee^{s\htbB_n} \, \dd s\biggr)^{-1}
              - \biggl(\int_0^1 \ee^{s\tbB} \, \dd s\biggr)^{-1}\biggr\}
 \]
 on the set
 \ $\{\omega \in \Omega
      : (\hvarrho_n(\omega), \hdelta_n(\omega)) \in \RR_{++}^2\}$.
\ Under the assumptions \ $\bc = \bzero$, \ $\bmu = \bzero$ \ and
 \ $\int_{\cU_2} (z_1 - z_2)^2 \, \nu(\dd\bz) > 0$, \ \eqref{rhodeltab1} implies
 \ $\hvarrho_n \stoch \varrho$ \ and \ $\hdelta_n \stoch \delta$ \ as
 \ $n \to \infty$,
 \ hence \eqref{htbB} follows.
By \eqref{int_ee^stbB}, we obtain
 \begin{align*}
  \bXi_n = \frac{n^{1/2}}{2}
             \Biggl(\frac{1}{\int_0^1 (\hvarrho_n)^s \, \dd s} - 1\Biggr)
             \begin{bmatrix} 1 & 1 \\ 1 & 1 \end{bmatrix}
             + \frac{n^{1/2}}{2}
               \Biggl(\frac{1}{\int_0^1 (\hdelta_n)^s \, \dd s}
                      - \frac{1}{\int_0^1 \delta^s \, \dd s}\Biggr)
               \begin{bmatrix} 1 & -1 \\ -1 & 1 \end{bmatrix}
 \end{align*}
 on the set
 \ $\{\omega \in \Omega
      : (\hvarrho_n(\omega), \hdelta_n(\omega)) \in \RR_{++}^2\}$.
\ Here we have
 \[
   n^{1/2} \Biggl(\frac{1}{\int_0^1 (\hvarrho_n(\omega))^s \, \dd s} - 1\Biggr)
   = \begin{cases}
      n^{1/2}
      \bigl(\frac{\log(\hvarrho_n(\omega))}{\hvarrho_n(\omega)-1} - 1\bigr)
      & \text{if \ $\hvarrho_n(\omega) \in \RR_{++} \setminus \{1\}$,} \\
      0 & \text{if \ $\hvarrho_n(\omega) = 1$.}
     \end{cases}
 \]
By the method of the proof of \eqref{h} as above, \eqref{rhodeltab1}
 implies
 \begin{equation}\label{int_varrho}
   n^{1/2} \Biggl(\frac{1}{\int_0^1 (\hvarrho_n)^s \, \dd s} - 1\Biggr)
   = h_n\bigl(n^{3/2}(\hvarrho_n - \varrho)\bigr)
   \stoch 0 , \qquad \text{as \ $n \to \infty$,}
 \end{equation}
 with the functions \ $h_n : \RR \to \RR$, \ $n \in \NN$, \ given by
 \[
   h_n(x)
   := \begin{cases}
       n^{1/2}
       \Bigl(\frac{n^{3/2}\log\bigl(1+\frac{x}{n^{3/2}}\bigr)}{x}
             -1\Bigr)
       & \text{if \ $x \in (-n^{3/2}, \infty) \setminus \{0\}$,} \\
       0 & \text{otherwise,}
      \end{cases}
 \]
 since \ $h_n(x_n) \to 0$ \ as \ $n \to \infty$ \ if \ $x_n \to x$ \ as
 \ $n \to \infty$ \ for all \ $x \in \RR$.
\ Indeed, for all \ $z \in \RR$ \ with \ $|z| \leq \frac{1}{2}$, \ we have
 \[
   \biggl| \log(1+z) - z + \frac{z^2}{2} \biggr|
   = \biggl| - \sum_{k=3}^\infty \frac{(-z)^k}{k} \biggr|
   \leq \sum_{k=3}^\infty \frac{|z|^k}{k}
   \leq \frac{1}{3} \sum_{k=3}^\infty |z|^k
   = \frac{|z|^3}{3(1-|z|)}
   \leq \frac{2}{3} |z|^3 .
 \]
Consequently, if \ $\bigl|\frac{x_n}{n^{3/2}}\bigr| \leq \frac{1}{2}$ \ and
 \ $\frac{x_n}{n^{3/2}} \in (-n^{3/2}, \infty) \setminus \{0\}$ \ then
 \begin{align*}
  \biggl|n^{1/2}
         \biggl(\frac{n^{3/2}\log\bigl(1+\frac{x_n}{n^{3/2}}\bigr)}{x_n}
                -1\biggr)
         - n^{1/2}
           \biggl(\frac{n^{3/2}
                  \bigl(\frac{x_n}{n^{3/2}}
                        - \frac{x_n^2}{2n^3}\bigr)}{x_n}
                  -1\biggr)
  \biggr|
  \leq n^{1/2} \frac{n^{3/2}}{|x_n|} \frac{2}{3}
       \left|\frac{x_n}{n^{3/2}}\right|^3
  \to 0
 \end{align*}
 as \ $n \to \infty$.
\ Recall that \ $n^{1/2} (\log(\hvarrho_n) - \log(\varrho)) \stoch 0$ \ as
 \ $n \to \infty$.
\ Consequently, by \eqref{htbB}, \eqref{rhodeltab1}, \eqref{int_varrho} and Slutsky's
 lemma,
 \begin{equation}\label{appr}
  \begin{bmatrix}
   n^{1/2} (\hgamma_n - \gamma) \\
   n^{1/2} (\hkappa_n - \kappa) \\
   n^{1/2} (\htBbeta_n - \tBbeta)
  \end{bmatrix}
  - \bK
    \begin{bmatrix}
     n^{1/2} (\log(\hdelta_n) - \log(\delta)) \\
     n^{1/2} \Bigl(\frac{1}{\int_0^1 (\hdelta_n)^s \, \dd s}
                    - \frac{1}{\int_0^1 \delta^s \, \dd s}\Bigr) \\[2mm]
     n^{1/2} (\hoBbeta_n - \oBbeta)
    \end{bmatrix}
  \stoch 0
 \end{equation}
 as \ $n \to \infty$ \ with
 \[
   \bK
   := \begin{bmatrix}
       \frac{1}{2} \begin{bmatrix} 1 & 0 \\ -1 & 0 \end{bmatrix}
        & \bzero \\
       \frac{\obeta_1-\obeta_2}{2}
       \begin{bmatrix} 0 & 1 \\ 0 & -1 \end{bmatrix}
        & \bigl(\int_0^1 \ee^{s\tbB} \, \dd s\bigr)^{-1}
      \end{bmatrix}
    = \begin{bmatrix}
       \frac{1}{2} \bv \be_1^\top & \bzero \\
       \frac{\obeta_1-\obeta_2}{2} \bv \be_2^\top
        & \bigl(\int_0^1 \ee^{s\tbB} \, \dd s\bigr)^{-1}
      \end{bmatrix} .
 \]
Indeed, the elements of the sequence in \eqref{appr} take the form
 \[
   \begin{bmatrix}
    \frac{1}{2} n^{1/2} (\log(\hvarrho_n) - \log(\varrho)) \\
    \frac{1}{2} n^{1/2} (\log(\hvarrho_n) - \log(\varrho)) \\
    \Bigl\{\Bigl(\int_0^1 \ee^{s\htbB_n} \, \dd s\Bigr)^{-1}
            - \Bigl(\int_0^1 \ee^{s\tbB} \, \dd s\Bigr)^{-1}\Bigr\}
    n^{1/2} (\hoBbeta_n - \oBbeta)
    + \frac{n^{1/2}}{2}
      \Bigl(\frac{1}{\int_0^1 (\hvarrho_n)^s \, \dd s} - 1\Bigr)
      (\obeta_1 + \obeta_2)
      \begin{bmatrix} 1 \\ 1 \end{bmatrix}
    \end{bmatrix} .
 \]
Moreover,
 \[
   n^{1/2} \Biggl(\frac{1}{\int_0^1 (\hdelta_n)^s \, \dd s}
                  - \frac{1}{\int_0^1 \delta^s \, \dd s}\Biggr)
   = \frac{n^{1/2}(\log(\hdelta_n)-\log(\delta))}{\hdelta_n-1}
     - \log(\delta) \,
       \frac{n^{1/2}(\hdelta_n-\delta)}{(\delta-1)(\hdelta_n-1)}
 \]
 on the set
 \ $\{\omega \in \Omega : \hdelta_n(\omega) \in \RR_{++} \setminus \{1\}\}$.
\ From \eqref{rhodeltab1} we conclude
 \[
   \begin{bmatrix}
    n^{1/2}(\hdelta_n - \delta) \\
    n^{1/2} (\hoBbeta_n - \oBbeta)
   \end{bmatrix}
   \distr
   \cN_3(\bzero, \tbS) ,
   \qquad \text{as \ $n \to \infty$,}
 \]
 with
 \begin{align*}
  \tbS
  &:=\begin{bmatrix}
      0 & 0 \\
      0 & \int_0^1 \int_{\cU_2}
                (\ee^{t\tbB} \bz) (\ee^{t\tbB} \bz)^\top \nu(\dd\bz) \, \dd t
     \end{bmatrix}
     + \frac{3}{4}
       \int_{\cU_2} (z_1 + z_2)^2 \, \nu(\dd\bz)
       \begin{bmatrix}
        0 \\
        \bu
       \end{bmatrix}
       \begin{bmatrix}
        0 \\
        \bu
       \end{bmatrix}^\top \\
  &\:\quad
     + (1 - \delta^2)
       \begin{bmatrix}
        1 \\
        -\frac{\tbeta_1-\tbeta_2}{2\log(\delta^{-1})} \bv
       \end{bmatrix}
       \begin{bmatrix}
        1 \\
        -\frac{\tbeta_1-\tbeta_2}{2\log(\delta^{-1})} \bv
       \end{bmatrix}^\top .
 \end{align*}
Again by the method of the proof of \eqref{h} as above, \eqref{rhodeltab1}
 implies
 \begin{equation}\label{help1}
   \begin{bmatrix}
    n^{1/2} (\log(\hdelta_n) - \log(\delta)) \\
    n^{1/2} \Bigl(\frac{1}{\int_0^1 (\hdelta_n)^s \, \dd s}
                   - \frac{1}{\int_0^1 \delta^s \, \dd s}\Bigr) \\[2mm]
    n^{1/2} (\hoBbeta_n - \oBbeta)
   \end{bmatrix}
   = g_n\left(\begin{bmatrix}
               n^{1/2}(\hdelta_n - \delta) \\
               n^{1/2} (\hoBbeta_n - \oBbeta)
              \end{bmatrix}\right)
   \distr
   \cN_4(\bzero, \ttbS)
 \end{equation}
 as \ $n \to \infty$, \ with the functions \ $g_n : \RR^3 \to \RR^4$,
 \ $n \in \NN$, \ given by
 \[
   g_n(\bx)
   := \begin{bmatrix}
       n^{1/2} \log\bigl(1+\frac{x_1}{\delta n^{1/2}}\bigr) \\[1mm]
       \frac{n^{1/2}\log\bigl(1+\frac{x_1}{\delta n^{1/2}}\bigr)}
            {\delta - 1 + \frac{x_1}{n^{1/2}}}
       - \frac{x_1\log(\delta)}
              {(\delta-1)\bigl(\delta - 1 + \frac{x_1}{n^{1/2}}\bigr)} \\
       x_2 \\
       x_3
      \end{bmatrix}
 \]
 for \ $\bx = (x_1, x_2, x_3)^\top \in \RR^3$ \ with
 \ $x_1 \in (-\delta n^{1/2}, \infty) \setminus \{(1-\delta)n^{1/2}\}$ \ and
 \ $g_n(\bx) := \bzero$ \ otherwise, and with
 \begin{align*}
  \ttbS
  &:=\begin{bmatrix}
      \bzero & \bzero \\
      \bzero & \int_0^1 \int_{\cU_2}
                (\ee^{t\tbB} \bz) (\ee^{t\tbB} \bz)^\top \nu(\dd\bz) \, \dd t
     \end{bmatrix}
     + \frac{3}{4}
       \int_{\cU_2} (z_1 + z_2)^2 \, \nu(\dd\bz)
       \begin{bmatrix}
        \bzero \\
        \bu
       \end{bmatrix}
       \begin{bmatrix}
        \bzero \\
        \bu
       \end{bmatrix}^\top \\
  &\:\quad
     + (1 - \delta^2)
       \begin{bmatrix}
        \frac{1}{\delta}
        \begin{bmatrix}
         1 \\
         \frac{\delta-1-\delta\log(\delta)}{(\delta-1)^2}
        \end{bmatrix} \\[4mm]
        -\frac{\tbeta_1-\tbeta_2}{2\log(\delta^{-1})} \bv
       \end{bmatrix}
       \begin{bmatrix}
        \frac{1}{\delta}
        \begin{bmatrix}
         1 \\
         \frac{\delta-1-\delta\log(\delta)}{(\delta-1)^2}
        \end{bmatrix} \\[4mm]
        -\frac{\tbeta_1-\tbeta_2}{2\log(\delta^{-1})} \bv
       \end{bmatrix}^\top .
 \end{align*}
 since \ $g_n(\bx_n) \to g(\bx)$ \ as \ $n \to \infty$ \ if
 \ $\bx_n \to \bx$ \ as \ $n \to \infty$ \ for all \ $\bx \in \RR^3$,
 \ where the function \ $g : \RR^3 \to \RR^4$ \ is given by
 \ $g(\bx) := \bL \bx$, \ $\bx \in \RR^3$, \ with
 \[
   \bL
   := \begin{bmatrix}
       \frac{1}{\delta} & 0 & 0 \\
       \frac{\delta-1-\delta\log(\delta)}{(\delta-1)^2\delta} & 0 & 0 \\
       0 & 1 & 0 \\
       0 & 0 & 1
      \end{bmatrix} ,
 \]
 and \ $\ttbS = \bL \tbS \bL^\top$.
\ Hence, by the continuous mapping theorem, Slutsky's lemma,  \eqref{appr}
 and \eqref{help1} imply
 \[
   \begin{bmatrix}
    n^{1/2} (\hgamma_n - \gamma) \\
    n^{1/2} (\hkappa_n - \kappa) \\
    n^{1/2} (\htBbeta_n - \tBbeta)
   \end{bmatrix}
   \distr
   \cN_4(\bzero, \tbR)
 \]
 as \ $n \to \infty$, \ where
 \begin{align*}
  \tbR = \bK \ttbS \bK^\top
  &= \begin{bmatrix}
      \bzero & \bzero \\
      \bzero & \bigl(\int_0^1 \ee^{s\tbB} \, \dd s\bigr)^{-1}
               \int_0^1 \int_{\cU_2}
                (\ee^{t\tbB} \bz) (\ee^{t\tbB} \bz)^\top \nu(\dd\bz) \, \dd t
               \bigl(\int_0^1 \ee^{s\tbB} \, \dd s\bigr)^{-1}
     \end{bmatrix} \\
  &\:\quad
     + \frac{3}{4}
       \int_{\cU_2} (z_1 + z_2)^2 \, \nu(\dd\bz)
       \begin{bmatrix}
        \bzero \\
        \bigl(\int_0^1 \ee^{s\tbB} \, \dd s\bigr)^{-1} \bu
       \end{bmatrix}
       \begin{bmatrix}
        \bzero \\
        \bigl(\int_0^1 \ee^{s\tbB} \, \dd s\bigr)^{-1} \bu
       \end{bmatrix}^\top \\
  &\:\quad
     + (1 - \delta^2)
       \begin{bmatrix}
        \frac{1}{2\delta} \bv \\
        \frac{(\obeta_1-\obeta_2)(\delta-1-\delta\log(\delta))}
             {2(\delta-1)^2\delta}
        \bv
        - \frac{\tbeta_1-\tbeta_2}{2\log(\delta^{-1})}
          \bigl(\int_0^1 \ee^{s\tbB} \, \dd s\bigr)^{-1} \bv
       \end{bmatrix} \\
  &\phantom{:=+}
       \times
       \begin{bmatrix}
        \frac{1}{2\delta} \bv \\
        \frac{(\obeta_1-\obeta_2)(\delta-1-\delta\log(\delta))}
             {2(\delta-1)^2\delta}
        \bv
        - \frac{\tbeta_1-\tbeta_2}{2\log(\delta^{-1})}
          \bigl(\int_0^1 \ee^{s\tbB} \, \dd s\bigr)^{-1} \bv
       \end{bmatrix}^\top .
 \end{align*}
By \eqref{Putzer} and \eqref{int_ee^stbB}, for all \ $t \in \RR_+$ \ and
 \ $\bz \in \RR^2$, \ we have
 \begin{align*}
  \biggl(\int_0^1 \ee^{s\tbB} \, \dd s\biggr)^{-1} \ee^{t\tbB} \bz
  &= \Bigl(\frac{1}{2} \bu \bu^\top
           + \frac{\log(\delta)}{2(\delta-1)} \bv \bv^\top\Bigr)
     \Bigl(\frac{1}{2} \bu \bu^\top
           + \frac{\delta^t}{2} \bv \bv^\top\Bigr) \bz \\
  &= \frac{1}{2} (z_1 + z_2) \bu
     + \frac{\delta^t\log(\delta)}{2(\delta-1)} (z_1 - z_2) \bv ,
 \end{align*}
 thus
 \begin{align*}
  &\biggl(\int_0^1 \ee^{s\tbB} \, \dd s\biggr)^{-1}
   \int_0^1 \int_{\cU_2}
   (\ee^{t\tbB} \bz) (\ee^{t\tbB} \bz)^\top \nu(\dd\bz) \, \dd t
   \biggl(\int_0^1 \ee^{s\tbB} \, \dd s\biggr)^{-1} \\
  &\qquad\qquad\qquad\qquad= \frac{1}{4} \int_{\cU_2} (z_1 + z_2)^2 \, \nu(\dd\bz)
      \begin{bmatrix} 1 & 1 \\ 1 & 1 \end{bmatrix}
      + \frac{1}{2} \int_{\cU_2} (z_1^2 - z_2^2) \, \nu(\dd\bz)
        \begin{bmatrix} 1 & 0 \\ 0 & -1 \end{bmatrix} \\
  &\qquad\qquad\qquad\qquad\quad
      + \frac{(\delta^2-1)\log(\delta)}{8(\delta-1)^2}
        \int_{\cU_2} (z_1 - z_2)^2 \, \nu(\dd\bz)
        \begin{bmatrix} 1 & -1 \\ -1 & 1 \end{bmatrix} .
 \end{align*}
Moreover, by \eqref{int_ee^stbB}, we obtain
 \[
   \biggl(\int_0^1 \ee^{s\tbB} \, \dd s\biggr)^{-1} \bu
   = \Bigl(\frac{1}{2} \bu \bu^\top
           + \frac{\log(\delta)}{2(\delta-1)} \bv \bv^\top\Bigr)
     \bu
   = \bu
 \]
 and
 \[
   \biggl(\int_0^1 \ee^{s\tbB} \, \dd s\biggr)^{-1} \bv
   = \Bigl(\frac{1}{2} \bu \bu^\top
           + \frac{\log(\delta)}{2(\delta-1)} \bv \bv^\top\Bigr)
     \bv
   = \frac{\log(\delta)}{\delta-1} \bv .
 \]
Further, by \eqref{ttBbeta}, we obtain
 \[
   \tbeta_1 - \tbeta_2
   = \bv^\top \tBbeta
   = \bv^\top
     \biggl(\int_0^1 \ee^{s\tbB} \, \dd s\biggr)^{-1}
     \oBbeta
   = \frac{\log(\delta)}{\delta-1} \bv^\top \oBbeta
   = \frac{\log(\delta)}{\delta-1} (\obeta_1 - \obeta_2) ,
 \]
 hence
 \ $\obeta_1 - \obeta_2 = \frac{\delta-1}{\log(\delta)} (\tbeta_1 - \tbeta_2)$.
\ Consequently,
 \[
   \frac{(\obeta_1-\obeta_2)(\delta-1-\delta\log(\delta))}
        {2(\delta-1)^2\delta}
   \bv
   - \frac{\tbeta_1-\tbeta_2}{2\log(\delta^{-1})}
     \biggl(\int_0^1 \ee^{s\tbB} \, \dd s\biggr)^{-1} \bv
   = \frac{\tbeta_1-\tbeta_2}{2\delta\log(\delta)} \bv .
 \]
Summarizing, we obtain \ $\tbR = \bR$, \ thus we conclude
 \eqref{gamma,kappa,tBbeta1}.
\proofend

\section{Decomposition of the process}
\label{section_deco}

Let us introduce the sequence
 \[
   U_k := \langle \bu, \bX_k \rangle = X_{k,1} + X_{k,2} , \qquad
   k \in \ZZ_+ ,
 \]
 where \ $\bX_k =: (X_{k,1}, X_{k,2})^\top$.
\ One can observe that \ $U_k \geq 0$ \ for all \ $k \in \ZZ_+$, \ and, by
 \eqref{regr},
 \begin{equation}\label{rec_U}
  U_k = \varrho U_{k-1} + \tvarrho \langle \bu, \tBbeta \rangle
        + \langle \bu, \bM_k \rangle ,
  \qquad k \in \NN ,
 \end{equation}
 where
 \[
   \tvarrho := \frac{1-\varrho}{\log(\varrho^{-1})} ,
 \]
 since
 \ $\langle \bu, \ee^\tbB \bX_{k-1} \rangle = \bu^\top \ee^\tbB \bX_{k-1}
    = \varrho \bu^\top \bX_{k-1} = \varrho U_{k-1}$
 and
 \ $\langle \bu, \oBbeta \rangle
    = \int_0^1 \langle \bu, \ee^{s\tbB} \tBbeta \rangle \, \dd s
    = \int_0^1 \varrho^s \langle \bu, \tBbeta \rangle \, \dd s
    = \frac{1-\varrho}{\log(\varrho^{-1})} \langle \bu, \tBbeta \rangle$,
 \ because \ $\bu$ \ is a left eigenvector of
 \ $\ee^{s \tbB}$, \ $s \in \RR_+$, \ belonging to the eigenvalue
 \ $\varrho^s$.
\ In case of \ $\varrho = 1$, \ $(U_k)_{k \in \ZZ_+}$ \ is a nonnegative unstable
 AR(1) process with
 positive drift \ $\langle \bu, \tBbeta \rangle$ \ and with
 heteroscedastic innovation \ $(\langle \bu, \bM_k \rangle)_{k \in \NN}$.
\ Note that in case of \ $\varrho = 1$, \ the solution of the recursion
 \eqref{rec_U} is
 \begin{equation}\label{rec_U_sol}
   U_k = \sum_{j=1}^k \langle \bu, \bM_j + \tBbeta \rangle , \qquad
   k \in \NN .
 \end{equation}
Moreover, let
 \[
   V_k := \langle \bv, \bX_k \rangle = X_{k,1} - X_{k,2} , \qquad
   k \in \ZZ_+ .
 \]
By \eqref{regr}, we have
 \begin{equation}\label{rec_V}
  V_k = \delta V_{k-1}
        + \tdelta \langle \bv, \tBbeta \rangle
        + \langle \bv, \bM_k \rangle ,
  \qquad k \in \NN ,
 \end{equation}
 where
 \[
   \tdelta := \frac{1-\delta}{\log(\delta^{-1})} ,
 \]
 since
 \ $\langle \bv, \ee^\tbB \bX_{k-1} \rangle = \bv^\top \ee^\tbB \bX_{k-1}
    = \delta \bv^\top \bX_{k-1} = \delta V_{k-1}$
 and
 \ $\langle \bv, \oBbeta \rangle
    = \int_0^1 \langle \bv, \ee^{s\tbB} \tBbeta \rangle \, \dd s
    = \int_0^1 \delta^s \langle \bv, \tBbeta \rangle \, \dd s
    = \frac{1-\delta}{\log(\delta^{-1})} \langle \bv, \tBbeta \rangle$,
 \ because \ $\bv$ \ is a left eigenvector of
 \ $\ee^{s \tbB}$, \ $s \in \RR_+$, \ belonging to the eigenvalue
 \ $\delta^s$.
\ Thus \ $(V_k)_{k \in \ZZ_+}$ \ is a stable AR(1) process with drift
 \ $\tdelta \langle \bv, \tBbeta \rangle$ \ and with
 heteroscedastic innovation \ $(\langle \bv, \bM_k \rangle)_{k \in \NN}$,
 \ since \ $\gamma + \kappa =0$, \ $\gamma \in \RR$ \ and
 \ $\kappa \in \RR_{++}$ \ yield
 \ $\delta = \ee^{\gamma-\kappa} = \ee^{-2\kappa} \in (0, 1)$.
\ Note that the solution of the recursion \eqref{rec_V} is
 \begin{equation}\label{rec_V_sol}
   V_k = \sum_{j=1}^k
          \delta^{k-j}
          \big\langle \bv, \bM_j + \tdelta \, \tBbeta \big\rangle ,
   \qquad k \in \NN .
 \end{equation}
Observe that
 \begin{equation}\label{XUV}
  X_{k,1} = (U_k + V_k)/2 , \qquad X_{k,2} = (U_k - V_k)/2 , \qquad k \in \ZZ_+ .
 \end{equation}

By \eqref{CLSEr}, \eqref{rec_U}, \eqref{CLSEd}, \eqref{rec_V}, \eqref{CLSEb}
 and \eqref{Xregr}, for each \ $n \in \NN$, \ we have
 \begin{gather}
  \hvarrho_n - \varrho
  = \frac{n \sum_{k=1}^n \langle \bu, \bM_k \rangle U_{k-1}
           - \sum_{k=1}^n \langle \bu, \bM_k \rangle \sum_{k=1}^n U_{k-1}}
          {n \sum_{k=1}^n U_{k-1}^2 - \bigl(\sum_{k=1}^n U_{k-1}\bigr)^2} ,
  \label{r-} \\
  \hdelta_n - \delta
  = \frac{n \sum_{k=1}^n \langle \bv, \bM_k \rangle V_{k-1}
           - \sum_{k=1}^n \langle \bv, \bM_k \rangle \sum_{k=1}^n V_{k-1}}
          {n \sum_{k=1}^n V_{k-1}^2 - \bigl(\sum_{k=1}^n V_{k-1}\bigr)^2} ,
  \label{d-} \\
  \hoBbeta_n - \oBbeta
  = \frac{1}{n} \sum_{k=1}^n \bM_k
     - \frac{1}{2n}
       \sum_{k=1}^n
        \begin{bmatrix}
         U_{k-1} & V_{k-1} \\
         U_{k-1} & - V_{k-1}
        \end{bmatrix}
        \begin{bmatrix}
         \hvarrho_n - \varrho \\
         \hdelta_n - \delta
        \end{bmatrix}
  \label{b-}
 \end{gather}
 on the sets \ $H_n$, \ $\tH_n$ \ and \ $H_n \cap \tH_n$, \ respectively.

Theorem \ref{main_rdb} will follow from the following statements by the
 continuous mapping theorem and by Slutsky's lemma, see below.

\begin{Thm}\label{main_Ad}
Under the assumptions of Theorem \ref{main}, we have
 \begin{gather*}
  n^{-3/2} \sum_{k=1}^n V_{k-1}
  \stoch 0 \qquad \text{as \ $n \to \infty$,} \\
   \sum_{k=1}^n
    \begin{bmatrix}
     n^{-2} U_{k-1} \\
     n^{-3} U_{k-1}^2 \\
     n^{-2} V_{k-1}^2 \\
     n^{-1} \bM_k \\
     n^{-2} \langle \bu, \bM_k \rangle U_{k-1} \\
     n^{-3/2} \langle \bv, \bM_k \rangle V_{k-1}
    \end{bmatrix}
   \distr
   \begin{bmatrix}
    \int_0^1 \cY_t \, \dd t \\[1mm]
    \int_0^1 \cY_t^2 \, \dd t \\[1mm]
    (1 - \delta^2)^{-1}
    \langle \tbC \bv, \bv \rangle
    \int_0^1 \cY_t \, \dd t \\
    \bcM_1 \\
    \int_0^1 \cY_t \, \dd \langle \bu, \bcM_t \rangle \\[1mm]
    (1 - \delta^2)^{-1/2}
    \langle \tbC \bv, \bv \rangle
    \int_0^1 \cY_t \, \dd \tcW_t
   \end{bmatrix} \qquad \text{as \ $n \to \infty$.}
 \end{gather*}
\end{Thm}

In case of \ $\langle \tbC \bv, \bv \rangle = 0$ \ the third and
 sixth coordinates of the limit vector in the second convergence of Theorem
 \ref{main_Ad} is 0, thus other scaling factors
 should be chosen for these coordinates, described in the following theorem.

\begin{Thm}\label{maint_Ad}
Suppose that the assumptions of Theorem \ref{main} hold.
If \ $\langle \tbC \bv, \bv \rangle = 0$, \ then
 \begin{gather}\nonumber
  n^{-1} \sum_{k=1}^n V_{k-1}
  \stoch \frac{\tdelta \langle \bv, \tBbeta \rangle}{1-\delta}
  \qquad \text{as \ $n \to \infty$,} \\
  n^{-1} \sum_{k=1}^n V_{k-1}^2
  \stoch \frac{\langle \bV_0 \bv, \bv \rangle}{1-\delta^2}
         + \frac{(\tdelta)^2
                 \langle \bv, \tBbeta \rangle^2}
                {(1-\delta)^2}
         =: M \qquad \text{as \ $n \to \infty$,} \label{M} \\
   \sum_{k=1}^n
    \begin{bmatrix}
     n^{-2} U_{k-1} \\
     n^{-3} U_{k-1}^2 \\
     n^{-1} \langle \bu, \bM_k \rangle \\
     n^{-2} \langle \bu, \bM_k \rangle U_{k-1} \\
     n^{-1/2} \langle \bv, \bM_k \rangle \\
     n^{-1/2} \langle \bv, \bM_k \rangle V_{k-1}
    \end{bmatrix}
   \distr
   \begin{bmatrix}
    \int_0^1 \cY_t \, \dd t \\
    \int_0^1 \cY_t^2 \, \dd t \\
    \langle \bu, \bcM_1 \rangle \\
    \int_0^1
     \cY_t \, \dd \langle \bu, \bcM_t \rangle \\
    \langle \bV_0 \bv, \bv \rangle^{1/2}
    \begin{bmatrix}
     1 & \frac{\tdelta\langle \bv, \tBbeta \rangle}{1-\delta} \\
     \frac{\tdelta\langle \bv, \tBbeta \rangle}{1-\delta} & M
    \end{bmatrix}^{1/2}
    \tbcW_1
   \end{bmatrix} \nonumber
 \end{gather}
 as \ $n \to \infty$, \ where \ $\tbcW_1$ \ is a 2-dimensional random vector
 with standard normal distribution, independent from \ $(\bcW_t)_{t\in\RR_+}$,
 \ and \ $\bV_0$ \ is defined in Proposition \ref{moment_formula_2}.
\end{Thm}

In case of \ $\langle \tbC \bu, \bu \rangle = 0$ \ the third
 and fourth coordinates of the limit vector of the third convergence
 in Theorem \ref{maint_Ad} is 0, since \ $(\cY_t)_{t \in \RR_+}$ \ is the
 deterministic function \ $\cY_t = \langle \bu, \tBbeta \rangle t$,
 \ $t \in \RR_+$ \ (see Remark \ref{REMARK_par}), hence other scaling
 factors should be chosen for these coordinates, as given in the following
 theorem.

\begin{Thm}\label{main1_Ad}
Suppose that the assumptions of Theorem \ref{main} hold.
If \ $\langle \tbC \bu, \bu \rangle = 0$, \ then
 \begin{gather*}
  n^{-2} \sum_{k=1}^n U_{k-1}
  \stoch \frac{\langle \bu, \tBbeta \rangle}{2} ,
  \qquad
  n^{-3} \sum_{k=1}^n U_{k-1}^2
  \stoch \frac{\langle \bu, \tBbeta \rangle^2}{3}
  \qquad \text{as \ $n \to \infty$,} \\
   \sum_{k=1}^n
    \begin{bmatrix}
     n^{-1/2} \langle \bu, \bM_k \rangle \\
     n^{-3/2} \langle \bu, \bM_k \rangle U_{k-1} \\
     n^{-1/2} \langle \bv, \bM_k \rangle \\
     n^{-1/2} \langle \bv, \bM_k \rangle V_{k-1}
    \end{bmatrix}
   \distr
   \cN_4\left( \bzero, \bSigma \right) \qquad \text{as \ $n \to \infty$}
 \end{gather*}
 with
 \begin{align*}
  \bSigma
  &:= \begin{bmatrix}
       \bu^\top \bV_0^{1/2} \\
       \frac{\langle \bu, \tBbeta \rangle}{2}
       \bu^\top \bV_0^{1/2} \\
       \bv^\top \bV_0^{1/2} \\
       \frac{\tdelta\langle \bv, \tBbeta \rangle}{1-\delta}
       \bv^\top \bV_0^{1/2}
      \end{bmatrix}
      \begin{bmatrix}
       \bu^\top \bV_0^{1/2} \\
       \frac{\langle \bu, \tBbeta \rangle}{2}
       \bu^\top \bV_0^{1/2} \\
       \bv^\top \bV_0^{1/2} \\
       \frac{\tdelta\langle \bv, \tBbeta \rangle}{1-\delta}
       \bv^\top \bV_0^{1/2}
      \end{bmatrix}^\top
      + \frac{\langle \bu, \tBbeta \rangle^2
              \bu^\top \bV_0 \bu}
             {12}
        \begin{bmatrix} 0 \\ 1 \\ 0 \\ 0 \end{bmatrix}
        \begin{bmatrix} 0 \\ 1 \\ 0 \\ 0 \end{bmatrix}^\top
      + \frac{(\bv^\top \bV_0 \bv)^2}{1-\delta^2}
        \begin{bmatrix} 0 \\ 0 \\ 0 \\ 1 \end{bmatrix}
        \begin{bmatrix} 0 \\ 0 \\ 0 \\ 1 \end{bmatrix}^\top .
 \end{align*}
\end{Thm}

\noindent
\textbf{Proof of Theorem \ref{main_rdb}.}
The statements about the existence of the estimators \ $\hvarrho_n$
 \ and \ $(\hvarrho_n, \hdelta_n, \hoBbeta_n)$ \ under the given conditions
 follow from Lemma \ref{LEMMA_CLSE_exist_discrete}.

By the continuous mapping theorem and Slutsky's lemma, Theorem \ref{main_Ad}, and
 \eqref{r-} imply \eqref{rho}.
Indeed, by \eqref{r-}, we have
 \begin{equation}\label{r--}
   n(\hvarrho_n - \varrho)
   = \frac{n^{-2} \sum_{k=1}^n \langle \bu, \bM_k \rangle U_{k-1}
           - n^{-1} \sum_{k=1}^n \langle \bu, \bM_k \rangle
             n^{-2} \sum_{k=1}^n U_{k-1}}
          {n^{-3} \sum_{k=1}^n U_{k-1}^2
           - \bigl(n^{-2} \sum_{k=1}^n U_{k-1}\bigr)^2}
 \end{equation}
 on the set \ $H_n$, \ and
 \ $\PP\bigl(\int_0^1 \cY^2_t \, \dd t - (\int_0^1 \cY_t \, \dd t)^2 > 0 \bigr) = 1$,
 \ see the proof of Theorem 3.4 in Barczy et al.\ \cite{BarKorPap}.
Consequently,
 \[
   n(\hvarrho_n - \varrho)
   \distr
   \frac{\int_0^1 \cY_t \, \dd\langle \bu, \bcM_t \rangle
         - \langle \bu, \bcM_1 \rangle \int_0^1 \cY_t \, \dd t}
        {\int_0^1 \cY^2_t \, \dd t - (\int_0^1 \cY_t \, \dd t)^2}
   \qquad \text{as \ $n \to \infty$,}
 \]
 and we obtain \eqref{rho}.

Again by the continuous mapping theorem and Slutsky's lemma, Theorem \ref{main1_Ad},
 and \eqref{r-} imply \eqref{rho1}.
Indeed, by \eqref{r-}, we have
 \begin{equation}\label{r---}
   n^{3/2}(\hvarrho_n - \varrho)
   = \frac{\begin{bmatrix} - n^{-2} \sum_{k=1}^n U_{k-1} \\ 1 \end{bmatrix}^\top
           \begin{bmatrix}
            n^{-1/2} \sum_{k=1}^n \langle \bu, \bM_k \rangle \\
            n^{-3/2} \sum_{k=1}^n \langle \bu, \bM_k \rangle U_{k-1}
           \end{bmatrix}}
          {n^{-3} \sum_{k=1}^n U_{k-1}^2
           - \bigl(n^{-2} \sum_{k=1}^n U_{k-1}\bigr)^2}
 \end{equation}
 on the set \ $H_n$.
\ The first two convergences in Theorem \ref{main1_Ad} imply
 \[
   n^{-3} \sum_{k=1}^n U_{k-1}^2
   - \biggl(n^{-2} \sum_{k=1}^n U_{k-1}\biggr)^2
   \stoch
   \frac{\langle \bu, \tBbeta \rangle^2}{3}
   - \biggl(\frac{\langle \bu, \tBbeta \rangle}{2}\biggr)^2
   = \frac{\langle \bu, \tBbeta \rangle^2}{12} \qquad
   \text{as \ $n \to \infty$.}
 \]
Moreover, the third convergence in Theorem \ref{main1_Ad} implies
 \[
   \begin{bmatrix}
    n^{-1/2} \sum_{k=1}^n \langle \bu, \bM_k \rangle \\
    n^{-3/2} \sum_{k=1}^n \langle \bu, \bM_k \rangle U_{k-1}
   \end{bmatrix}
   \distr
   \cN_2(\bzero, \bSigma_{1,1}) \qquad
   \text{as \ $n \to \infty$,}
 \]
 with
 \begin{align*}
  \bSigma_{1,1}
  := \begin{bmatrix}
       \bu^\top \bV_0^{1/2} \\
       \frac{\langle \bu, \tBbeta \rangle}{2}
       \bu^\top \bV_0^{1/2}
      \end{bmatrix}
      \begin{bmatrix}
       \bu^\top \bV_0^{1/2} \\
       \frac{\langle \bu, \tBbeta \rangle}{2}
       \bu^\top \bV_0^{1/2}
      \end{bmatrix}^\top
      + \frac{\langle \bu, \tBbeta \rangle^2
              \bu^\top \bV_0 \bu}
             {12}
        \begin{bmatrix} 0 \\ 1 \end{bmatrix}
        \begin{bmatrix} 0 \\ 1 \end{bmatrix}^\top
   = \bu^\top \bV_0 \bu
     \begin{bmatrix}
      1 & \frac{\langle \bu, \tBbeta \rangle}{2} \\
      \frac{\langle \bu, \tBbeta \rangle}{2}
       & \frac{\langle \bu, \tBbeta \rangle^2}{3}
     \end{bmatrix} .
 \end{align*}
Consequently,
 \[
   n^{3/2}(\hvarrho_n - \varrho)
   \distr
   \cN(0, \sigma^2)
   \qquad \text{as \ $n \to \infty$}
 \]
 with
 \[
   \sigma^2
   := \biggl(\frac{12}{\langle \bu, \tBbeta \rangle^2}\biggr)^2
      \begin{bmatrix} -\frac{\langle \bu, \tBbeta \rangle}{2} \\ 1 \end{bmatrix}^\top
      \bSigma_{1,1}
      \begin{bmatrix} -\frac{\langle \bu, \tBbeta \rangle}{2} \\ 1 \end{bmatrix}
   = \frac{12}{\langle \bu, \tBbeta \rangle^2} \bu^\top \bV_0 \bu ,
 \]
 and we obtain \eqref{rho1}, since, by Proposition \ref{moment_formula_2},
 \begin{equation}\label{bV_0bu_LEFT,bu_LEFT}
  \begin{aligned}
   \langle \bV_0 \bu, \bu \rangle
   &= \bu^\top
      \int_0^1
       \ee^{u\tbB}
       \biggl( \int_{\cU_2} \bz \bz^\top \nu(\dd \bz) \biggr)
       \ee^{u\tbB^\top}
       \dd u
      \, \bu \\
   &\quad
      + \bu^\top
        \sum_{\ell=1}^2
         \int_0^1
          \biggl( \int_0^{1-u}
                   \langle \ee^{v\tbB} \tBbeta, \be_ \ell \rangle
                   \, \dd v \biggr)
          \ee^{u\tbB} \bC_\ell \ee^{u\tbB^\top}
          \dd u
        \, \bu \\
   &= \int_{\cU_2} \bu^\top \bz \bz^\top \bu \, \nu(\dd \bz)
    = \int_{\cU_2} \langle \bu, \bz \rangle^2 \, \nu(\dd \bz)
    = \int_{\cU_2} (z_1 + z_2)^2 \, \nu(\dd \bz) .
  \end{aligned}
 \end{equation}

In order to prove \eqref{rhodeltab}, first note that
 \begin{equation}\label{tbCbv_LEFT,bv_LEFT=0}
   \langle \tbC \bv , \bv \rangle = 0 \qquad
   \text{if and only if} \qquad
   \|\bc\|^2 + \sum_{i=1}^2 \int_{\cU_2} (z_1 - z_2)^2 \, \mu_i(\dd\bz) = 0 .
 \end{equation}
Indeed, by the spectral mapping theorem, \ $\bv$ \ is a left
 eigenvector of \ $\ee^{s\tbB}$, \ $s \in \RR_+$, \ belonging to the eigenvalue
 \ $\delta^s$ \ and \ $\tbu$ \ is a right eigenvector of \ $\ee^{s\tbB}$,
 \ $s \in \RR_+$, \ belonging to the eigenvalue 1, hence
 \begin{align*}
  \langle \tbC \bv, \bv \rangle
  &= \sum_{i=1}^2 \langle \be_i, \tbu \rangle
      \bv^\top \bV_i \, \bv
   = \sum_{i=1}^2 \langle \be_i, \tbu \rangle
      \sum_{\ell=1}^2
       \int_0^1
        \langle \ee^{(1-u)\tbB} \be_i, \be_\ell \rangle
        \bv^\top \ee^{u\tbB} \bC_\ell \, \ee^{u\tbB^\top} \!\! \bv
        \, \dd u \\
  &= \sum_{i=1}^2 \langle \be_i, \tbu \rangle
      \sum_{\ell=1}^2
       \int_0^1
        \langle \ee^{(1-u)\tbB} \be_i, \be_\ell \rangle
        \delta^{2u} \bv^\top \bC_\ell \, \bv \, \dd u \\
  &= \sum_{\ell=1}^2
      \int_0^1
       \be_\ell^\top \ee^{(1-u)\tbB} \sum_{i=1}^2 \be_i \be_i^\top \tbu
       \delta^{2u} \langle \bC_\ell \bv, \bv \rangle \, \dd u
   = \sum_{\ell=1}^2
      \int_0^1
       \be_\ell^\top \ee^{(1-u)\tbB} \tbu
       \delta^{2u} \langle \bC_\ell \bv, \bv \rangle \, \dd u \\
  &= \sum_{\ell=1}^2
      \be_\ell^\top \tbu \langle \bC_\ell \bv, \bv \rangle
     \int_0^1 \delta^{2u} \, \dd u
   = \langle \obC \bv, \bv \rangle
     \int_0^1 \delta^{2u} \, \dd u
   = \frac{1-\delta^2}{2\log(\delta^{-1})} \,
     \langle \obC \bv, \bv \rangle .
 \end{align*}
Thus \ $\langle \tbC \bv, \bv \rangle = 0$ \ if and only if
 \ $\langle \obC \bv, \bv \rangle = 0$.
\ Recalling
 \[
   \langle \obC \bv, \bv \rangle
   = \sum_{k=1}^2
      \langle \be_k, \tbu \rangle
      \langle \bC_k \bv, \bv \rangle ,
 \]
 one can observe that \ $\langle \obC \bv, \bv \rangle = 0$ \ if and
 only if
 \ $\langle \bC_k \bv, \bv \rangle
    = 2 c_k + \int_{\cU_2} \langle \bv, \bz \rangle^2 \, \mu_k(\dd\bz) = 0$
 \ for each \ $k \in \{1, 2\}$, \ which is equivalent to \ $\bc = \bzero$
 \ and \ $\int_{\cU_2} (z_1 - z_2)^2 \, \mu_k(\dd\bz) = 0$ \ for each
 \ $k \in \{1, 2\}$.

By the continuous mapping theorem and Slutsky's lemma, Theorem \ref{main_Ad},
 \eqref{r-}, \eqref{d-} and \eqref{b-} imply \eqref{rhodeltab}.
Indeed, by \eqref{d-} and \eqref{b-}, we have
 \[
   n^{1/2}(\hdelta_n - \delta)
   = \frac{n^{-3/2} \sum_{k=1}^n \langle \bv, \bM_k \rangle V_{k-1}
           - n^{-1} \sum_{k=1}^n \langle \bv, \bM_k \rangle
             n^{-3/2} \sum_{k=1}^n V_{k-1}}
          {n^{-2} \sum_{k=1}^n V_{k-1}^2
           - \bigl(n^{-3/2} \sum_{k=1}^n V_{k-1}\bigr)^2}
 \]
 on the set \ $\tH_n$, \ and
 \[
   \hoBbeta_n - \oBbeta
   = \frac{1}{n} \sum_{k=1}^n \bM_k
     - \frac{1}{2}
       \begin{bmatrix}
        n^{-2} \sum_{k=1}^n U_{k-1} & n^{-3/2} \sum_{k=1}^n V_{k-1} \\
        n^{-2} \sum_{k=1}^n U_{k-1} & - n^{-3/2} \sum_{k=1}^n V_{k-1}
       \end{bmatrix}
       \begin{bmatrix}
        n(\hvarrho_n - \varrho) \\
        n^{1/2}(\hdelta_n - \delta)
       \end{bmatrix}
 \]
 on the set \ $H_n \cap \tH_n$.
\ Recalling \eqref{r--} and taking into account
 \ $n^{-3/2} \sum_{k=1}^n V_{k-1} \stoch 0$ \ as \ $n \to \infty$
 \ (see Lemma \ref{main_Vt}),
 \ $\PP\bigl(\int_0^1 \cY^2_t \, \dd t - (\int_0^1 \cY_t \, \dd t)^2 > 0 \bigr) = 1$
 \ (see the proof of Theorem 3.4 in Barczy et al.\ \cite{BarKorPap}),
 \ $\langle \tbC \bv , \bv \rangle > 0$ \ and
 \ $\PP\bigl(\int_0^1 \cY_t \, \dd t > 0 \bigr) = 1$, \ we obtain
 \eqref{rhodeltab}.

In order to prove \eqref{rhodeltabt}, first note that, under the
 additional condition \ $\langle \tbC \bv , \bv \rangle = 0$,
 \ we have
 \begin{equation}\label{bV_0bv_LEFT,bv_LEFT=0}
  \langle \bV_0 \bv , \bv \rangle = 0 \qquad
  \text{if and only if} \qquad
  \int_{\cU_2} (z_1 - z_2)^2 \, \nu(\dd\bz) = 0 ,
 \end{equation}
 since, by Proposition \ref{moment_formula_2},
 \begin{equation}\label{bV_0bv_LEFT,bv_LEFT}
  \begin{aligned}
   \langle \bV_0 \bv, \bv \rangle
   &= \bv^\top
      \int_0^1
       \ee^{u\tbB}
       \biggl( \int_{\cU_2} \bz \bz^\top \nu(\dd \bz) \biggr)
       \ee^{u\tbB^\top}
       \dd u
      \, \bv \\
   &\quad
      + \bv^\top
        \sum_{\ell=1}^2
         \int_0^1
          \biggl( \int_0^{1-u}
                   \langle \ee^{v\tbB} \tBbeta, \be_ \ell \rangle
                   \, \dd v \biggr)
          \ee^{u\tbB} \bC_\ell \ee^{u\tbB^\top}
          \dd u
        \, \bv \\
   &= \biggl( \int_0^1 \delta^{2u} \, \dd u \biggr)
      \int_{\cU_2} \bv^\top \bz \bz^\top \bv \, \nu(\dd \bz) \\
   &= \frac{1-\delta^2}{2\log(\delta^{-1})}
      \int_{\cU_2} \langle \bv, \bz \rangle^2 \, \nu(\dd \bz)
    = \frac{1-\delta^2}{2\log(\delta^{-1})}
      \int_{\cU_2} (z_1 - z_2)^2 \, \nu(\dd \bz) .
  \end{aligned}
 \end{equation}
By the continuous mapping theorem and Slutsky's lemma, Theorem \ref{maint_Ad},
 \eqref{r-}, \eqref{d-} and \eqref{b-} imply \eqref{rhodeltabt}.
Indeed, by \eqref{d-} and \eqref{b-}, we have
 \begin{equation}\label{d--}
   n^{1/2}(\hdelta_n - \delta)
   = \frac{\begin{bmatrix} - n^{-1} \sum_{k=1}^n V_{k-1} \\ 1 \end{bmatrix}^\top
           \begin{bmatrix}
            n^{-1/2} \sum_{k=1}^n \langle \bv, \bM_k \rangle \\
            n^{-1/2} \sum_{k=1}^n \langle \bv, \bM_k \rangle V_{k-1}
           \end{bmatrix}}
          {n^{-1} \sum_{k=1}^n V_{k-1}^2
           - \bigl(n^{-1} \sum_{k=1}^n V_{k-1}\bigr)^2}
 \end{equation}
 on the set \ $\tH_n$, \ and
 \begin{align*}
  \hoBbeta_n - \oBbeta
  &= \frac{1}{2}
     \begin{bmatrix} 1 & 1 \\ 1 & -1 \end{bmatrix}
     \begin{bmatrix}
      n^{-1} \sum_{k=1}^n \langle \bu, \bM_k \rangle \\
      n^{-1} \sum_{k=1}^n \langle \bv, \bM_k \rangle
     \end{bmatrix} \\
  &\quad
     - \frac{1}{2}
       \begin{bmatrix}
        n^{-2} \sum_{k=1}^n U_{k-1} & n^{-3/2} \sum_{k=1}^n V_{k-1} \\
        n^{-2} \sum_{k=1}^n U_{k-1} & - n^{-3/2} \sum_{k=1}^n V_{k-1}
       \end{bmatrix}
       \begin{bmatrix}
        n(\hvarrho_n - \varrho) \\
        n^{1/2}(\hdelta_n - \delta)
       \end{bmatrix}
 \end{align*}
 on the set \ $H_n \cap \tH_n$.
\ Recalling \eqref{r--} and taking into account the first two convergences in
 Theorem \ref{maint_Ad},
 \ $\PP\bigl(\int_0^1 \cY^2_t \, \dd t - (\int_0^1 \cY_t \, \dd t)^2 > 0 \bigr) = 1$
 \ and \ $\langle \bV_0 \bv , \bv \rangle > 0$, \ we obtain
 \[
   \begin{bmatrix}
    n (\hvarrho_n - \varrho) \\
    n^{1/2} (\hdelta_n - \delta) \\
    \hoBbeta_n - \oBbeta
   \end{bmatrix}
   \distr
   \begin{bmatrix}
    \cI \\
    \cJ \\
    \frac{1}{2} \langle \bu, \bcM_1 \rangle
    \begin{bmatrix} 1 \\ 1 \end{bmatrix}
    - \frac{1}{2} \cI \int_0^1 \cY_t \, \dd t
      \begin{bmatrix} 1 \\ 1 \end{bmatrix}
   \end{bmatrix} ,
 \]
 as \ $n \to \infty$ \ with
 \[
   \cJ := \frac{1-\delta^2}{\langle \bV_0 \bv, \bv \rangle}
          \begin{bmatrix}
           - \frac{\tdelta\langle \bv, \tBbeta \rangle}{1-\delta} \\
           1
          \end{bmatrix}^\top
          \langle \bV_0 \bv, \bv \rangle^{1/2}
          \begin{bmatrix}
           1 & \frac{\tdelta\langle \bv, \tBbeta \rangle}{1-\delta} \\
           \frac{\tdelta\langle \bv, \tBbeta \rangle}{1-\delta} & M
          \end{bmatrix}^{1/2}
          \tbcW_1 .
 \]
Calculating the variance, it is easy to check that
 \ $\cJ \distre \sqrt{1-\delta^2} \, \tcW_1$, \ hence we conclude
 \eqref{rhodeltabt}.

In order to prove \eqref{rhodeltab1}, first note that
 \ $\langle \tbC \bu , \bu \rangle = 0$ \ if and only if
 \ $\bc = \bzero$ \ and \ $\bmu = \bzero$.
\ Indeed, by Remark \ref{REMARK_par}, we have
 \ $\langle \tbC \bu , \bu \rangle
    = \langle \obC \bu , \bu \rangle$,
 \ and \ $\langle \obC \bu , \bu \rangle = 0$ \ if and only if
 \ $\bc = \bzero$ \ and \ $\bmu = \bzero$.
\ Hence, \ $\langle \tbC \bu , \bu \rangle = 0$ \ or
 \ $\langle \obC \bu , \bu \rangle = 0$ \ implies
 \ $\langle \tbC \bv , \bv \rangle = 0$ \ and
 \ $\langle \obC \bv , \bv \rangle = 0$ \ as well.

Consequently, under the additional condition
 \ $\langle \tbC \bu , \bu \rangle = 0$, \ we have
 \ $\langle \bV_0 \bv, \bv \rangle = 0$ \ if and only if
 \ $\int_{\cU_2} (z_1 - z_2)^2 \, \nu(\dd \bz) = 0$, \ see
 \eqref{bV_0bv_LEFT,bv_LEFT=0}.

By \eqref{b-}, we have
 \begin{align*}
  n^{1/2}(\hoBbeta_n - \oBbeta)
  &= \frac{1}{2}
     \begin{bmatrix} 1 & 1 \\ 1 & -1 \end{bmatrix}
     \begin{bmatrix}
      n^{-1/2} \sum_{k=1}^n \langle \bu, \bM_k \rangle \\
      n^{-1/2} \sum_{k=1}^n \langle \bv, \bM_k \rangle
     \end{bmatrix} \\
  &\quad
     - \frac{1}{2}
       \begin{bmatrix}
        n^{-2} \sum_{k=1}^n U_{k-1} & n^{-1} \sum_{k=1}^n V_{k-1} \\
        n^{-2} \sum_{k=1}^n U_{k-1} & - n^{-1} \sum_{k=1}^n V_{k-1}
       \end{bmatrix}
       \begin{bmatrix}
        n^{3/2}(\hvarrho_n - \varrho) \\
        n^{1/2}(\hdelta_n - \delta)
       \end{bmatrix}
 \end{align*}
 on the set \ $H_n \cap \tH_n$.
\ Recalling \eqref{r---} and \eqref{d--}, we obtain
 \[
   \begin{bmatrix}
    n^{3/2} (\hvarrho_n - \varrho) \\
    n^{1/2} (\hdelta_n - \delta) \\
    n^{1/2} (\hoBbeta_n - \oBbeta)
   \end{bmatrix}
   = \begin{bmatrix}
      \bA_{1,1}^{(n)} & \bA_{1,2}^{(n)} \\
      \bA_{2,1}^{(n)} & \bA_{2,2}^{(n)}
     \end{bmatrix}
     \begin{bmatrix}
      n^{-1/2} \sum_{k=1}^n \langle \bu, \bM_k \rangle \\
      n^{-3/2} \sum_{k=1}^n \langle \bu, \bM_k \rangle U_{k-1} \\
      n^{-1/2} \sum_{k=1}^n \langle \bv, \bM_k \rangle \\
      n^{-1/2} \sum_{k=1}^n \langle \bv, \bM_k \rangle V_{k-1}
     \end{bmatrix} ,
 \]
 with
 \begin{gather*}
  \bA_{1,1}^{(n)}
  := \frac{\begin{bmatrix}
            - n^{-2} \sum\limits_{k=1}^n U_{k-1} & 1 \\
            0 & 0
           \end{bmatrix}}
          {n^{-3} \sum\limits_{k=1}^n U_{k-1}^2
           - \bigl(n^{-2} \sum\limits_{k=1}^n U_{k-1}\bigr)^2} , \qquad
  \bA_{1,2}^{(n)}
  := \frac{\begin{bmatrix}
            0 & 0 \\
            - n^{-1} \sum\limits_{k=1}^n V_{k-1} & 1
           \end{bmatrix}}
          {n^{-1} \sum\limits_{k=1}^n V_{k-1}^2
           - \bigl(n^{-1} \sum\limits_{k=1}^n V_{k-1}\bigr)^2} , \\
  \bA_{2,1}^{(n)}
  := \frac{\begin{bmatrix}
            n^{-3} \sum\limits_{k=1}^n U_{k-1}^2
             & - n^{-2} \sum\limits_{k=1}^n U_{k-1} \\
            n^{-3} \sum\limits_{k=1}^n U_{k-1}^2
             & - n^{-2} \sum\limits_{k=1}^n U_{k-1}
           \end{bmatrix}}
           {2 n^{-3} \sum\limits_{k=1}^n U_{k-1}^2
            - 2 \bigl(n^{-2} \sum\limits_{k=1}^n U_{k-1}\bigr)^2} , \qquad
  \bA_{2,2}^{(n)}
  := \frac{\begin{bmatrix}
            n^{-1} \sum\limits_{k=1}^n V_{k-1}^2
             & - n^{-1} \sum\limits_{k=1}^n V_{k-1} \\
            - n^{-1} \sum\limits_{k=1}^n V_{k-1}^2
             & n^{-1} \sum\limits_{k=1}^n V_{k-1}
           \end{bmatrix}}
          {2 n^{-1} \sum\limits_{k=1}^n V_{k-1}^2
           - 2 \bigl(n^{-1} \sum\limits_{k=1}^n V_{k-1}\bigr)^2} .
 \end{gather*}
The first two convergences in Theorem \ref{maint_Ad} hold, since the
 assumption \ $\langle \tbC \bu, \bu \rangle = 0$ \ implies
 \ $\langle \tbC \bv, \bv \rangle = 0$.
\ Since \ $\tBbeta \ne \bzero$, \ we have
 \ $\langle \bu, \tBbeta \rangle > 0$.
\ Moreover, as we have already proved, the assumption
 \ $\|\bc\|^2 + \sum_{i=1}^2 \int_{\cU_2} (z_1 - z_2)^2 \, \mu_i(\dd\bz) > 0$
  implies \ $\langle \tbC \bv , \bv \rangle > 0$.
\ Hence Theorem \ref{main1_Ad} and Lemmas \ref{main_Vt} and \ref{main_VVt} imply
 \[
   \begin{bmatrix}
    \bA_{1,1}^{(n)} & \bA_{1,2}^{(n)} \\
    \bA_{2,1}^{(n)} & \bA_{2,2}^{(n)}
   \end{bmatrix}
   \stoch
   \begin{bmatrix}
    \bA_{1,1} & \bA_{1,2} \\
    \bA_{2,1} & \bA_{2,2}
   \end{bmatrix} \qquad \text{as \ $n \to \infty$,}
 \]
 with
 \begin{gather*}
  \bA_{1,1} := \frac{6}{\langle \bu, \tBbeta \rangle^2}
               \begin{bmatrix}
                - \langle \bu, \tBbeta \rangle & 2 \\
                0 & 0
               \end{bmatrix} , \qquad
  \bA_{1,2} := \frac{1-\delta^2}{\bv^\top \bV_0 \bv}
               \begin{bmatrix}
                0 & 0 \\
                -\frac{\tdelta\langle \bv, \tBbeta \rangle}
                      {1-\delta}
                 & 1
               \end{bmatrix} , \\
  \bA_{2,1} := \frac{1}{\langle \bu, \tBbeta \rangle}
               \begin{bmatrix}
                2 \langle \bu, \tBbeta \rangle & -3 \\
                2 \langle \bu, \tBbeta \rangle & -3
               \end{bmatrix} , \qquad
  \bA_{2,2} := \frac{1-\delta^2}{2 \bv^\top \bV_0 \bv}
                \begin{bmatrix}
                 M & -\frac{\tdelta\langle \bv, \tBbeta \rangle}
                           {1-\delta} \\[2mm]
                 -M & \frac{\tdelta\langle \bv, \tBbeta \rangle}
                           {1-\delta}
   \end{bmatrix} .
 \end{gather*}
By Theorem \ref{main1_Ad} and the continuous mapping theorem, we have
 \[
   \begin{bmatrix}
    n^{3/2} (\hvarrho_n - \varrho) \\
    n^{1/2} (\hdelta_n - \delta) \\
    n^{1/2} (\hoBbeta_n - \oBbeta)
   \end{bmatrix}
   \distr
   \cN_4(\bzero, \bT) \qquad \text{as \ $n \to \infty$,}
 \]
 with
 \[
   \bT
   := \begin{bmatrix}
       \bA_{1,1} & \bA_{1,2} \\
       \bA_{2,1} & \bA_{2,2}
      \end{bmatrix}
      \bSigma
      \begin{bmatrix}
       \bA_{1,1} & \bA_{1,2} \\
       \bA_{2,1} & \bA_{2,2}
      \end{bmatrix}^\top .
 \]
Here \ $\bT$ \ takes the form
 \begin{align*}
  \bT
  = \begin{bmatrix}
       \bzero & \bzero \\
       \bzero & \bV_0
      \end{bmatrix}
      + \frac{3\bu^\top \bV_0 \bu}{4}
        \begin{bmatrix}
         \frac{4}{\langle \bu, \tBbeta \rangle} \be_1 \\
         -\bu
        \end{bmatrix}
        \begin{bmatrix}
         \frac{4}{\langle \bu, \tBbeta \rangle} \be_1 \\
         -\bu
        \end{bmatrix}^\top
      + (1 - \delta^2)
        \begin{bmatrix}
         \be_2 \\
         -\frac{\tdelta\langle \bv, \tBbeta \rangle}{2(1-\delta)}
          \bv
        \end{bmatrix}
        \begin{bmatrix}
         \be_2 \\
         -\frac{\tdelta\langle \bv, \tBbeta \rangle}{2(1-\delta)}
          \bv
        \end{bmatrix}^\top ,
 \end{align*}
 where \ $\be_1 = \begin{bmatrix} 1 \\ 0 \end{bmatrix}$ \ and
 \ $\be_2 = \begin{bmatrix} 0 \\ 1 \end{bmatrix}$.
\ Indeed,
 \[
   \begin{bmatrix}
    \bA_{1,1} & \bA_{1,2} \\
    \bA_{2,1} & \bA_{2,2}
   \end{bmatrix}
   \begin{bmatrix}
    \bu^\top \bV_0^{1/2} \\
    \frac{\langle \bu, \tBbeta \rangle}{2}
    \bu^\top \bV_0^{1/2} \\
    \bv^\top \bV_0^{1/2} \\
    \frac{\tdelta\langle \bv, \tBbeta \rangle}{1-\delta}
    \bv^\top \bV_0^{1/2}
   \end{bmatrix}
   = \begin{bmatrix}
      \bzero^\top \\
      \bzero^\top \\
      \frac{1}{2} \bu^\top \bV_0^{1/2}
      + \frac{1}{2} \bv^\top \bV_0^{1/2} \\
      \frac{1}{2} \bu^\top \bV_0^{1/2}
      - \frac{1}{2} \bv^\top \bV_0^{1/2}
     \end{bmatrix}
   = \begin{bmatrix}
      \bzero^\top \\
      \bzero^\top \\
      \be_1^\top \bV_0^{1/2} \\
      \be_2^\top \bV_0^{1/2}
     \end{bmatrix} ,
 \]
 hence
 \begin{align*}
  &\begin{bmatrix}
    \bA_{1,1} & \bA_{1,2} \\
    \bA_{2,1} & \bA_{2,2}
   \end{bmatrix}
   \begin{bmatrix}
    \bu^\top \bV_0^{1/2} \\
    \frac{\langle \bu, \tBbeta \rangle}{2}
    \bu^\top \bV_0^{1/2} \\
    \bv^\top \bV_0^{1/2} \\
    \frac{\tdelta\langle \bv, \tBbeta \rangle}{1-\delta}
    \bv^\top \bV_0^{1/2}
   \end{bmatrix}
   \begin{bmatrix}
    \bu^\top \bV_0^{1/2} \\
    \frac{\langle \bu, \tBbeta \rangle}{2}
    \bu^\top \bV_0^{1/2} \\
    \bv^\top \bV_0^{1/2} \\
    \frac{\tdelta\langle \bv, \tBbeta \rangle}{1-\delta}
    \bv^\top \bV_0^{1/2}
   \end{bmatrix}^\top
   \begin{bmatrix}
    \bA_{1,1} & \bA_{1,2} \\
    \bA_{2,1} & \bA_{2,2}
   \end{bmatrix}^\top \\
  &\qquad\qquad
   = \begin{bmatrix}
      \bzero^\top \\
      \bzero^\top \\
      \be_1^\top \bV_0^{1/2} \\
      \be_2^\top \bV_0^{1/2}
     \end{bmatrix}
     \begin{bmatrix}
      \bzero^\top \\
      \bzero^\top \\
      \be_1^\top \bV_0^{1/2} \\
      \be_2^\top \bV_0^{1/2}
     \end{bmatrix}^\top
   = \begin{bmatrix}
      \bzero & \bzero \\
      \bzero & \bV_0
     \end{bmatrix} .
 \end{align*}
Moreover,
 \[
   \begin{bmatrix}
    \bA_{1,1} & \bA_{1,2} \\
    \bA_{2,1} & \bA_{2,2}
   \end{bmatrix}
   \begin{bmatrix} 0 \\ 1 \\ 0 \\ 0 \end{bmatrix}
   = \begin{bmatrix}
      \frac{12}{\langle \bu, \tBbeta \rangle^2} \be_1 \\
      - \frac{3}{\langle \bu, \tBbeta \rangle} \bu
     \end{bmatrix} , \qquad
   \begin{bmatrix}
    \bA_{1,1} & \bA_{1,2} \\
    \bA_{2,1} & \bA_{2,2}
   \end{bmatrix}
   \begin{bmatrix} 0 \\ 0 \\ 0 \\ 1 \end{bmatrix}
   = \begin{bmatrix}
      \frac{1-\delta^2}{\bv^\top\bV_0\bv} \be_2 \\[2mm]
      - \frac{(1-\delta^2)\tdelta\langle\bv,\tBbeta\rangle}
             {2(1-\delta)\bv^\top\bV_0\bv}
        \bv
     \end{bmatrix} .
 \]
As in the proof of the formula for \ $\bu^\top \bV_0 \bu$ \ in
 \eqref{bV_0bu_LEFT,bu_LEFT}, using \ $\bC_1 = \bzero$ \ and
 \ $\bC_2 = \bzero$, \ we obtain
 \[
   \bV_0 = \int_0^1 \int_{\cU_2}
            (\ee^{t\tbB} \bz) (\ee^{t\tbB} \bz)^\top \nu(\dd\bz) \, \dd t ,
 \]
 hence \ $\bT = \bS$, \ and we conclude \eqref{rhodeltab1}.
\proofend

\section{Proof of Theorem \ref{main_Ad}}
\label{section_proof_main}

Consider the sequence of stochastic processes
 \[
   \bcZ^{(n)}_t
   := \begin{bmatrix}
       \bcM_t^{(n)} \\
       \bcN_t^{(n)} \\
       \bcP_t^{(n)}
      \end{bmatrix}
   := \sum_{k=1}^\nt
       \bZ^{(n)}_k ,
 \]
 with
 \[
   \bZ^{(n)}_k
   := \begin{bmatrix}
       n^{-1} \bM_k \\
       n^{-2} \bM_k U_{k-1} \\
       n^{-3/2} \bM_k V_{k-1}
      \end{bmatrix}
    = \begin{bmatrix}
       n^{-1} \\
       n^{-2} U_{k-1} \\
       n^{-3/2} V_{k-1}
      \end{bmatrix}
      \otimes \bM_k
 \]
 for \ $t \in \RR_+$ \ and \ $k, n \in \NN$, \ where \ $\otimes$ \ denotes
 Kronecker product of matrices.
Theorem \ref{main_Ad} follows from Lemma \ref{main_VV} and the following
 theorem (this will be explained after Theorem \ref{main_conv}).

\begin{Thm}\label{main_conv}
Under the assumptions of Theorem \ref{main}, we have
 \begin{equation}\label{conv_Z}
   \bcZ^{(n)} \distr \bcZ , \qquad \text{as \ $n\to\infty$,}
 \end{equation}
 where the process \ $(\bcZ_t)_{t \in \RR_+}$ \ with values in \ $(\RR^2)^3$ \ is
 the pathwise unique strong solution of the SDE
 \begin{equation}\label{ZSDE}
  \dd \bcZ_t
  = \gamma(t, \bcZ_t) \begin{bmatrix} \dd \bcW_t \\ \dd \tbcW_t \end{bmatrix} ,
  \qquad t \in \RR_+ ,
 \end{equation}
 with initial value \ $\bcZ_0 = \bzero$, \ where \ $(\bcW_t)_{t \in \RR_+}$ \ and
 \ $(\tbcW_t)_{t \in \RR_+}$ \ are independent 2-dimensional standard Wiener
 processes, and \ $\gamma : \RR_+ \times (\RR^2)^3 \to (\RR^{2\times2})^{3\times2}$
 \ is defined by
 \[
   \gamma(t, \bx)
   := \begin{bmatrix}
       (\langle \bu, \bx_1 + t \tBbeta \rangle^+)^{1/2} \, \tbC^{1/2}
        & \bzero \\
       (\langle \bu, \bx_1 + t \tBbeta \rangle^+)^{3/2} \, \tbC^{1/2}
        & \bzero \\
       \bzero & \bigl(\frac{\langle \tbC \bv, \bv \rangle}
                           {1 - \delta^2}\bigr)^{1/2}
           \langle \bu, \bx_1 + t \tBbeta \rangle \, \tbC^{1/2}
      \end{bmatrix}
 \]
 for \ $t \in \RR_+$ \ and
 \ $\bx = (\bx_1^\top , \bx_2^\top , \bx_3^\top)^\top \in (\RR^2)^3$.
\end{Thm}
(Note that the statement of Theorem \ref{main_conv} holds even if
 \ $\langle \tbC \bv, \bv \rangle = 0$, \ when the last
 2-dimensional coordinate process of the pathwise unique strong solution
 \ $(\bcZ_t)_{t\in\RR_+}$ \ is \ $\bzero$.)

The SDE \eqref{ZSDE} has the form
 \begin{align}\label{MNPSDE}
  \dd \bcZ_t
  =: \begin{bmatrix}
      \dd \bcM_t \\
      \dd \bcN_t \\
      \dd \bcP_t
     \end{bmatrix}
  = \begin{bmatrix}
     (\langle \bu, \bcM_t + t \tBbeta \rangle^+)^{1/2} \, \tbC^{1/2} \,
     \dd \bcW_t \\[1mm]
     (\langle \bu, \bcM_t + t \tBbeta \rangle^+)^{3/2} \, \tbC^{1/2} \,
     \dd \bcW_t \\[1mm]
     \bigl(\frac{\langle \tbC \bv, \bv \rangle}
                {1 - \delta^2}\bigr)^{1/2}
     \langle \bu, \bcM_t + t \tBbeta \rangle \, \tbC^{1/2} \,
     \dd \tbcW_t
    \end{bmatrix} , \qquad t \in \RR_+ .
 \end{align}
One can prove that the first 2-dimensional equation of the SDE
 \eqref{MNPSDE} has a pathwise unique strong solution
 \ $(\bcM_t^{(\by_0)})_{t\in\RR_+}$ \ with arbitrary initial value
 \ $\bcM_0^{(\by_0)} = \by_0 \in \RR^2$.
\ Indeed, it is equivalent to the existence of a pathwise unique strong
 solution of the SDE
 \begin{equation}\label{SDE_P_Q}
  \begin{cases}
   \dd \cS_t
   = \langle \bu, \tBbeta \rangle \, \dd t
     + (\cS_t^+)^{1/2} \, \bu^\top \tbC^{1/2} \, \dd \bcW_t, \\[2mm]
   \dd \bcQ_t
   = - \bPi \tBbeta \, \dd t
     + (\cS_t^+)^{1/2} \, \bigl( \bI_2 - \bPi \bigr) \tbC^{1/2} \, \dd \bcW_t ,
  \end{cases}
  \qquad t \in \RR_+ ,
 \end{equation}
 with initial value
 \ $\bigl(\cS_0^{(\by_0)}, \, \bcQ_0^{(\by_0)}\bigr)
    =\big(\langle \bu, \by_0 \rangle, \, (\bI_2 - \bPi)\by_0\bigr)
    \in \RR \times \RR^2$,
 \ where \ $\bI_2$ \ denotes the $2$-dimensional unit matrix and
 \ $\bPi := \tbu \bu^\top$, \ since we have the correspondences
 \begin{gather*}
  \cS_t^{(\by_0)} = \bu^\top (\bcM_t^{(\by_0)} + t \tBbeta) , \qquad
  \bcQ_t^{(\by_0)} = \bcM_t^{(\by_0)} - \cS_t^{(\by_0)} \tbu \\
  \bcM_t^{(\by_0)} = \bcQ_t^{(\by_0)} + \cS_t^{(\by_0)} \tbu ,
 \end{gather*}
 see the proof of Isp\'any and Pap \cite[Theorem 3.1]{IspPap2}.
By Remark \ref{REMARK_SDE}, \ $\cS_t^+$ \  may be replaced by
 \ $\cS_t$ \ for all \ $t \in \RR_+$ \ in the first equation of
 \eqref{SDE_P_Q} provided that \ $\langle \bu, \by_0 \rangle \in \RR_+$,
 \ hence \ $\langle \bu, \bcM_t + t \tBbeta \rangle^+$ \ may be
 replaced by \ $\langle \bu, \bcM_t + t \tBbeta \rangle$ \ for
 all \ $t \in \RR_+$ \ in \eqref{MNPSDE}.
Thus the SDE \eqref{ZSDE} has a pathwise unique strong solution with initial
 value \ $\bcZ_0 = \bzero$, \ and we have
 \[
   \bcZ_t
   = \begin{bmatrix}
      \bcM_t \\
      \bcN_t \\
      \bcP_t
     \end{bmatrix}
   = \begin{bmatrix}
      \int_0^t
       \langle \bu, \bcM_s + s \tBbeta \rangle^{1/2} \, \tbC^{1/2}
       \, \dd \bcW_s \\[1mm]
      \int_0^t
       \langle \bu, \bcM_s + s \tBbeta \rangle \, \dd \bcM_s \\[1mm]
      \left(\frac{\langle \tbC \bv, \bv \rangle}
                 {1- \delta^2}\right)^{1/2}
      \int_0^t
       \langle \bu, \bcM_s + s \tBbeta \rangle \, \tbC^{1/2}
       \, \dd \tbcW_s
     \end{bmatrix} , \qquad t\in\RR_+ .
 \]
By the method of the proof of \ $\cX^{(n)} \distr \cX$ \ in Theorem 3.1 in
 Barczy et al.\ \cite{BarIspPap0}, applying Lemma \ref{Conv2Funct}, one can
 easily derive
 \begin{align}\label{convXZ}
  \begin{bmatrix} \bcX^{(n)} \\ \bcZ^{(n)} \end{bmatrix}
  \distr \begin{bmatrix} \tbcX \\ \bcZ \end{bmatrix} , \qquad
  \text{as \ $n \to \infty$,}
 \end{align}
 where
 \[
   \bcX^{(n)}_t = n^{-1} \bX_\nt , \qquad
   \tbcX_t: = \langle \bu, \bcM_t + t \tBbeta \rangle \tbu ,
   \qquad t \in \RR_+ , \qquad n\in \NN .
 \]
Now, with the process
 \begin{equation}\label{tcY}
  \tcY_t := \langle \bu, \tbcX_t \rangle
          = \langle \bu, \bcM_t + t \tBbeta \rangle ,
  \qquad t \in \RR_+ ,
 \end{equation}
 we have
 \[
   \tbcX_t = \tcY_t \tbu , \qquad t \in \RR_+ ,
 \]
 since \ $\langle \bu, \tbu \rangle = 1$.
\ By It\^o's formula and the first 2-dimensional equation of the SDE
 \eqref{MNPSDE} we obtain
 \[
   \dd \tcY_t
   = \langle \bu, \tBbeta \rangle \, \dd t
     + (\tcY_t^+)^{1/2} \bu^\top \tbC^{1/2} \dd \bcW_t , \qquad
   t \in \RR_+ .
 \]
If
 \ $\langle \tbC \bu, \bu \rangle
    = \|\bu^\top \tbC^{1/2}\|^2 = 0$
 \ then \ $\bu^\top \tbC^{1/2} = \bzero$,
 \ hence \ $\dd \tcY_t = \langle \bu, \tBbeta \rangle \dd t$,
 \ $t \in \RR_+$,
 \ implying that the process \ $(\tcY_t)_{t\in\RR_+}$ \ satisfies the SDE
 \eqref{SDE_Y}.
If \ $\langle \tbC \bu, \bu \rangle \ne 0$ \ then the process
 \[
   \ttcW_t := \frac{\langle \tbC^{1/2} \bu, \bcW_t \rangle}
                   {\langle \tbC \bu, \bu \rangle^{1/2}},
   \qquad t \in \RR_+ ,
 \]
 is a (one-dimensional) standard Wiener process, hence the process
 \ $(\tcY_t)_{t\in\RR_+}$ \ satisfies the SDE \eqref{SDE_Y}.
Consequently, \ $\tcY = \cZ$ \ (due to pathwise uniqueness), and hence
 \ $\tbcX = \bcX$.
\ Next, similarly to the proof of \eqref{seged2}, by Lemma \ref{Marci},
 convergence \eqref{convXZ} with
 \ $U_{k-1} = \langle \bu, \bX_{k-1} \rangle$ \ and Lemma \ref{main_VV}
 imply
 \[
   \sum_{k=1}^n
    \begin{bmatrix}
     n^{-2} U_{k-1} \\
     n^{-3} U_{k-1}^2 \\
     n^{-2} V_{k-1}^2 \\
     n^{-1} \bM_k \\
     n^{-2} \langle \bu, \bM_k \rangle U_{k-1} \\
     n^{-3/2} \langle \bv, \bM_k \rangle V_{k-1}
    \end{bmatrix}
   \distr \begin{bmatrix}
           \int_0^1 \langle \bu, \bcX_t \rangle \, \dd t \\[1mm]
           \int_0^1 \langle \bu, \bcX_t \rangle^2 \, \dd t \\[1mm]
           \frac{\langle \tbC \bv, \bv \rangle}
                {1 - \delta^2}
           \int_0^1 \langle \bu, \bcX_t \rangle \, \dd t \\
           \bcM_1 \\
           \int_0^1 \cY_t \, \dd \langle \bu, \bcM_t \rangle \\[1mm]
           \bigl(\frac{\langle \tbC \bv, \bv \rangle}
                      {1 - \delta^2}\bigr)^{1/2}
           \int_0^1 \cY_t \, \dd \langle \bv, \tbC^{1/2} \tbcW_t \rangle
          \end{bmatrix} ,
 \]
 as \ $n \to \infty$.
\ This limiting random vector can be written in the form as given in Theorem
 \ref{main_Ad}, since \ $\langle \bu, \bcX_t \rangle = \cY_t$,
 \ $\langle \bu, \bcM_t \rangle
    = \langle \bu, \bcX_t \rangle
      - \langle \bu, \tBbeta \rangle \, t
    = \cY_t - \langle \bu, \tBbeta \rangle \, t$
 \ (using \eqref{tcY}) and
 \ $\langle \bv, \tbC^{1/2} \tbcW_t \rangle
    = \langle \tbC \bv, \bv \rangle^{1/2} \, \tcW_t$
 \ for all \ $t \in \RR_+$ \ with a (one-dimensional) standard Wiener process
 \ $(\tcW_t)_{t\in\RR_+}$.

\noindent
\textbf{Proof of Theorem \ref{main_conv}.}
In order to show convergence \ $\bcZ^{(n)} \distr \bcZ$, \ we apply Theorem
 \ref{Conv2DiffThm} with the special choices \ $\bcU := \bcZ$,
 \ $\bU^{(n)}_k := \bZ^{(n)}_k$, \ $n, k \in \NN$,
 \ $(\cF_k^{(n)})_{k\in\ZZ_+} := (\cF_k)_{k\in\ZZ_+}$ \ and the function \ $\gamma$
 \ which is defined in Theorem \ref{main_conv}.
Note that the discussion after Theorem \ref{main_conv} shows that the SDE
 \eqref{ZSDE} admits a pathwise unique strong solution
 \ $(\bcZ_t^\bz)_{t\in\RR_+}$ \ for all initial values
 \ $\bcZ_0^\bz = \bz \in (\RR^2)^3$.
\ Applying Cauchy--Schwarz inequality and Corollary \ref{EEX_EEU_EEV}, one can
 check that \ $\EE(\|\bU^{(n)}_k\|^2) < \infty$ \ for all \ $n, k \in \NN$.

Now we show that conditions (i) and (ii) of Theorem \ref{Conv2DiffThm} hold.
The conditional variance takes the form
 \[
   \var\bigl(\bZ^{(n)}_k \mid \cF_{k-1}\bigr)
   =\begin{bmatrix}
    n^{-2}
    & n^{-3} U_{k-1}
    & n^{-5/2} V_{k-1} \\
    n^{-3} U_{k-1}
    & n^{-4} U_{k-1}^2
    & n^{-7/2} U_{k-1} V_{k-1} \\
    n^{-5/2} V_{k-1}
    & n^{-7/2} U_{k-1} V_{k-1}
    & n^{-3} V_{k-1}^2
   \end{bmatrix}
   \otimes \bV_{\!\!\bM_k}
 \]
 for \ $n \in \NN$, \ $k \in \{1, \ldots, n\}$, \ with
 \ $\bV_{\!\!\bM_k} := \var(\bM_k \mid \cF_{k-1})$.
\ Moreover,
 \[
   \gamma(s,\bcZ_s^{(n)}) \gamma(s,\bcZ_s^{(n)})^\top
   =\begin{bmatrix}
    \langle \bu, \bcM_s^{(n)} + s \tBbeta \rangle
     & \langle \bu, \bcM_s^{(n)} + s \tBbeta \rangle^2 & \bzero \\
    \langle \bu, \bcM_s^{(n)} + s \tBbeta \rangle^2
     & \langle \bu, \bcM_s^{(n)} + s \tBbeta \rangle^3 & \bzero \\
    \bzero & \bzero
     & \frac{\langle \tbC \bv, \bv \rangle}
            {1 - \delta^2}
       \langle \bu, \bcM_s^{(n)} + s \tBbeta \rangle^2
   \end{bmatrix}
   \otimes \tbC
 \]
 for \ $s\in\RR_+$, \ where we used that
 \ $\langle \bu, \bcM_s^{(n)} + s \tBbeta \rangle^+
    = \langle \bu, \bcM_s^{(n)} + s \tBbeta \rangle$, \ $s \in \RR_+$,
 \ $n \in \NN$.
\ Indeed, by \eqref{Mk}, we get
 \begin{align}\label{M+}
  \begin{aligned}
   \langle \bu, \bcM_s^{(n)} + s \tBbeta \rangle
   &= \frac{1}{n}
      \sum_{k=1}^\ns
       \langle \bu, \bX_k - \ee^\tbB \bX_{k-1} - \tBbeta \rangle
      + \langle \bu, s \tBbeta \rangle \\
   &= \frac{1}{n}
      \langle \bu, \bX_\ns \rangle
      + \frac{ns - \ns}{n} \langle \bu, \tBbeta \rangle
    = \frac{1}{n} U_\ns + \frac{ns - \ns}{n} \langle \bu, \tBbeta \rangle
      \in \RR_+
  \end{aligned}
 \end{align}
 for \ $s \in \RR_+$, \ $n \in \NN$, \ since
 \ $\bu^\top \ee^\tbB = \bu^\top$ \ implies
 \ $\langle \bu, \ee^\tbB \bX_{k-1} \rangle
    = \bu^\top \ee^\tbB \bX_{k-1} = \bu^\top \bX_{k-1}
    = \langle \bu, \bX_{k-1} \rangle$.

In order to check condition (i) of Theorem \ref{Conv2DiffThm}, we need  to
 prove that for each \ $T > 0$, \ as \ $n \to \infty$,
 \begin{gather}
  \sup_{t\in[0,T]}
   \bigg\| \frac{1}{n^2} \sum_{k=1}^{\nt} \bV_{\!\!\bM_k}
           - \int_0^t
              \langle \bu, \bcM_s^{(n)} + s \tBbeta \rangle \, \tbC \,
              \dd s \bigg\|
  \stoch 0 , \label{Zcond1} \\
  \sup_{t\in[0,T]}
   \bigg\| \frac{1}{n^3}
           \sum_{k=1}^{\nt} U_{k-1} \bV_{\!\!\bM_k}
           - \int_0^t
              \langle \bu, \bcM_s^{(n)} + s \tBbeta \rangle^2 \, \tbC \,
              \dd s \bigg\|
  \stoch 0 , \label{Zcond2} \\
  \sup_{t\in[0,T]}
   \bigg\| \frac{1}{n^4}
           \sum_{k=1}^{\nt} U_{k-1}^2 \bV_{\!\!\bM_k}
           - \int_0^t
              \langle \bu, \bcM_s^{(n)} + s \tBbeta \rangle^3 \, \tbC \,
              \dd s \bigg\|
  \stoch 0 , \label{Zcond3} \\
  \sup_{t\in[0,T]}
   \bigg\| \frac{1}{n^3}
           \sum_{k=1}^{\nt} V_{k-1}^2 \bV_{\!\!\bM_k}
           - \frac{\langle \tbC \bv, \bv \rangle}
                  {1 - \delta^2}
             \int_0^t
              \langle \bu, \bcM_s^{(n)} + s \tBbeta \rangle^2 \, \tbC \,
              \dd s \bigg\|
  \stoch 0 , \label{Zcond4} \\
  \sup_{t\in[0,T]}
   \bigg\| \frac{1}{n^{5/2}}
           \sum_{k=1}^{\nt} V_{k-1} \bV_{\!\!\bM_k} \bigg\|
  \stoch 0 , \label{Zcond5} \\
  \sup_{t\in[0,T]}
   \bigg\| \frac{1}{n^{7/2}}
           \sum_{k=1}^{\nt}
            U_{k-1} V_{k-1} \bV_{\!\!\bM_k} \bigg\|
  \stoch 0 . \label{Zcond6}
 \end{gather}

First we show \eqref{Zcond1}.
By \eqref{M+},
 \ $\int_0^t \langle \bu, \bcM_s^{(n)} + s \tBbeta \rangle \, \dd s$
 \ has the form
 \[
   \frac{1}{n^2} \sum_{k=1}^{\nt-1} U_k
   + \frac{nt - \nt}{n^2} U_\nt
   + \frac{\nt + (nt - \nt)^2}{2 n^2} \langle \bu, \tBbeta \rangle .
 \]
Using Proposition \ref{moment_formula_2}, formula \eqref{XUV} and
 \ $\tbC = (\bV_1 + \bV_2) / 2$, \ we obtain
 \begin{equation}\label{VMk}
  \begin{aligned}
   \bV_{\!\!\bM_k}
   &= \var(\bM_{k} \mid \cF_{k-1})
    = \frac{1}{2} U_{k-1} (\bV_1 + \bV_2) + \frac{1}{2} V_{k-1} (\bV_1 - \bV_2)
      + \bV_0 \\
   &= U_{k-1} \tbC + \frac{1}{2} V_{k-1} (\bV_1 - \bV_2) + \bV_0 .
  \end{aligned}
 \end{equation}
Thus, in order to show \eqref{Zcond1}, it suffices to prove
 \begin{gather}
  n^{-2} \sum_{k=1}^{\nT} |V_k| \stoch 0 , \qquad
  n^{-2} \sup_{t \in [0,T]} U_\nt \stoch 0, \label{1supsumV_1supU} \\
  n^{-2} \sup_{t \in [0,T]} \left[ \nt + (nt - \nt)^2 \right] \to 0 ,
  \label{1supnt}
 \end{gather}
 as \ $n \to \infty$.
\ Using \eqref{seged_UV_UNIFORM1} with \ $(\ell, i, j) = (2, 0, 1)$ \ and
 \eqref{seged_UV_UNIFORM2} with \ $(\ell, i, j) = (2, 1, 0)$, \ we have
 \eqref{1supsumV_1supU}.
Clearly, \eqref{1supnt} follows from \ $|nt - \nt| \leq 1$, \ $n \in \NN$,
 \ $t \in \RR_+$, \ thus we conclude \eqref{Zcond1}.

Next we turn to prove \eqref{Zcond2}.
By \eqref{M+},
 \begin{align*}
  \int_0^t \langle \bu, \bcM_s^{(n)} + s \tBbeta \rangle^2 \, \dd s
  &= \frac{1}{n^3} \sum_{k=1}^{\nt-1} U_k^2
     + \frac{1}{n^3} \langle \bu, \tBbeta \rangle \sum_{k=1}^{\nt-1} U_k
     + \frac{nt - \nt}{n^3} U_\nt^2 \\
  &\phantom{\quad}
     + \frac{(nt - \nt)^2}{n^3} \langle \bu, \tBbeta \rangle U_\nt
     + \frac{\nt + (nt - \nt)^3}{3n^3} \langle \bu, \tBbeta \rangle^2 .
 \end{align*}
Recalling formula \eqref{VMk}, we obtain
 \begin{equation}\label{UM2F}
  \sum_{k=1}^{\nt} U_{k-1} \bV_{\!\!\bM_k}
  = \sum_{k=1}^{\nt} U_{k-1}^2 \tbC
     + \frac{1}{2} \sum_{k=1}^{\nt} U_{k-1} V_{k-1} (\bV_1 - \bV_2)
     + \sum_{k=1}^{\nt} U_{k-1} \bV_0 .
 \end{equation}
Thus, in order to show \eqref{Zcond2}, it suffices to prove
 \begin{gather}
  n^{-3} \sum_{k=1}^{\nT} |U_k V_k| \stoch 0 , \label{2supsumUV} \\
  n^{-3} \sum_{k=1}^{\nT} U_k \stoch 0 , \label{2supsumU} \\
  n^{-3/2} \sup_{t \in [0,T]} U_\nt \stoch 0, \label{2supU} \\
  n^{-3} \sup_{t \in [0,T]} \left[ \nt + (nt - \nt)^3 \right] \to 0
    \label{2supnt}
 \end{gather}
 as \ $n \to \infty$.
\ Using \eqref{seged_UV_UNIFORM1} with \ $(\ell, i, j) = (2, 1, 1)$ \ and
 \ $(\ell, i, j) = (2, 1, 0)$, \ we have \eqref{2supsumUV} and
 \eqref{2supsumU}, respectively.
By \eqref{seged_UV_UNIFORM2} with \ $(\ell, i, j) = (3, 1, 0)$, \ we have
 \eqref{2supU}.
Clearly, \eqref{2supnt} follows from \ $|nt - \nt| \leq 1$, \ $n \in \NN$,
 \ $t \in \RR_+$, \ thus we conclude \eqref{Zcond2}.

Now we turn to check \eqref{Zcond3}.
Again by \eqref{M+}, we have
 \begin{align*}
  \int_0^t \langle \bu, \bcM_s^{(n)} + s \tBbeta \rangle^3 \, \dd s
  &= \frac{1}{n^4} \sum_{k=1}^{\nt-1} U_k^3
     + \frac{3}{2n^4} \langle \bu, \tBbeta \rangle \sum_{k=1}^{\nt-1} U_k^2
     + \frac{1}{n^4} \langle \bu, \tBbeta \rangle^2
       \sum_{k=1}^{\nt-1} U_k \\
  &\phantom{\quad}
     + \frac{nt - \nt}{n^4} U_\nt^3
     + \frac{3 (nt - \nt)^2}{2n^4} \langle \bu, \tBbeta \rangle U_\nt^2 \\
  &\phantom{\quad}
     + \frac{(nt - \nt)^3}{n^4} \langle \bu, \tBbeta \rangle^2 \, U_\nt
     + \frac{\nt + (nt - \nt)^4}{4n^4} \langle \bu, \tBbeta \rangle^3 .
 \end{align*}
Recalling formula \eqref{VMk}, we obtain
 \begin{equation}\label{U2M2F}
  \sum_{k=1}^{\nt} U_{k-1}^2 \bV_{\!\!\bM_k}
  = \sum_{k=1}^{\nt} U_{k-1}^3 \tbC
    + \frac{1}{2} \sum_{k=1}^{\nt} U_{k-1}^2 V_{k-1} (\bV_1 - \bV_2)
    + \sum_{k=1}^{\nt} U_{k-1}^2 \bV_0 .
 \end{equation}
Thus, in order to show \eqref{Zcond3}, it suffices to prove
 \begin{gather}
  n^{-4} \sum_{k=1}^{\nT} | U_k^2 V_k | \stoch 0 , \label{3supsumUUV} \\
  n^{-4} \sum_{k=1}^{\nT} U_k^2 \stoch 0 , \label{3supsumUU} \\
  n^{-4} \sum_{k=1}^{\nT} U_k \stoch 0 , \label{3supsumU} \\
  n^{-4/3} \sup_{t \in [0,T]} U_\nt \stoch 0 , \label{3supU} \\
  n^{-4} \sup_{t \in [0,T]} \left[ \nt + (nt - \nt)^4 \right] \to 0
    \label{3supnt}
 \end{gather}
 as \ $n \to \infty$.
\ Using \eqref{seged_UV_UNIFORM1} with \ $(\ell, i, j) = (4, 2, 1)$,
 \ $(\ell, i, j) = (4, 2, 0)$, \ and \ $(\ell, i, j) = (2, 1, 0)$, \ we have
 \eqref{3supsumUUV}, \eqref{3supsumUU} and \eqref{3supsumU}, respectively.
By \eqref{seged_UV_UNIFORM2} with \ $(\ell, i, j) = (4, 1, 0)$, \ we have
 \eqref{3supU}.
Clearly, \eqref{3supnt} follows again from \ $|nt - \nt| \leq 1$,
 \ $n \in \NN$, \ $t \in \RR_+$, \ thus we conclude \eqref{Zcond3}.
Note that the proof of \eqref{Zcond1}--\eqref{Zcond3} is essentially the same
 as the proof of (5.5)--(5.7) in Isp\'any et al.\ \cite{IspKorPap}.

Next we turn to prove \eqref{Zcond4}.
First we show that
 \begin{align}\label{Zcond4a}
  n^{-3}
  \sup_{t \in [0,T]}
   \left\| \sum_{k=1}^{\nt} V_{k-1}^2 \bV_{\!\!\bM_k}
           - \frac{\langle \tbC \bv, \bv \rangle}
                  {1 - \delta^2}
             \sum_{k=1}^{\nt} U_{k-1}^2 \tbC \right\|
  \stoch 0 ,
 \end{align}
 as \ $n \to \infty$ \ for all \ $T > 0$.
\ By \eqref{VMk},
 \[
   \sum_{k=1}^{\nt} V_{k-1}^2 \bV_{\!\!\bM_k}
   = \sum_{k=1}^{\nt} U_{k-1} V_{k-1}^2 \tbC
     + \frac{1}{2} \sum_{k=1}^{\nt} V_{k-1}^3 (\bV_1 - \bV_2)
     + \sum_{k=1}^{\nt} V_{k-1}^2 \bV_0 .
 \]
Using \eqref{seged_UV_UNIFORM1} with \ $(\ell, i, j) = (6, 0, 3)$ \ and
 \ $(\ell, i, j) = (4, 0, 2)$, \ we have
 \begin{align*}
  n^{-3} \sum_{k=1}^{\nT} | V_k |^3 \stoch 0 , \qquad
  n^{-3} \sum_{k=1}^{\nT} V_k^2 \stoch 0 , \qquad \text{as \ $n\to\infty$,}
 \end{align*}
 hence \eqref{Zcond4a} will follow from
 \begin{align}\label{Zcond4b}
  n^{-3}
  \sup_{t \in [0,T]}
   \left| \sum_{k=1}^{\nt} U_{k-1} V_{k-1}^2
          - \frac{\langle \tbC \bv, \bv \rangle}
                 {1 - \delta^2}
            \sum_{k=1}^{\nt} U_{k-1}^2 \right|
  \stoch 0 ,
 \end{align}
 as \ $n \to \infty$ \ for all \ $T>0$.
\ By the method of the proof of Lemma \ref{main_VV}
 (see also Isp\'any et al.\
 \cite[page 16 of arXiv version]{IspKorPap}), applying
 Proposition
 \ref{moment_formula_3} with \ $q = 3$, \ we obtain a decomposition of
 \ $\sum_{k=1}^{\nt} U_{k-1} V_{k-1}^2$, \ namely,
 \begin{align*}
  \sum_{k=1}^{\nt} U_{k-1} V_{k-1}^2
  &= \frac{1}{1-\delta^2}
     \sum_{k=2}^{\nt}
      \big[U_{k-1} V_{k-1}^2 - \EE(U_{k-1} V_{k-1}^2 \mid \cF_{k-2}) \big] \\
  &\quad
   + \frac{\langle \tbC \bv, \bv \rangle}
          {1 - \delta^2}
       \sum_{k=2}^{\nt} U_{k-2}^2
     - \frac{\delta^2}
            {1 - \delta^2}
       U_{\nt - 1} V_{\nt - 1}^2
     + \OO(n) \\
  &\quad
   + \text{lin.~comb.~of \ $\sum_{k=2}^{\nt} U_{k-2} V_{k-2}$,
             \ $\sum_{k=2}^{\nt} V_{k-2}^2$, \ $\sum_{k=2}^{\nt} U_{k-2}$
             \ and \ $\sum_{k=2}^{\nt} V_{k-2}$.}
 \end{align*}
Note that Proposition \ref{moment_formula_3} with \ $q = 3$ \ is needed
 above in order to express products
 \ $\EE(M_{k-1,i_1} M_{k-1,i_2} M_{k-1,i_2} \mid \cF_{k-2})$,
 \ $i_1, i_2, i_3 \in \{1, 2\}$, \ as a first order polynomial of \ $\bX_{k-2}$,
 \ and hence, by \eqref{XUV}, as a linear combination of \ $U_{k-2}$,
 \ $V_{k-2}$ \ and 1.
Using \eqref{seged_UV_UNIFORM4} with \ $(\ell, i, j) = (8, 1, 2)$ \ we have
 \[
   n^{-3}\sup_{t \in [0,T]}\,
         \Biggl\vert \sum_{k=2}^\nt
                      \big[U_{k-1} V_{k-1}^2
                           - \EE(U_{k-1} V_{k-1}^2 \mid \cF_{k-2}) \big]
         \Biggr\vert
   \stoch 0 , \qquad \text{as \ $n\to\infty$.}
 \]
In order to show \eqref{Zcond4b}, it suffices to prove
 \begin{gather}
  n^{-3} \sum_{k=1}^{\nT} | U_k V_k | \stoch 0 , \qquad
  n^{-3} \sum_{k=1}^{\nT} V_k^2 \stoch 0 , \label{4supsumUV_4supsumVV} \\
  n^{-3} \sum_{k=1}^{\nT} U_k \stoch 0 , \qquad
  n^{-3} \sum_{k=1}^{\nT} |V_k| \stoch 0 , \label{4supsumU_4supsumV} \\
  n^{-3} \sup_{t \in [0,T]} U_\nt V_\nt^2 \stoch 0 , \qquad
  n^{-3/2} \sup_{t \in [0,T]} U_{\nt} \stoch 0 , \label{4supUVV_4supUU}
   \end{gather}
 as \ $n \to \infty$.
\ Using \eqref{seged_UV_UNIFORM1} with \ $(\ell, i, j) = (2, 1, 1)$,
 \ $(\ell, i, j) = (4, 0, 2)$, \ $(\ell, i, j) = (2, 1, 0)$ \ and
 \ $(\ell, i, j) = (2, 0, 1)$, \ we have \eqref{4supsumUV_4supsumVV} and
 \eqref{4supsumU_4supsumV}.
By \eqref{seged_UV_UNIFORM2} with \ $(\ell, i, j) = (4, 1, 2)$ \ and
 \ $(\ell, i, j) = (3, 1, 0)$, \ we have \eqref{4supUVV_4supUU}.
Thus we conclude \eqref{Zcond4a}.
By \eqref{VMk} and \eqref{seged_UV_UNIFORM1} with \ $(\ell, i, j) = (2, 1, 1)$
 \ and \ $(\ell, i, j) = (2, 1, 0)$, \ we get
 \begin{align}\label{Zcond2a}
  n^{-3}
  \sup_{t \in [0,T]}
   \left\| \sum_{k=1}^{\nt} U_{k-1} \bV_{\!\!\bM_k}
           - \sum_{k=1}^{\nt} U_{k-1}^2 \tbC \right\|
  \stoch 0 ,
 \end{align}
 as \ $n \to \infty$ \ for all \ $T>0$.
\ As a last step, using \eqref{Zcond2}, we obtain \eqref{Zcond4}.
Convergences \eqref{Zcond5} and \eqref{Zcond6} can be proved similarly
 (see also the same considerations in
 Isp\'any et al.\ \cite[pages 17-20 of arXiv version]{IspKorPap}).

Finally, we check condition (ii) of Theorem \ref{Conv2DiffThm}, that is, the
 conditional Lindeberg condition
 \begin{equation}\label{Zcond3_new}
   \sum_{k=1}^{\lfloor nT \rfloor}
     \EE \big( \|\bZ^{(n)}_k\|^2 \bbone_{\{\|\bZ^{(n)}_k\| > \theta\}}
               \bmid \cF_{k-1} \big)
    \stoch 0 , \qquad \text{as \ $n\to\infty$}
 \end{equation}
 for all \ $\theta>0$ \ and \ $T>0$.
\ We have
 \ $\EE \big( \|\bZ^{(n)}_k\|^2 \bbone_{\{\|\bZ^{(n)}_k\| > \theta\}}
              \bmid \cF_{k-1} \big)
    \leq \theta^{-2} \EE \big( \|\bZ^{(n)}_k\|^4 \bmid \cF_{k-1} \big)$
 \ and
 \[
   \|\bZ^{(n)}_k\|^4
   \leq 3 \left( n^{-4} + n^{-8} U_{k-1}^4 + n^{-6} V_{k-1}^4 \right)
        \|\bM_k\|^4 .
 \]
 Hence, for all \ $\theta>0$ \ and \ $T>0$, \ we have
 \[
   \sum_{k=1}^{\nT}
    \EE \big( \|\bZ^{(n)}_k\|^2 \bbone_{\{\|\bZ^{(n)}_k\| > \theta\}} \big)
   \to 0 ,
   \qquad \text{as \ $n\to\infty$,}
 \]
 since \ $\EE( \|\bM_k\|^4 ) = \OO(k^2)$,
 \ $\EE( \|\bM_k\|^4 U_{k-1}^4 ) \leq \sqrt{\EE(\|\bM_k\|^8) \EE(U_{k-1}^8)}
    = \OO(k^6)$
 \ and
 \ $\EE( \|\bM_k\|^4 V_{k-1}^4 ) \leq \sqrt{\EE(\|\bM_k\|^8) \EE(V_{k-1}^8)}
    = \OO(k^4)$
 \ by Corollary \ref{EEX_EEU_EEV}.
This yields \eqref{Zcond3_new}.
\proofend

We call the attention that our moment conditions \eqref{moment_condition_m}
 with \ $q = 8$ \ are used for applying Corollaries \ref{EEX_EEU_EEV} and
 \ref{LEM_UV_UNIFORM}.

\section{Proof of Theorem \ref{maint_Ad}}
\label{section_proof_maint}

The first and second convergences follow from Lemmas \ref{main_Vt} and
 \ref{main_VVt}.
The proof of the third convergence in Theorem \ref{maint_Ad} is similar to
 the proof of Theorem \ref{main_Ad}.
Consider the sequence of stochastic processes
 \[
   \bcZ^{(n)}_t
   := \begin{bmatrix}
       \cM_t^{(n)} \\
       \cN_t^{(n)} \\
       \cP_t^{(n)} \\
       \cR_t^{(n)}
      \end{bmatrix}
   := \sum_{k=1}^\nt
       \bZ^{(n)}_k
   \qquad \text{with} \qquad
   \bZ^{(n)}_k
   := \begin{bmatrix}
       n^{-1} \langle \bu, \bM_k \rangle \\
       n^{-2} \langle \bu, \bM_k \rangle U_{k-1} \\
       n^{-1/2} \langle \bv, \bM_k \rangle \\
       n^{-1/2} \langle \bv, \bM_k \rangle V_{k-1}
      \end{bmatrix}
 \]
 for \ $t \in \RR_+$ \ and \ $k, n \in \NN$.

\begin{Thm}\label{maint_conv}
Suppose that the assumptions of Theorem \ref{main} hold.
If \ $\langle \tbC \bv, \bv \rangle = 0$, \ then
 \begin{equation}\label{conv_Zt}
  \bcZ^{(n)} \distr \bcZ \qquad \text{as \ $n\to\infty$,}
 \end{equation}
 where the process \ $(\bcZ_t)_{t \in \RR_+}$ \ with values in
 \ $\RR^4$ \ is the pathwise unique strong solution of the SDE
 \begin{equation}\label{ZSDEt}
  \dd \bcZ_t
  = \gamma(t, \bcZ_t) \begin{bmatrix} \dd \cW_t \\ \dd \tbcW_t \end{bmatrix} ,
  \qquad t \in \RR_+ ,
 \end{equation}
 with initial value \ $\bcZ_0 = \bzero$, \ where \ $(\cW_t)_{t \in \RR_+}$ \ and
 \ $(\tbcW_t)_{t \in \RR_+}$ \ are independent standard Wiener processes of
 dimension 1 and 2, respectively, and
 \ $\gamma : \RR_+ \times \RR^4 \to \RR^{4\times3}$
 \ is defined by
 \[
   \gamma(t, \bx)
   := \begin{bmatrix}
       ((x_1 + \langle \bu, \tBbeta \rangle t)^+)^{1/2}
       \langle \tbC \bu, \bu \rangle^{1/2}
        & \bzero^\top \\
       ((x_1 + \langle \bu, \tBbeta \rangle t)^+)^{3/2}
       \langle \tbC \bu, \bu \rangle^{1/2}
        & \bzero^\top \\
       \bzero
        & \langle \bV_0 \bv, \bv \rangle^{1/2}
          \begin{bmatrix}
           1 & \frac{\tdelta\langle \bv, \tBbeta \rangle}{1-\delta} \\
           \frac{\tdelta\langle \bv, \tBbeta \rangle}{1-\delta} & M
          \end{bmatrix}^{1/2}
      \end{bmatrix}
 \]
 for \ $t \in \RR_+$ \ and \ $\bx = (x_1, x_2, x_3, x_4)^\top \in \RR^4$.
\end{Thm}

The SDE \eqref{ZSDEt} has the form
 \begin{align}\label{MNPRSDEt}
  \dd \bcZ_t
  =: \begin{bmatrix}
      \dd \cM_t \\
      \dd \cN_t \\
      \dd \cP_t \\
      \dd \cR_t
     \end{bmatrix}
  = \begin{bmatrix}
     ((\cM_t + \langle \bu, \tBbeta \rangle t)^+)^{1/2}
     \langle \tbC \bu, \bu \rangle^{1/2} \,
     \dd \cW_t \\[1mm]
     ((\cM_t + \langle \bu, \tBbeta \rangle t)^+)^{3/2}
     \langle \tbC \bu, \bu \rangle^{1/2} \,
     \dd \cW_t \\[1mm]
     \langle \bV_0 \bv, \bv \rangle^{1/2}
     \begin{bmatrix}
      1 & \frac{\tdelta\langle \bv, \tBbeta \rangle}{1-\delta} \\
      \frac{\tdelta\langle \bv, \tBbeta \rangle}{1-\delta} & M
     \end{bmatrix}^{1/2}
     \dd \tbcW_t
    \end{bmatrix} , \qquad t \in \RR_+ .
 \end{align}
One can prove that the first equation of the SDE
 \eqref{MNPRSDEt} has a pathwise unique strong solution
 \ $(\cM_t^{(y_0)})_{t\in\RR_+}$ \ with arbitrary initial value
 \ $\cM_0^{(y_0)} = y_0 \in \RR$.
\ Indeed, it is equivalent to the existence of a pathwise unique strong
 solution of the SDE
 \[
   \dd \cS_t
   = \langle \bu, \tBbeta \rangle \, \dd t
     + (\cS_t^+)^{1/2} \, \langle \tbC \bu, \bu \rangle^{1/2} \,
       \dd \cW_t ,
 \]
 with initial value
 \ $\cS_0^{(y_0)} = y_0 \in \RR$, \ since we have the correspondence
 \ $\cS_t^{(y_0)} = \cM_t^{(y_0)} + \langle \bu, \tBbeta \rangle \, t$,
 \ $t \in \RR_+$.
\ Thus the SDE \eqref{ZSDEt} has a pathwise unique strong solution with initial
 value \ $\bcZ_0 = \bzero$, \ and we have
 \[
   \bcZ_t
   = \begin{bmatrix}
      \cM_t \\
      \cN_t \\
      \cP_t \\
      \cR_t
     \end{bmatrix}
   = \begin{bmatrix}
      \langle \tbC \bu, \bu \rangle^{1/2}
      \int_0^t \tcY_s^{1/2} \, \dd \cW_s \\[1mm]
      \int_0^t \tcY_s \, \dd \cM_s \\[1mm]
      \langle \bV_0 \bv, \bv \rangle^{1/2}
      \begin{bmatrix}
       1 & \frac{\tdelta\langle \bv, \tBbeta \rangle}{1-\delta} \\
       \frac{\tdelta\langle \bv, \tBbeta \rangle}{1-\delta} & M
      \end{bmatrix}^{1/2}
      \tbcW_t
     \end{bmatrix} , \qquad t\in\RR_+ ,
 \]
 where, by It\^o's formula, the process
 \ $\tcY_t := \cM_t + \langle \bu, \tBbeta \rangle t$,
 \ $t \in \RR_+$, \ satisfies the SDE \eqref{SDE_Y}.
\ By Remark \ref{REMARK0.7},
 \ $(\tcY_t)_{t\in\RR_+} \distre (\cY_t)_{t\in\RR_+}$, \ hence
 \ $(\cM_t)_{t\in\RR_+}
    \distre (\cY_t - \langle \bu, \tBbeta \rangle t)_{t\in\RR_+}
    = (\langle \bu, \bcM_t \rangle)_{t\in\RR_+}$,
 \ thus we obtain convergence of the last four coordinates of the third
 convergence in Theorem \ref{maint_Ad}.
One can again easily derive
 \begin{align}\label{convXZt}
  \begin{bmatrix}
   \langle \bu, \bcX^{(n)} \rangle \\
   \bcZ^{(n)}
  \end{bmatrix}
  \distr
  \begin{bmatrix}
   \langle \bu, \bcX \rangle \\
   \bcZ
  \end{bmatrix} \qquad
  \text{as \ $n \to \infty$,}
 \end{align}
 where
 \[
   \langle \bu, \bcX^{(n)}_t \rangle
   = n^{-1} \langle \bu, \bX_\nt \rangle , \qquad
   \langle \bu, \bcX_t \rangle
   = \langle \bu, \cY_t \tbu \rangle
   = \cY_t ,
 \]
 for \ $t \in \RR_+$, \ $n \in \NN$.
\ Next, similarly to the proof of \eqref{seged2}, by Lemma \ref{Marci},
 convergence \eqref{convXZt} and Lemma \ref{main_VVt} imply the third
 convergence in Theorem \ref{maint_Ad}.

\noindent
\textbf{Proof of Theorem \ref{maint_conv}.}
It is similar to the proof of Theorem \ref{main_conv}.
The conditional variance \ $\var\bigl(\bZ^{(n)}_k \mid \cF_{k-1}\bigr)$
 \ has the form
 \[
   \begin{bmatrix}
    \bu^\top \bV_{\!\!\bM_k} \bu
    \begin{bmatrix}
     n^{-2} & n^{-3} U_{k-1} \\
     n^{-3} U_{k-1} & n^{-4} U_{k-1}^2
    \end{bmatrix}
     & \bu^\top \bV_{\!\!\bM_k} \bv
       \begin{bmatrix}
        n^{-3/2} & n^{-3/2} V_{k-1} \\
        n^{-5/2} U_{k-1} & n^{-5/2} U_{k-1} V_{k-1}
       \end{bmatrix} \\[4mm]
    \bv^\top \bV_{\!\!\bM_k} \bu
       \begin{bmatrix}
        n^{-3/2} & n^{-5/2} U_{k-1} \\
        n^{-3/2} V_{k-1} & n^{-5/2} U_{k-1} V_{k-1}
       \end{bmatrix}
     & \bv^\top \bV_{\!\!\bM_k} \bv
       \begin{bmatrix}
        n^{-1} & n^{-1} V_{k-1} \\
        n^{-1} V_{k-1} & n^{-1} V_{k-1}^2
       \end{bmatrix}
   \end{bmatrix}
 \]
 for \ $n \in \NN$, \ $k \in \{1, \ldots, n\}$, \ with
 \ $\bV_{\!\!\bM_k} = \var( \bM_k \mid \cF_{k-1} )$, \ and
 \ $\gamma(s,\bcZ_s^{(n)}) \gamma(s,\bcZ_s^{(n)})^\top$ \ has the form
 \[
   \begin{bmatrix}
    \bu^\top \tbC \bu
    \begin{bmatrix}
     \cM_s^{(n)} + \langle \bu, \tBbeta \rangle s
      & (\cM_s^{(n)} + \langle \bu, \tBbeta \rangle s)^2 \\
     (\cM_s^{(n)} + \langle \bu, \tBbeta \rangle s)^2
      & (\cM_s^{(n)} + \langle \bu, \tBbeta \rangle s)^3
    \end{bmatrix}
    & \bzero \\
    \bzero
     & \bv^\top \bV_0 \bv
       \begin{bmatrix}
        1 & \frac{\tdelta\langle \bv, \tBbeta \rangle}{1-\delta} \\
        \frac{\tdelta\langle \bv, \tBbeta \rangle}{1-\delta} & M
       \end{bmatrix}
   \end{bmatrix}
 \]
 for \ $s \in \RR_+$, \ where we used that
 \ $(\cM_s^{(n)} + \langle \bu, \tBbeta \rangle s)^+
    = \cM_s^{(n)} + \langle \bu, \tBbeta \rangle s$, \ $s \in \RR_+$,
 \ $n \in \NN$.

In order to check condition (i) of Theorem \ref{Conv2DiffThm},
 we need to prove only that for each \ $T > 0$,
 \begin{gather}
  \sup_{t\in[0,T]}
   \bigg| \frac{1}{n}
           \sum_{k=1}^{\nt}
            \bv^\top \bV_{\!\!\bM_k} \bv
           - t \langle \bV_0 \bv, \bv \rangle \bigg|
  \stoch 0 , \label{tZcond41} \\
  \sup_{t\in[0,T]}
   \bigg| \frac{1}{n}
           \sum_{k=1}^{\nt}
            V_{k-1} \bv^\top \bV_{\!\!\bM_k} \bv
           - t \langle \bV_0 \bv, \bv \rangle
             \frac{\tdelta\langle \bv, \tBbeta \rangle}{1-\delta} \bigg|
  \stoch 0 , \label{tZcond42} \\
  \sup_{t\in[0,T]}
   \bigg| \frac{1}{n}
           \sum_{k=1}^{\nt}
            V_{k-1}^2 \bv^\top \bV_{\!\!\bM_k} \bv
           - t \langle \bV_0 \bv, \bv \rangle M \bigg|
  \stoch 0 , \label{tZcond43} \\
  \sup_{t\in[0,T]}
   \bigg| \frac{1}{n^{3/2}}
           \sum_{k=1}^{\nt} \bu^\top \bV_{\!\!\bM_k} \bv \bigg|
  \stoch 0 , \label{tZcond51} \\
  \sup_{t\in[0,T]}
   \bigg| \frac{1}{n^{3/2}}
           \sum_{k=1}^{\nt} V_{k-1} \bu^\top \bV_{\!\!\bM_k} \bv \bigg|
  \stoch 0 , \label{tZcond52} \\
  \sup_{t\in[0,T]}
   \bigg| \frac{1}{n^{5/2}}
           \sum_{k=1}^{\nt}
            U_{k-1} \bu^\top
            \bV_{\!\!\bM_k} \bv \bigg|
  \stoch 0  \label{tZcond61} \\
  \sup_{t\in[0,T]}
   \bigg| \frac{1}{n^{5/2}}
           \sum_{k=1}^{\nt}
            U_{k-1} V_{k-1} \bu^\top
            \bV_{\!\!\bM_k} \bv \bigg|
  \stoch 0  \label{tZcond62}
 \end{gather}
 as \ $n \to \infty$, \ since the rest follows from \eqref{Zcond1}, \eqref{Zcond2}
 and \eqref{Zcond3}.

First we show \eqref{tZcond41}, \eqref{tZcond42} and \eqref{tZcond43}.
The assumption \ $\langle \tbC \bv, \bv \rangle = 0$ \ yields
 \ $\langle \bV_i \bv, \bv \rangle = 0$, \ $i \in \{1, 2\}$, \ see
 the proof of Lemma \ref{main_VVt}.
Thus, by Proposition \ref{moment_formula_2},
 \begin{gather}
  \sum_{k=1}^{\nt}\bv^\top \bV_{\!\!\bM_k} \bv
  = \sum_{k=1}^{\nt} \langle \bV_0 \bv, \bv \rangle , \\
  \sum_{k=1}^{\nt} V_{k-1} \bv^\top \bV_{\!\!\bM_k} \bv
  = \sum_{k=1}^{\nt} V_{k-1} \langle \bV_0 \bv, \bv \rangle , \\
  \sum_{k=1}^{\nt} V_{k-1}^2 \bv^\top \bV_{\!\!\bM_k} \bv
  = \sum_{k=1}^{\nt} V_{k-1}^2 \langle \bV_0 \bv, \bv \rangle ,
 \end{gather}
 hence \eqref{tZcond41}, \eqref{tZcond42} and \eqref{tZcond43} follow from Lemmas
 \ref{main_Vt} and \ref{main_VVt}.

Next we turn to check \eqref{tZcond51}, \eqref{tZcond52}, \eqref{tZcond61} and
 \eqref{tZcond62}.
For each \ $i \in \{1, 2\}$,
 \ $\langle \bV_i \bv, \bv \rangle = 0$ \ yields
 \ $\bV_i \bv = \bzero$, \ since
 \ $0 = \langle \bV_i \bv, \bv \rangle = \bv^\top \bV_i \bv
      = (\bV_i^{1/2} \bv)^\top (\bV_i^{1/2} \bv)
      = \|\bV_i^{1/2} \bv\|^2$
 \ implies \ $\bV_i^{1/2} \bv = \bzero$, \ and hence
 \ $\bV_i \bv = \bV_i^{1/2} (\bV_i^{1/2} \bv) = \bzero$.
 \
Consequently, by \eqref{VMk},
 \begin{gather}
  \sum_{k=1}^{\nt} \bu^\top \bV_{\!\!\bM_k} \bv
  = \sum_{k=1}^{\nt} \langle \bV_0 \bu, \bv \rangle , \\
  \sum_{k=1}^{\nt} V_{k-1} \bu^\top \bV_{\!\!\bM_k} \bv
  = \sum_{k=1}^{\nt} V_{k-1} \langle \bV_0 \bu, \bv \rangle , \\
  \sum_{k=1}^{\nt} U_{k-1} \bu^\top \bV_{\!\!\bM_k} \bv
  = \sum_{k=1}^{\nt} U_{k-1} \langle \bV_0 \bu, \bv \rangle , \\
  \sum_{k=1}^{\nt} U_{k-1} V_{k-1} \bu^\top \bV_{\!\!\bM_k} \bv
  = \sum_{k=1}^{\nt} U_{k-1} V_{k-1} \langle \bV_0 \bu, \bv \rangle .
 \end{gather}
Thus \eqref{tZcond51} is trivial.
Using \eqref{seged_UV_UNIFORM1_mod} with \ $(\ell, i, j) = (2, 0, 1)$,
 \ $(\ell, i, j) = (3, 1, 0)$ \ and \ $(\ell, i, j) = (3, 1, 1)$, \ we
 conclude \eqref{tZcond52}, \eqref{tZcond61} and \eqref{tZcond62}.

Finally, we check condition (ii) of Theorem \ref{Conv2DiffThm}, that is, the
 conditional Lindeberg condition for all \ $\theta > 0$ \ and \ $T>0$.
\ We have
 \ $\EE \big( \|\bZ^{(n)}_k\|^2 \bbone_{\{\|\bZ^{(n)}_k\| > \theta\}}
              \bmid \cF_{k-1} \big)
    \leq \theta^{-2} \EE \big( \|\bZ^{(n)}_k\|^4 \bmid \cF_{k-1} \big)$
 \ and
 \[
   \|\bZ^{(n)}_k\|^4
   \leq 4 n^{-4} \langle \bu, \bM_k \rangle^4
        + 4 n^{-8} \langle \bu, \bM_k \rangle^4 U_{k-1}^4
        + 4 n^{-2} \langle \bv, \bM_k \rangle^4
        + 4 n^{-2} \langle \bv, \bM_k \rangle^4 V_{k-1}^4 .
 \]
 Hence, for all \ $\theta > 0$ \ and \ $T > 0$, \ we have
 \[
   \sum_{k=1}^{\nT}
    \EE \big( \|\bZ^{(n)}_k\|^2 \bbone_{\{\|\bZ^{(n)}_k\| > \theta\}} \big)
   \to 0 ,
   \qquad \text{as \ $n\to\infty$,}
 \]
 since
 \ $\EE(\langle \bu, \bM_k \rangle^4) \leq 4 \EE(\|\bM_k\|^4) = \OO(k^2)$
 \ and
 \ $\EE(\langle \bu, \bM_k \rangle^4 U_{k-1}^4)
    \leq 4 \EE(\|\bM_k\|^4 U_{k-1}^4)) = \OO(k^6)$,
 \ as in the proof of Theorem \ref{main_conv}, and
 \ $\EE(\langle \bv, \bM_k \rangle^4) = \OO(1)$ \ and
 \ $\EE(\langle \bv, \bM_k \rangle^4 V_{k-1}^4)
    \leq \sqrt{\EE(V_{k-1}^8) \EE(\langle \bv, \bM_k \rangle^8)}
    = \OO(1)$
 \ by Corollary \ref{EEX_EEU_EEV}.
\proofend

\section{Proof of Theorem \ref{main1_Ad}}
\label{section_proof_main1}

The first and second convergences in Theorem \ref{main1_Ad} follow from the
 following two approximations.

\begin{Lem}\label{main1_U}
Suppose that the assumptions of Theorem \ref{main} hold.
If \ $\langle \tbC \bu, \bu \rangle = 0$, \ then for each \ $T > 0$,
 \begin{equation}\label{1Zcond3mU}
  \sup_{t\in[0,T]}
   \bigg| \frac{1}{n^2}
          \sum_{k=1}^{\nt}
           U_{k-1} - \langle \bu, \tBbeta \rangle \, \frac{t^2}{2} \bigg|
  \stoch 0 , \qquad \text{as \ $n \to \infty$.}
 \end{equation}
\end{Lem}

\noindent
\textbf{Proof.}
We have
 \begin{align*}
  \bigg| \frac{1}{n^2} \sum_{k=1}^{\nt} U_{k-1}
         - \langle \bu, \tBbeta \rangle \, \frac{t^2}{2} \bigg|
  \leq \frac{1}{n^2}
        \sum_{k=1}^{\nt}
         \Bigl| U_{k-1}
                - \langle \bu, \tBbeta \rangle (k-1) \Bigr|
        + \langle \bu, \tBbeta \rangle
          \bigg| \frac{1}{n^2} \sum_{k=1}^{\nt} (k-1) - \frac{t^2}{2} \bigg| ,
 \end{align*}
 where
 \[
   \sup_{t\in[0,T]}
    \bigg| \frac{1}{n^2} \sum_{k=1}^{\nt} (k-1) - \frac{t^2}{2} \bigg|
   \to 0 , \qquad \text{as \ $n \to \infty$,}
 \]
 hence, in order to show \eqref{1Zcond3mU}, it suffices to prove
 \begin{equation}\label{1Zcond3mmU}
  \frac{1}{n^2}
  \sum_{k=1}^{\nT}
   \Bigl| U_k - \langle \bu, \tBbeta \rangle k \Bigr|
  \stoch 0 , \qquad \text{as \ $n \to \infty$.}
 \end{equation}
Recursion \eqref{rec_U} yields
 \ $\EE(U_k) = \EE(U_{k-1}) + \langle \bu, \tBbeta \rangle$,
 \ $k \in \NN$, \ with intital value \ $\EE(U_0) = 0$, \ hence
 \ $\EE(U_k) = \langle \bu, \tBbeta \rangle k$, \ $k \in \NN$.
\ For the sequence
 \begin{equation}\label{tU}
  \tU_k := U_k - \EE(U_k) = U_k - \langle \bu, \tBbeta \rangle k ,
  \qquad k \in \NN ,
 \end{equation}
 by \eqref{rec_U}, we get a recursion
 \ $\tU_k = \tU_{k-1} + \langle \bu, \bM_k \rangle$, \ $k \in \NN$, \ with
 intital value \ $\tU_0 = 0$.
\ Applying Doob's maximal inequality (see, e.g., Revuz and Yor
 \cite[Chapter II, Theorem 1.7]{RevYor}) for the martingale
 \ $\tU_n = \sum_{k=1}^n \langle \bu, \bM_k \rangle$, \ $n \in \NN$,
 \begin{align*}
  \EE\Biggl(\sup_{t\in[0,T]}
             \Biggl|\sum_{k=1}^\nt
                     \langle \bu, \bM_k \rangle\Biggr|^2\Biggr)
  \leq 4 \EE\Biggl(\Biggl|\sum_{k=1}^\nT
                           \langle \bu, \bM_k \rangle\Biggr|^2\Biggr)
  = 4 \sum_{k=1}^\nT \EE\bigl(\langle \bu, \bM_k \rangle^2\bigr)
  = \OO(n) ,
 \end{align*}
 where we applied Corollary \ref{EEX_EEU_EEV}.
Consequently,
 \begin{equation}\label{supU}
   n^{-1} \max_{k\in\{1,\ldots,\nT\}} |U_k - \langle \bu, \tBbeta \rangle k|
   = n^{-1} \max_{k\in\{1,\ldots,\nT\}} |\tU_k|
   \stoch 0 \qquad \text{as \ $n \to \infty$,}
 \end{equation}
 thus
 \[
   \frac{1}{n^2}
   \sum_{k=1}^{\nT}
    \bigl| U_k - k \langle \bu, \tBbeta \rangle \bigr|
   \leq \frac{\nT}{n^2}
        \max_{k\in\{1,\ldots,\nT\}}
         \bigl| U_k - k \langle \bu, \tBbeta \rangle \bigr|
   \stoch 0 ,
 \]
 as \ $n \to \infty$, \ hence we conclude \eqref{1Zcond3mmU}, thus
 \eqref{1Zcond3mU}.
\proofend

\begin{Lem}\label{main1_UU}
Suppose that the assumptions of Theorem \ref{main} hold.
If \ $\langle \tbC \bu, \bu \rangle = 0$, \ then for each \ $T > 0$,
 \begin{equation}\label{1Zcond3m}
  \sup_{t\in[0,T]}
   \bigg| \frac{1}{n^3}
          \sum_{k=1}^{\nt}
           U_{k-1}^2 - \langle \bu, \tBbeta \rangle^2 \, \frac{t^3}{3} \bigg|
  \stoch 0 , \qquad \text{as \ $n \to \infty$.}
 \end{equation}
\end{Lem}

\noindent
\textbf{Proof.}
We have
 \begin{align*}
  \bigg| \frac{1}{n^3} \sum_{k=1}^{\nt} U_{k-1}^2
         - \langle \bu, \tBbeta \rangle^2 \, \frac{t^3}{3} \bigg|
  \leq \frac{1}{n^3}
        \sum_{k=1}^{\nt}
         \Bigl| U_{k-1}^2
                - \langle \bu, \tBbeta \rangle^2 (k-1)^2 \Bigr|
        + \langle \bu, \tBbeta \rangle^2
          \bigg| \frac{1}{n^3} \sum_{k=1}^{\nt} (k-1)^2 - \frac{t^3}{3} \bigg| ,
 \end{align*}
 where
 \[
   \sup_{t\in[0,T]}
    \bigg| \frac{1}{n^3} \sum_{k=1}^{\nt} (k-1)^2 - \frac{t^3}{3} \bigg|
   \to 0 , \qquad \text{as \ $n \to \infty$,}
 \]
 hence, in order to show \eqref{1Zcond3m}, it suffices to prove
 \begin{equation}\label{1Zcond3mm}
  \frac{1}{n^3}
  \sum_{k=1}^{\nT}
   \Bigl| U_k^2 - \langle \bu, \tBbeta \rangle^2 k^2 \Bigr|
  \stoch 0 , \qquad \text{as \ $n \to \infty$.}
 \end{equation}
We have
 \[
   | U_k^2 - k^2 \langle \bu, \tBbeta \rangle^2 |
   \leq | U_k - k \langle \bu, \tBbeta \rangle |^2
        + 2 k \langle \bu, \tBbeta \rangle
          | U_k - k \langle \bu, \tBbeta \rangle | ,
 \]
 hence, by \eqref{supU},
 \begin{align*}
  n^{-2}
  \max_{k\in\{1,\ldots,\nT\}}
   | U_k^2 - k^2 \langle \bu, \tBbeta \rangle^2 |
  &\leq \Bigl(n^{-1}
              \max_{k\in\{1,\ldots,\nT\}}
               | U_k - k \langle \bu, \tBbeta \rangle | \Bigr)^2 \\
  &\:\quad
        + \frac{2 \nT}{n^2} \langle \bu, \tBbeta \rangle
          \max_{k\in\{1,\ldots,\nT\}}
          | U_k - k \langle \bu, \tBbeta \rangle |
  \stoch 0 ,
 \end{align*}
 as \ $n \to \infty$.
\ Thus,
 \[
   \frac{1}{n^3}
   \sum_{k=1}^{\nT}
    \bigl| U_k^2 - k^2 \langle \bu, \tBbeta \rangle^2 \bigr|
   \leq \frac{\nT}{n^3}
        \max_{k\in\{1,\ldots,\nT\}}
         \bigl| U_k^2 - k^2 \langle \bu, \tBbeta \rangle^2 \bigr|
   \stoch 0 ,
 \]
 as \ $n \to \infty$, \ thus we conclude \eqref{1Zcond3mm}, and hence
 \eqref{1Zcond3m}.
\proofend

For the last convergence in Theorem \ref{main1_Ad} we need the following
 approximation.

\begin{Lem}\label{main_UVt}
Suppose that the assumptions of Theorem \ref{main} hold.
If \ $\langle \tbC \bu, \bu \rangle = 0$, \ then for each \ $T > 0$,
 \[
   \sup_{t\in[0,T]}
    \Biggl| \frac{1}{n^2} \sum_{k=1}^\nt U_k V_k
            - \langle \bu, \tBbeta \rangle
              \langle \bv, \tBbeta \rangle
              \frac{\tdelta t^2}{2(1-\delta)} \Biggr|
   \stoch 0 \qquad
   \text{as \ $n \to \infty$.}
 \]
\end{Lem}

\noindent
\textbf{Proof.}
First we show, by the method of the proof of Lemma \ref{main1_UU}, convergence
 \begin{equation}\label{sumU}
  \sup_{t\in[0,T]}
   \Biggl|\frac{1}{n^2} \sum_{k=1}^\nt U_k
          - \langle \bu, \tBbeta \rangle \frac{t^2}{2} \Biggr|
  \stoch 0 \qquad \text{as \ $n \to \infty$}
 \end{equation}
 for each \ $T > 0$.
\ We have
 \begin{align*}
  \bigg| \frac{1}{n^2} \sum_{k=1}^{\nt} U_{k-1}
         - \langle \bu, \tBbeta \rangle \frac{t^2}{2} \bigg|
  \leq \frac{1}{n^2}
        \sum_{k=1}^{\nt}
         \Bigl| U_{k-1} - \langle \bu, \tBbeta \rangle (k-1) \Bigr|
        + \langle \bu, \tBbeta \rangle
          \bigg| \frac{1}{n^2} \sum_{k=1}^{\nt} (k-1) - \frac{t^2}{2} \bigg| ,
 \end{align*}
 where
 \[
   \sup_{t\in[0,T]}
    \bigg| \frac{1}{n^2} \sum_{k=1}^{\nt} (k-1) - \frac{t^2}{2} \bigg|
   \to 0 , \qquad \text{as \ $n \to \infty$,}
 \]
 hence, in order to show \eqref{sumU}, it suffices to prove
 \begin{equation}\label{sumUm}
  \frac{1}{n^2}
  \sum_{k=1}^{\nT} \Bigl| U_k - \langle \bu, \tBbeta \rangle k \Bigr|
  \stoch 0 , \qquad \text{as \ $n \to \infty$.}
 \end{equation}
Using \eqref{supU}, we obtain
 \[
   \frac{1}{n^2}
   \sum_{k=1}^{\nT}
    \Bigl| U_k - \langle \bu, \tBbeta \rangle k \Bigr|
   \leq \frac{\nT}{n^2}
        \max_{k\in\{1,\ldots,\nT\}}
         \Bigl| U_k - \langle \bu, \tBbeta \rangle k \Bigr|
   \stoch 0 ,
 \]
 as \ $n \to \infty$, \ thus we conclude \eqref{sumUm}, and hence
 \eqref{sumU}.

In order to prove the statement of the lemma, we derive a decomposition of
 \ $\sum_{k=1}^\nt U_k V_k$.
\ Using recursions \eqref{rec_U} and \eqref{rec_V}, we obtain
 \begin{align*}
  \EE(U_k V_k \mid \cF_{k-1})
  &=\EE\Bigl[\bigl(U_{k-1}
                   + \langle \bu, \bM_k + \tBbeta \rangle \bigr)
             \bigl(\delta V_{k-1}
                   + \langle \bv, \bM_k + \tdelta \, \tBbeta \rangle \bigr)
             \,\big|\, \cF_{k-1} \Bigr] \\
  &= \delta U_{k-1} V_{k-1}
     + \tdelta \langle \bv, \tBbeta \rangle U_{k-1}
     + \delta \langle \bu, \tBbeta \rangle V_{k-1}
     + \tdelta \langle \bu, \tBbeta \rangle \langle \bv, \tBbeta \rangle
     + \langle \bV_0 \bu, \bv \rangle ,
 \end{align*}
 since, by \ $\langle \tbC \bu, \bu \rangle = 0$ \ and
 \ $\tbC = (\bV_1 + \bV_2)/2$, \ we conclude
 \ $\langle \bV_i \bu, \bu \rangle  = 0$, \ $i \in \{1, 2\}$, \ thus
 by \eqref{VMk},
 \[
   \EE( \langle \bu, \bM_k \rangle \langle \bv, \bM_k \rangle
        \mid \cF_{k-1})
   = \bu^\top \EE(\bM_k \bM_k^\top \mid \cF_{k-1}) \bv
   = \langle \bV_0 \bu, \bv \rangle .
 \]
Consequently,
 \begin{align*}
  \sum_{k=1}^\nt U_k V_k
  &= \sum_{k=1}^\nt \big[ U_k V_k - \EE(U_k V_k \mid \cF_{k-1}) \big]
     + \sum_{k=1}^\nt \EE(U_k V_k \mid \cF_{k-1}) \\
  &= \sum_{k=1}^\nt \big[ U_k V_k - \EE(U_k V_k \mid \cF_{k-1}) \big]
     + \delta \sum_{k=1}^\nt U_{k-1} V_{k-1}
     + \tdelta \langle \bv, \tBbeta \rangle \sum_{k=1}^\nt U_{k-1} \\
  &\quad
     + \delta \langle \bu, \tBbeta \rangle \sum_{k=1}^\nt V_{k-1}
     + \tdelta \langle \bu, \tBbeta \rangle \langle \bv, \tBbeta \rangle
       \nt
     + \langle \bV_0 \bu, \bv \rangle \nt ,
 \end{align*}
 and we obtain
 \begin{align*}
   \sum_{k=1}^\nt U_k V_k
   &= \frac{1}{1 - \delta}
      \sum_{k=1}^\nt \big[ U_k V_k - \EE(U_k V_k \mid \cF_{k-1}) \big]
      - \frac{\delta}{1 - \delta} U_\nt V_\nt
      + \frac{ \tdelta \langle \bv, \tBbeta \rangle}{1 - \delta}
        \sum_{k=1}^\nt U_{k-1} \\
   &\quad
      + \frac{\delta \langle \bu, \tBbeta \rangle}{1 - \delta}
        \sum_{k=1}^\nt V_{k-1}
      + \frac{\tdelta \langle \bu, \tBbeta \rangle
              \langle \bv, \tBbeta \rangle
              + \langle \bV_0 \bu, \bv \rangle}
             {1 - \delta}
        \nt .
 \end{align*}
Using \eqref{seged_UV_UNIFORM4_mod} with \ $(\ell, i, j) = (4, 1, 1)$ \ we
 obtain
 \begin{align*}
  n^{-2}
  \sup_{t\in[0,T]}
   \Biggl|\sum_{k=1}^\nt \big[ U_k V_k - \EE(U_k V_k \mid \cF_{k-1}) \big]\Biggr|
  \stoch 0 \qquad \text{as \ $n \to \infty$.}
 \end{align*}
Using \eqref{seged_UV_UNIFORM2_mod} with \ $(\ell, i, j) = (3, 1, 1)$ \ we
 obtain \ $n^{-2} \sup_{t\in[0,T]} |U_\nt V_\nt| \stoch 0$ \ as \ $n \to \infty$.
\ The assumption \ $\langle \tbC \bu, \bu \rangle = 0$ \ implies
 \ $\langle \tbC \bv, \bv \rangle = 0$, \ hence, by \eqref{main_Vt}, we
 obtain
 \[
   n^{-2} \sup_{t\in[0,T]} \Biggl|\sum_{k=1}^\nt V_{k-1}\Biggr|
   \stoch 0 .
 \]
Consequently,
 \[
   n^{-2}
   \sup_{t\in[0,T]}
    \Biggl|\sum_{k=1}^\nt U_k V_k
           - \frac{\tdelta \langle \bv, \tBbeta \rangle}{1 - \delta}
             \sum_{k=1}^\nt U_{k-1}\Biggr|
   \stoch 0
 \]
 as \ $n \to \infty$.
\ Using \eqref{seged_UV_UNIFORM2_mod} with \ $(\ell, i, j) = (2, 1, 0)$ \ we
 obtain \ $n^{-2} \sup_{t\in[0,T]} U_\nt \stoch 0$ \ as \ $n \to \infty$.
\ Thus, by \eqref{sumU}, we conclude the statement of the lemma.
\proofend

The proof of the last convergence in Theorem \ref{main1_Ad} is similar to the
 proof of Theorems \ref{main_Ad} and \ref{maint_Ad}.
Consider the sequence of stochastic processes
 \[
   \bcZ^{(n)}_t
   := \begin{bmatrix}
       \bcM_t^{(n)} \\
       \cN_t^{(n)} \\
       \cP_t^{(n)}
      \end{bmatrix}
   := \sum_{k=1}^\nt
       \bZ^{(n)}_k
   \qquad \text{with} \qquad
   \bZ^{(n)}_k
   := \begin{bmatrix}
       n^{-1/2} \bM_k \\
       n^{-3/2} \langle \bu, \bM_k \rangle U_{k-1} \\
       n^{-1/2} \langle \bv, \bM_k \rangle V_{k-1}
      \end{bmatrix}
 \]
 for \ $t \in \RR_+$ \ and \ $k, n \in \NN$.
\ The last convergence in  Theorem \ref{main1_Ad} follows from the following
 theorem.

\begin{Thm}\label{main1_conv}
If \ $\langle \tbC \bu, \bu \rangle = 0$ \ then
 \begin{equation}\label{conv_Z1}
  \bcZ^{(n)} \distr \bcZ , \qquad \text{as \ $n\to\infty$,}
 \end{equation}
 where the process \ $(\bcZ_t)_{t \in \RR_+}$ \ with values in \ $\RR^4$ \ is the
 pathwise unique strong solution of the SDE
 \begin{equation}\label{ZSDE1}
  \dd \bcZ_t = \gamma(t) \, \dd\tbcW_t , \qquad t \in \RR_+ ,
 \end{equation}
 with initial value \ $\bcZ_0 = \bzero$, \ where \ $(\tbcW_t)_{t \in \RR_+}$ \ is
 a 4-dimensional standard Wiener process, and
 \ $\gamma : \RR_+ \to \RR^{4\times4}$ \ is defined by
 \ $\gamma(t) := (\tbS(t))^{1/2}$, \ $t \in \RR_+$, \ with
 \[
   \tbS(t)
   := \begin{bmatrix}
       \bV_0 & \langle \bu, \tBbeta \rangle \bV_0 \bu t
        & \frac{\tdelta \langle \bv, \tBbeta \rangle}{1-\delta}
          \bV_0 \bv \\
       \langle \bu, \tBbeta \rangle \bu^\top \bV_0 t
        & \langle \bu, \tBbeta \rangle^2 \bu^\top \bV_0 \bu t^2
        & \frac{\tdelta \langle \bu, \tBbeta \rangle
                \langle \bv, \tBbeta \rangle}
               {1-\delta}
          \bu^\top \bV_0 \bv t \\[1mm]
       \frac{\tdelta \langle \bv, \tBbeta \rangle}{1-\delta}
       \bv^\top \bV_0
        & \frac{\tdelta \langle \bu, \tBbeta \rangle
                \langle \bv, \tBbeta \rangle}
               {1-\delta}
          \bv^\top \bV_0 \bu t
        & M \bv^\top \bV_0 \bv
      \end{bmatrix} , \qquad
   t \in \RR_+ .
 \]
\end{Thm}

The pathwise unique strong solution of the SDE \eqref{ZSDE1} with initial value
 \ $\bcZ_0 = \bzero$ \ has the form \ $\bcZ_t = \int_0^t \gamma(s) \, \dd \tbcW_s$,
 \ $t \in \RR_+$.
\ Consequently, \ $\bcZ_1 \distre \cN_4(\bzero, \tbS)$ \ with
 \begin{align*}
  \tbS &:= \int_0^1 \gamma(s) \gamma(s)^\top \dd s
        = \begin{bmatrix}
           \bV_0 & \frac{\langle \bu, \tBbeta \rangle}{2} \bV_0 \bu
            & \frac{\tdelta \langle \bv, \tBbeta \rangle}{1-\delta}
              \bV_0 \bv \\[1mm]
           \frac{\langle \bu, \tBbeta \rangle}{2} \bu^\top \bV_0
            & \frac{\langle \bu, \tBbeta \rangle^2}{3}
              \bu^\top \bV_0 \bu
            & \frac{\tdelta \langle \bu, \tBbeta \rangle
                    \langle \bv, \tBbeta \rangle}
                   {2(1-\delta)}
              \bu^\top \bV_0 \bv \\[1mm]
           \frac{\tdelta \langle \bv, \tBbeta \rangle}{1-\delta}
           \bv^\top \bV_0
            & \frac{\tdelta \langle \bu, \tBbeta \rangle
                    \langle \bv, \tBbeta \rangle}
                   {2(1-\delta)}
              \bv^\top \bV_0 \bu
            & M \bv^\top \bV_0 \bv
          \end{bmatrix} \\
      &= \begin{bmatrix}
           \bV_0^{1/2} \\
           \frac{\langle \bu, \tBbeta \rangle}{2}
           \bu^\top \bV_0^{1/2} \\[1mm]
           \frac{\tdelta \langle \bv, \tBbeta \rangle}{1-\delta}
           \bv^\top \bV_0^{1/2}
          \end{bmatrix}
          \begin{bmatrix}
           \bV_0^{1/2} \\
           \frac{\langle \bu, \tBbeta \rangle}{2}
           \bu^\top \bV_0^{1/2} \\[1mm]
           \frac{\tdelta \langle \bv, \tBbeta \rangle}{1-\delta}
           \bv^\top \bV_0^{1/2}
          \end{bmatrix}^\top
        + \frac{\langle \bu, \tBbeta \rangle^2
                  \bu^\top \bV_0 \bu}
                 {12}
            \begin{bmatrix} \bzero \\ 1 \\ 0 \end{bmatrix}
            \begin{bmatrix} \bzero \\ 1 \\ 0 \end{bmatrix}^\top
          + \frac{(\bv^\top \bV_0 \bv)^2}{1-\delta^2}
            \begin{bmatrix} \bzero \\ 0 \\ 1 \end{bmatrix}
            \begin{bmatrix} \bzero \\ 0 \\ 1 \end{bmatrix}^\top .
 \end{align*}
We have
 \[
   \begin{bmatrix}
    n^{-1/2} \langle \bu, \bM_k \rangle \\
    n^{-3/2} \langle \bu, \bM_k \rangle U_{k-1} \\
    n^{-1/2} \langle \bv, \bM_k \rangle \\
    n^{-1/2} \langle \bv, \bM_k \rangle V_{k-1}
   \end{bmatrix}
   = \bA \begin{bmatrix}
          n^{-1/2} \bM_k \\
          n^{-3/2} \langle \bu, \bM_k \rangle U_{k-1} \\
          n^{-1/2} \langle \bv, \bM_k \rangle V_{k-1}
         \end{bmatrix} \qquad \text{with} \qquad
   \bA := \begin{bmatrix}
           \bu^\top & 0 & 0 \\
           \bzero^\top & 1 & 0 \\
           \bv^\top & 0 & 0 \\
           \bzero^\top & 0 & 1
          \end{bmatrix} ,
 \]
 thus the last convergence in Theorem \ref{main1_Ad} follows by the
 continuous mapping theorem, since \ $\bSigma = \bA \tbS \bA^\top$.

\noindent
\textbf{Proof of Theorem \ref{main1_conv}.}
First we observe that the matrix
 \[
  \tbS(t)
  = \begin{bmatrix}
     \bV_0^{1/2} \\
     \langle \bu, \tBbeta \rangle \bu^\top \bV_0^{1/2} t \\[1mm]
     \frac{\tdelta \langle \bv, \tBbeta \rangle}{1-\delta}
     \bv^\top \bV_0^{1/2}
    \end{bmatrix}
    \begin{bmatrix}
     \bV_0^{1/2} \\
     \langle \bu, \tBbeta \rangle \bu^\top \bV_0^{1/2} t \\[1mm]
     \frac{\tdelta \langle \bv, \tBbeta \rangle}{1-\delta}
     \bv^\top \bV_0^{1/2}
    \end{bmatrix}^\top
    + \frac{(\bv^\top \bV_0 \bv)^2}{1-\delta^2}
      \begin{bmatrix} \bzero \\ 0 \\ 1 \end{bmatrix}
      \begin{bmatrix} \bzero \\ 0 \\ 1 \end{bmatrix}^\top
 \]
 is positive semidefinite for all \ $t \in \RR_+$.

We follow again the method of the proof of Theorem \ref{main_conv}.
The conditional variance
 \ $\var\bigl(\bZ^{(n)}_k \mid \cF_{k-1}\bigr)$ \ has the form
 \[
   \begin{bmatrix}
    n^{-1} \bV_{\!\!\bM_k} & n^{-2} U_{k-1} \bV_{\!\!\bM_k} \bu
     & n^{-1} V_{k-1} \bV_{\!\!\bM_k} \bv \\
    n^{-2} U_{k-1} \bu^\top \bV_{\!\!\bM_k}
     & n^{-3} U_{k-1}^2 \bu^\top \bV_{\!\!\bM_k} \bu
     & n^{-2} U_{k-1} V_{k-1} \bu^\top \bV_{\!\!\bM_k} \bv \\
    n^{-1} V_{k-1} \bv^\top \bV_{\!\!\bM_k}
     & n^{-2} U_{k-1} V_{k-1} \bv^\top \bV_{\!\!\bM_k} \bu
     & n^{-1} V_{k-1}^2 \bv^\top \bV_{\!\!\bM_k} \bv
   \end{bmatrix}
 \]
 for \ $n \in \NN$, \ $k \in \{1, \ldots, n\}$, \ with
 \ $\bV_{\!\!\bM_k} = \var( \bM_k \mid \cF_{k-1} )$.
\ Moreover, \ $\gamma(s) \gamma(s)^\top = \tbS(s)$, \ $s \in \RR_+$.

In order to check condition (i) of Theorem \ref{Conv2DiffThm}, we need to
 prove only that for each \ $T > 0$,
 \begin{gather}
  \sup_{t\in[0,T]}
   \bigg| \frac{1}{n} \sum_{k=1}^{\nt} \bV_{\!\!\bM_k}
          - \int_0^t \bV_0 \, \dd s \bigg|
  \stoch 0 , \label{1Zcond1} \\
  \sup_{t\in[0,T]}
   \bigg| \frac{1}{n^2} \sum_{k=1}^{\nt} U_{k-1} \bV_{\!\!\bM_k} \bu
          - \int_0^t
             \langle \bu, \tBbeta \rangle \bV_0 \bu s \, \dd s \bigg|
  \stoch 0 , \label{1Zcond2} \\
  \sup_{t\in[0,T]}
   \bigg| \frac{1}{n^3}
          \sum_{k=1}^{\nt} U_{k-1}^2 \bu^\top \bV_{\!\!\bM_k} \bu
          - \int_0^t
             \langle \bu, \tBbeta \rangle^2 \bu^\top \bV_0 \bu
             s^2 \, \dd s \bigg|
  \stoch 0 , \label{1Zcond3} \\
  \sup_{t\in[0,T]}
   \bigg| \frac{1}{n} \sum_{k=1}^{\nt} V_{k-1} \bV_{\!\!\bM_k} \bv
          - \int_0^t
             \frac{\tdelta\langle \bv, \tBbeta \rangle}{1-\delta}
             \bV_0 \bv \, \dd s \bigg|
  \stoch 0 , \label{1Zcond4} \\
  \sup_{t\in[0,T]}
   \bigg| \frac{1}{n^2}
          \sum_{k=1}^{\nt}
           U_{k-1} V_{k-1} \bu^\top \bV_{\!\!\bM_k} \bv
          - \int_0^t
             \frac{\tdelta\langle \bu, \tBbeta \rangle
                   \langle \bv, \tBbeta \rangle}
                  {1-\delta}
             \bu^\top \bV_0 \bv \, s \, \dd s \bigg|
  \stoch 0 , \label{1Zcond5} \\
  \sup_{t\in[0,T]}
   \bigg| \frac{1}{n^2}
          \sum_{k=1}^{\nt}
           V_{k-1}^2 \bu^\top \bV_{\!\!\bM_k} \bv
          - \int_0^t
             M \bv^\top \bV_0 \bv \, s \, \dd s \bigg|
  \stoch 0 , \label{1Zcond6}
 \end{gather}
 as \ $n \to \infty$.

The assumption \ $\langle \tbC \bu, \bu \rangle = 0$ \ yields
 \ $\bV_{\!\!\bM_k} = \bV_0$, \ hence \eqref{1Zcond2}--\eqref{1Zcond6}
 follow from Lemmas \ref{main1_U}, \ref{main1_UU}, \ref{main_Vt},
 \ref{main_VVt} and \ref{main_UVt}, respectively.

Condition (ii) of Theorem \ref{Conv2DiffThm} can be checked as in the proof of
 Theorem \ref{main_conv}.
\proofend



\appendix

\vspace*{5mm}

\noindent{\bf\Large Appendices}

\section{SDE for multi-type CBI processes}
\label{section_SDE}

For handling \ $\bM_k$, \ $k \in \NN$, \ we need a representation of multi-type
 CBI processes as pathwise unique strong solutions of certain SDEs with jumps.
In what follows we recall some notations and results from Barczy et al.\
 \cite{BarLiPap2}.

Let \ $\cR := \bigcup_{j=0}^d \cR_j$, \ where
 \ $\cR_j$, \ $j \in \{0, 1, \ldots, d\}$, \ are disjoint sets given by
 \[
   \cR_0 := \cU_d \times \{ (\bzero, 0) \}^d
         \subset \RR_+^d \times (\RR_+^d \times \RR_+)^d ,
 \]
 and
 \[
   \cR_j := \{\bzero\} \times \cH_{j,1} \times \cdots \times \cH_{j,d}
         \subset \RR_+^d \times (\RR_+^d \times \RR_+)^d , \qquad
   j \in \{1, \ldots, d\} ,
 \]
 where
 \[
   \cH_{j,i} := \begin{cases}
              \cU_d \times \cU_1
               & \text{if \ $i = j$,} \\
              \{ (\bzero, 0) \} & \text{if \ $i \ne j$.}
             \end{cases}
 \]
Recall that \ $\cU_d = \RR_+^d \setminus \{\bzero\}$, \ and  hence
 \ $\cU_1 = \RR_{++}$.
\ Let \ $m$ \ be the uniquely defined measure on
 \ $\cV := \RR_+^d \times (\RR_+^d \times \RR_+)^d$ \ such that
 \ $m(\cV \setminus \cR) = 0$ \ and its restrictions on \ $\cR_j$,
 \ $j \in \{0, 1, \ldots, d\}$, \ are
 \begin{equation}\label{m}
  m|_{\cR_0}(\dd\br) = \nu(\dd\br) , \qquad
  m|_{\cR_j}(\dd\bz, \dd u) = \mu_j(\dd\bz) \, \dd u ,
  \quad j \in \{1, \ldots, d\} ,
 \end{equation}
 where we identify \ $\cR_0$ \ with \ $\cU_d$ \ and
 \ $\cR_1$, \ldots, $\cR_d$ \ with \ $\cU_d \times \cU_1$ \ in a
 natural way.
Using again this identification, let
 \ $f : \RR^d \times \cV \to \RR_+^d$, \ and
 \ $g : \RR^d \times \cV \to \RR_+^d$, \ be defined by
 \[
   f(\bx, \br)
   := \begin{cases}
       \bz \bbone_{\{\|\bz\|<1\}} \bbone_{\{u \leq x_j\}} ,
        & \text{if \ $\bx = (x_1, \ldots, x_d)^\top \in \RR^d$,
                \ $\br = (\bz, u) \in \cR_j$, \ $j \in \{1, \ldots, d\}$,} \\
       \bzero , & \text{otherwise,}
      \end{cases}
 \]
 \[
   g(\bx, \br)
   := \begin{cases}
       \br, & \text{if \ $\bx \in \RR^d$, \ $\br \in \cR_0$,} \\
       \bz \bbone_{\{\|\bz\|\geq1\}} \bbone_{\{u \leq x_j\}} ,
        & \text{if \ $\bx = (x_1, \ldots, x_d)^\top \in \RR^d$,
                \ $\br = (\bz, u) \in \cR_j$, \ $j \in \{1, \ldots, d\}$,} \\
       \bzero , & \text{otherwise.}
      \end{cases}
 \]
Consider the disjoint decomposition \ $\cR = \cV_0 \cup \cV_1$,
 \ where
 \ $\cV_0 := \bigcup_{j=1}^d \cR_{j,0}$ \ and
 \ $\cV_1 := \cR_0 \cup \bigl( \bigcup_{j=1}^d \cR_{j,1} \bigr)$ \ are
 disjoint
 decompositions with
 \ $\cR_{j,k}
    := \{\bzero\} \times \cH_{j,1,k} \times \cdots \times \cH_{j,d,k}$,
 \ $j \in \{1, \ldots, d\}$, \ $k \in \{0, 1\}$, \ and
 \[
   \cH_{j,i,k} := \begin{cases}
                \cU_{d,k} \times \cU_1 & \text{if \ $i = j$,} \\
                \{ (\bzero, 0) \} & \text{if \ $i \ne j$,}
               \end{cases} \qquad
   \cU_{d,k} := \begin{cases}
              \{ \bz \in \cU_d
                 : \|\bz\| < 1 \} & \text{if \ $k = 0$,} \\
              \{ \bz \in \cU_d
                 : \|\bz\| \geq 1 \} & \text{if \ $k = 1$.}
             \end{cases}
 \]
Note that \ $f(\bx, \br) = \bzero$ \ if \ $\br \in \cV_1$,
 \ $g(\bx, \br) = \bzero$ \ if \ $\br \in \cV_0$, \ hence
 \ $\be_i^\top f(\bx, \br) g(\bx, \br) \be_j = 0$ \ for all
 \ $(\bx, \br) \in \RR^d \times \cV$ \ and \ $i, j \in \{1, \ldots, d\}$.

Consider the following objects:
 \begin{enumerate}
  \item[(E1)]
   a probability space \ $(\Omega, \cF, \PP)$;
  \item[(E2)]
   a \ $d$-dimensional standard Brownian motion \ $(\bW_t)_{t\in\RR_+}$;
  \item[(E3)]
   a stationary Poisson point process \ $p$ \ on \ $\cV$ \ with
    characteristic
    measure \ $m$;
  \item[(E4)]
   a random vector \ $\bxi$ \ with values in \ $\RR_+^d$, \ independent of
    \ $\bW$ \ and \ $p$.
 \end{enumerate}

\begin{Rem}\label{dRM_strong}
Note that if objects (E1)--(E4) are given, then \ $\bxi$, \ $\bW$ \ and \ $p$
 \ are automatically mutually independent according to Remark 3.4 in
 Barczy et al. \cite{BarLiPap1}.
For a short review on point measures and point processes needed for this
 paper, see, e.g., Barczy et al. \cite[Section 2]{BarLiPap1}.
\proofend
\end{Rem}

Provided that the objects (E1)--(E4) are given, let
 \ $(\cF^{\bxi,\bW\!,\,p}_t)_{t\in\RR_+}$ \ denote the augmented filtration
 generated by \ $\bxi$, \ $\bW$ \ and \ $p$, \ see
 Barczy et al.\ \cite{BarLiPap1}.

Let us consider the \ $d$-dimensional SDE
 \begin{equation}\label{SDE_X}
  \begin{aligned}
   \bX_t
   &= \bX_0 + \int_0^t (\Bbeta + \bD \bX_s) \, \dd s
      + \sum_{i=1}^d \be_i \int_0^t \sqrt{2 c_i X_{s,i}^+} \, \dd W_{s,i} \\
   &\quad
      + \int_0^t \int_{\cV_0} f(\bX_{s-}, \br) \, \tN(\dd s, \dd \br)
      + \int_0^t \int_{\cV_1} g(\bX_{s-}, \br) \, N(\dd s, \dd \br) , \qquad
   t \in \RR_+ ,
  \end{aligned}
 \end{equation}
 where \ $\bX_t = (X_{t,1}, \ldots, X_{t,d})^\top$,
 \ $\bD := (d_{i,j})_{i,j\in\{1,\ldots,d\}}$ \ given by
 \[
   d_{i,j} := \tb_{i,j}
             - \int_{\cU_d} z_i \bbone_{\{\|\bz\|\geq1\}} \, \mu_j(\dd\bz) ,
 \]
 $N(\dd s, \dd\br)$ \ is the counting measure of \ $p$ \ on
 \ $\RR_{++} \times \cV$, \ and
 \ $\tN(\dd s, \dd\br) := N(\dd s, \dd\br) - \dd s \, m(\dd\br)$.

\begin{Def}\label{Def_strong_solution2}
Suppose that the objects \textup{(E1)--(E4)} are given.
An \ $\RR_+^d$-valued strong solution of the SDE \eqref{SDE_X} on
 \ $(\Omega, \cF, \PP)$ \ and with respect to the standard Brownian motion
 \ $\bW$, \ the stationary Poisson point process \ $p$ \ and initial value
 \ $\bxi$, \ is an \ $\RR_+^d$-valued
 \ $(\cF^{\bxi,\bW\!,\,p}_t)_{t\in\RR_+}$-adapted c\`{a}dl\`{a}g process
 \ $(\bX_t)_{t \in \RR_+}$ \ such that \ $\PP(\bX_0 = \bxi) = 1$,
 \[
   \PP\biggl( \int_0^t
               \int_{\cV_0} \|f(\bX_s, \br)\|^2 \, \dd s \, m(\dd \br)
                  < \infty \biggr)
        = 1 , \qquad
   \PP\biggl( \int_0^t \int_{\cV_1} \|g(\bX_{s-}, \br)\| \, N(\dd s, \dd \br)
                   < \infty \biggr)
         = 1
 \]
 for all \ $t \in \RR_+$, \ and equation \eqref{SDE_X} holds \ $\PP$-a.s.
\end{Def}

Further, note that the integrals \ $\int_0^t (\Bbeta + \bD \bX_s) \, \dd s$ \ and
 \ $\int_0^t \sqrt{2 c_i X_{s,i}^+} \, \dd W_{s,i}$, \ $i \in \{1, \ldots, d\}$,
 \ exist, since \ $\bX$ \ is c\`{a}dl\`{a}g.
For the following result, see Theorem 4.6 in Barczy et al.\ \cite{BarLiPap2}.

\begin{Thm}\label{strong_solution}
Let \ $(d, \bc, \Bbeta, \bB, \nu, \bmu)$ \ be a set of admissible parameters
 such that the moment condition \eqref{moment_condition_m_new} holds.
Suppose that objects \textup{(E1)--(E4)} are given.
If \ $\EE(\|\bxi\|) < \infty$, \ then there is a pathwise unique
 \ $\RR_+^d$-valued strong solution to the SDE \eqref{SDE_X} with initial value
 \ $\bxi$, \ and the solution is a multi-type CBI process with parameters
 \ $(d, \bc, \Bbeta, \bB, \nu, \bmu)$.
\end{Thm}

We note that the SDE \eqref{SDE_X} can be written in other forms,
 see Barczy et al. \cite[Section 5]{BarLiPap2} for \ $d\in\{1,2\}$ \
 or \eqref{SDE_atirasa_dim2} for \ $d=2$.

Further, one can rewrite the SDE \eqref{SDE_X} in a form which does not contain
 integrals with respect to non-compensated Poisson random measures, and then
 one can perform a linear transformation in order to remove randomness from the
 drift as follows, see Lemma 4.1 in Barczy et al.\ \cite{BarLiPap3}.
This form is very useful for handling \ $\bM_k$, \ $k \in \NN$.

\begin{Lem}\label{SDE_transform_sol}
Let \ $(d, \bc, \Bbeta, \bB, \nu, \bmu)$ \ be a set of admissible parameters
 such that the moment condition \eqref{moment_condition_m_new} holds.
Suppose that objects \textup{(E1)--(E4)} are given with
 \ $\EE(\|\bxi\|) < \infty$.
\ Let \ $(\bX_t)_{t\in\RR_+}$ \ be a pathwise unique \ $\RR_+^d$-valued strong
 solution to the SDE \eqref{SDE_X} with initial value \ $\bxi$.
\ Then
 \[
   \ee^{-t\tbB} \bX_t
   = \bX_0 + \int_0^t \ee^{-s\tbB} \tBbeta \, \dd s
     + \sum_{\ell=1}^d
        \int_0^t \ee^{-s\tbB} \be_\ell \sqrt{2 c_\ell X_{s,\ell}} \, \dd W_{s,\ell}
     + \int_0^t \int_{\cV} \ee^{-s\tbB} h(\bX_{s-}, \br) \, \tN(\dd s, \dd\br)
 \]
 for all \ $t \in \RR_+$, \ where the function
 \ $h : \RR^d \times \cV \to \RR^d$
  \ is defined by \ $h := f + g$, \ hence
 \begin{align*}
  \bX_t
  &= \ee^{(t-s)\tbB} \bX_s + \int_s^t \ee^{(t-u)\tbB} \tBbeta \, \dd u
     + \sum_{\ell=1}^d
        \int_s^t \ee^{(t-u)\tbB}
         \be_\ell \sqrt{2 c_\ell X_{u,\ell}} \, \dd W_{u,\ell} \\
  &\quad
     + \int_s^t \int_{\cV} \ee^{(t-u)\tbB} h(\bX_{u-}, \br) \, \tN(\dd u, \dd\br)
 \end{align*}
 for all \ $s, t \in \RR_+$, \ with \ $s \leq t$.
\ Consequently,
 \[
   \bM_k
   = \sum_{\ell=1}^d
      \int_{k-1}^k
       \ee^{(k-u)\tbB} \be_\ell \sqrt{2 c_\ell X_{u,\ell}} \, \dd W_{u,\ell}
     + \int_{k-1}^k
        \int_{\cV} \ee^{(k-u)\tbB} h(\bX_{u-}, \br) \, \tN(\dd u, \dd\br)
 \]
 for all \ $k \in \NN$.
\end{Lem}

\noindent
\textbf{Proof.}
The last statement follows from \eqref{Mk}, since
 \ $\int_{k-1}^k \ee^{(k-u)\tbB} \tBbeta \, \dd u
    = \int_0^1 \ee^{(1-u)\tbB} \tBbeta \, \dd u = \oBbeta$.
\proofend

Note that the formulas for \ $(\bX_t)_{t\in\RR_+}$ \ and \ $(\bM_k)_{k\in\NN}$
 \ in Lemma \ref{SDE_transform_sol} are generalizations of formulas (3.1) and
 (3.3) in Xu \cite{Xu}, the first displayed formula in the proof of Lemma 2.1
 in Huang et al.\ \cite{HuaMaZhu}, and formulas (1.5) and (1.7) in Li and Ma
 \cite{LiMa}, respectively.

\begin{Lem}\label{SDE_transform_sol_1t}
Let \ $(\bX_t)_{t\in\RR_+}$ \ be a 2-type CBI process with parameters
 \ $(2, \bc, \Bbeta, \bB, \nu, \bmu)$ \ such that \ $\bX_0 = \bzero$,
 \ $\Bbeta \ne \bzero$ \ or \ $\nu \ne 0$, \ and \eqref{tbB2} holds with some
 \ $\gamma \in \RR$ \ and \ $\kappa \in \RR_{++}$ \ such that
 \ $s = \gamma + \kappa = 0$ \ (hence it is irreducible and critical).
Suppose that the moment conditions \eqref{moment_condition_m} hold with
 \ $q = 2$.

If, in addition, \ $\langle \tbC \bv, \bv \rangle = 0$, \ then
 \ $\langle \bv, \bM_k \rangle \ase \langle \bv, \Beta_k \rangle$,
 \ $k \in \NN$, \ with
 \[
   \Beta_k := \int_{k-1}^k \int_{\cR_0} \ee^{(k-s)\tbB} \br \, \tN(\dd s,\dd\br) ,
   \qquad k \in \NN .
 \]

If, in addition, \ $\langle \tbC \bu, \bu \rangle = 0$, \ then
 \ $\langle \bu, \bM_k \rangle \ase \langle \bu, \Beta_k \rangle$,
 \ $k \in \NN$.

The sequence \ $(\Beta_k)_{k\in\NN}$ \ consists of independent and identically
 distributed random vectors.
\end{Lem}

\noindent
\textbf{Proof.}
The assumption \ $\langle \tbC \bv, \bv \rangle = 0$ \ implies
 \ $c_\ell = 0$ \ for each \ $\ell \in \{1, 2\}$
 \ (see the beginning of the proof of Theorem \ref{main_rdb}), thus
 \[
   \langle \bv, \bM_k \rangle
   = \int_{k-1}^k \int_{\cV}
      \langle \bv, \ee^{(k-s)\tbB} h(\bX_{s-}, \br) \rangle
      \, \tN(\dd s, \dd\br)
   = \langle \bv, \Beta_k \rangle + \zeta_{k,1} + \zeta_{k,2}
 \]
 with
 \[
   \zeta_{k,j}
   := \int_{k-1}^k \int_{\cR_j}
       \langle \bv, \ee^{(k-s)\tbB} \bz \rangle \bbone_{\{u\leq X_{s-,j}\}}
       \, \tN(\dd s, \dd\br) , \qquad k \in \NN , \quad j \in \{1, 2\} .
 \]
We have
 \ $\ee^{(k-s)\tbB^\top} \!\! \bv = \ee^{(\gamma-\kappa)(k-s)} \bv$, \ since
 \ $\bv$ \ is a left eigenvector of \ $\ee^{(k-s)\tbB}$ \ belonging to the
 eigenvalue \ $\ee^{(\gamma-\kappa)(k-s)}$, \ hence
 \[
   \zeta_{k,j}
   = \int_{k-1}^k \int_{\cR_j}
      \ee^{(\gamma-\kappa)(k-s)} \langle \bv, \bz \rangle
      \bbone_{\{u\leq X_{s-,j}\}}
      \, \tN(\dd s, \dd\br) ,
   \qquad k \in \NN , \quad j \in \{1, 2\} .
 \]
We have \ $\zeta_{k,j} = I_{k,j} - I_{k-1,j}$, \ $k \in \NN$, \ with
 \ $I_{t,j}
    := \int_0^t \int_{\cR_j}
         \ee^{(\gamma-\kappa)(k-s)} \langle \bv, \bz \rangle
         \bbone_{\{u\leq X_{s-,j}\}}
         \, \tN(\dd s, \dd\br)$,
 \ $t \in \RR_+$.
\ The process \ $(I_{t,j})_{t\in\RR_+}$ \ is a martingale, since
 \begin{align*}
  &\EE\biggl(\int_{k-1}^k \int_{\cU_2} \int_{\cU_1}
              |\ee^{(\gamma-\kappa)(k-s)} \langle \bv, \bz \rangle
               \bbone_{\{u\leq X_{s-,j}\}}|^2
              \, \dd s \, \mu_j(\dd\bz) \, \dd u\biggr) \\
  &\qquad\qquad
   = \int_{k-1}^k \ee^{2(\gamma-\kappa)(k-s)} \EE(X_{s,j}) \, \dd s
     \int_{\cU_2}
      |\langle \bv, \bz \rangle|^2
      \, \mu_j(\dd\bz) \\
  &\qquad\qquad
   \leq \|\bv\|^2
        \int_{k-1}^k \ee^{2(\gamma-\kappa)(k-s)} \EE(X_{s,j}) \, \dd s
        \int_{\cU_2} \|\bz\|^2 \, \mu_j(\dd\bz)
   < \infty ,
 \end{align*}
 see Ikeda and Watanabe \cite[Chapter II, page 62]{IkeWat}, part (vi) of
 Definition \ref{Def_admissible} and moment condition
 \eqref{moment_condition_m} with \ $q = 2$.
\ Consequently, for each \ $k \in \NN$ \ and $j \in \{1, 2\}$, \ we have
 \ $\EE(\zeta_{k,j}) = 0$.

Moreover, the assumption
 \ $\langle \tbC \bv, \bv \rangle = 0$ \ implies
 \ $\int_{\cU_2} \langle \bv, \bz \rangle^2 \, \mu_\ell(\dd\bz) = 0$
 \ for each \ $\ell \in \{1, 2\}$ \ (see the beginning of the proof of Theorem
 \ref{main_rdb}), thus
 \begin{align*}
  \EE(\zeta_{k,j}^2)
  &= \EE\biggl(\int_{k-1}^k \int_{\cU_2} \int_{\cU_1}
              \ee^{2(\gamma-\kappa)(k-s)} \langle \bv, \bz \rangle^2
              \bbone_{\{u\leq X_{s-,j}\}}
              \, \dd s \, \mu_j(\dd\bz) \, \dd u\biggr) \\
  &= \int_{k-1}^k \ee^{2(\gamma-\kappa)(k-s)} \EE(X_{s,j}) \, \dd s
     \int_{\cU_2} \langle \bv, \bz \rangle^2 \, \mu_j(\dd\bz)
   = 0
 \end{align*}
 by Ikeda and Watanabe \cite[Chapter II, Proposition 2.2]{IkeWat}.
Consequently, \ $\zeta_{k,j} \ase 0$, \ and we obtain
 \ $\langle \bv, \bM_k \rangle \ase \langle \bv, \Beta_k \rangle$,
 \ $k \in \NN$.

In a similar way, \ $\langle \tbC \bu, \bu \rangle = 0$ \ implies
 \ $\langle \bu, \bM_k \rangle \ase \langle \bu, \Beta_k \rangle$,
 \ $k \in \NN$.

The Poisson point process \ $p$ \ admits independent increments, hence
 \ $\Beta_k$, \ $k \in \NN$, \ are independent.

For each \ $k \in \NN$, \ the Laplace transform of the random vector
 \ $\Beta_k$ \ has the form
 \begin{align*}
  \EE(\ee^{-\langle\btheta,\Beta_k\rangle})
  &= \exp\biggl\{- \int_{k-1}^k \int_{\cU_2}
                    \left(1 - \ee^{-\langle\btheta,\ee^{(k-s)\tbB}\br\rangle}\right)
                    \dd s \, \nu(\br)\biggr\} \\
  &= \exp\biggl\{- \int_0^1 \int_{\cU_2}
                    \left(1 - \ee^{-\langle\btheta,\ee^{(1-u)\tbB}\br\rangle}\right)
                    \dd u \, \nu(\br)\biggr\}
   = \EE(\ee^{-\langle\btheta,\Beta_1\rangle})
 \end{align*}
 for all \ $\btheta \in \RR_+^2$, \ see, i.e., Kyprianou \cite[page 44]{Kyp},
 hence \ $\Beta_k$, \ $k \in \NN$, \ are identically distributed.
\proofend

\section{On moments of multi-type CBI processes}
\label{section_moments}

In the proof of Theorem \ref{main}, good bounds for moments of the random
 vectors and variables \ $(\bM_k)_{k\in\ZZ_+}$, \ $(\bX_k)_{k\in\ZZ_+}$,
 \ $(U_k)_{k\in\ZZ_+}$ \ and \ $(V_k)_{k\in\ZZ_+}$ \ are extensively used.
The following estimates are proved in Barczy and Pap
 \cite[Lemmas B.2 and B.3]{BarPap}.

\begin{Lem}\label{moment_estimations_X_critical}
Let \ $(\bX_t)_{t\in\RR_+}$ \ be a multi-type CBI process with parameters
 \ $(d, \bc, \Bbeta, \bB, \nu, \bmu)$ \ such that \ $\EE(\|\bX_0\|^q) < \infty$
 \ and the moment conditions \eqref{moment_condition_m} hold with some
 \ $q \in \NN$.
\ Suppose that \ $(\bX_t)_{t\in\RR_+}$ \ is irreducible and critical.
Then
 \begin{equation}\label{moment_ic}
  \sup_{t\in\RR_+} \frac{\EE(\|\bX_t\|^q)}{(1 + t)^q} < \infty .
 \end{equation}
In particular, \ $\EE(\|\bX_t\|^q) = \OO(t^q)$ \ as \ $t \to \infty$ \ in the
 sense that \ $\limsup_{t\to\infty} t^{-q} \EE(\|\bX_t\|^q) < \infty$.
\end{Lem}

\begin{Lem}\label{moment_estimations_1_2}
Let \ $(\bX_t)_{t\in\RR_+}$ \ be a multi-type CBI process with parameters
 \ $(d, \bc, \Bbeta, \bB, \nu, \bmu)$ \ such that \ $\EE(\|\bX_0\|^q) < \infty$
 \ and the moment conditions \eqref{moment_condition_m} hold, where
 \ $q = 2 p$ \ with some \ $p \in \NN$.
\ Suppose that \ $(\bX_t)_{t\in\RR_+}$ \ is irreducible and critical.
Then, for the martingale differences
 \ $\bM_n = \bX_n - \EE(\bX_n \mid \bX_{n-1})$, \ $n \in \NN$, \ we have
 \ $\EE(\|\bM_n\|^{2p}) = \OO(n^p)$ \ as \ $n \to \infty$ \ that is,
 \ $\sup_{n\in \NN} n^{-p} \EE(\|\bM_n\|^{2p}) <\infty$.
\end{Lem}

We have \ $\var(\bM_k \mid \cF_{k-1}) = \var(\bX_k \mid \bX_{k-1})$ \ and
 \ $\var(\bX_k \mid \bX_{k-1} = \bx) = \var(\bX_1 \mid \bX_0 = \bx)$ \ for all
 \ $\bx \in \RR_+^d$, \ since \ $(\bX_t)_{t\in\RR_+}$ \ is a time-homogeneous
 Markov process.
Hence Proposition 4.8 in Barczy et al. \cite{BarLiPap3} implies the following
 formula for \ $\var(\bM_k \mid \cF_{k-1})$.

\begin{Pro}\label{moment_formula_2}
Let \ $(\bX_t)_{t\in\RR_+}$ \ be a multi-type CBI process with parameters
 \ $(d, \bc, \Bbeta, \bB, \nu, \bmu)$ \ such that \ $\EE(\|\bX_0\|^2) < \infty$
 \ and the moment conditions \eqref{moment_condition_m} hold with \ $q = 2$.
\ Then for all \ $k \in \NN$, \ we have
 \[
   \var(\bM_k \mid \cF_{k-1}) = \sum_{i=1}^d (\be_i^\top \bX_{k-1}) \bV_i + \bV_0 ,
 \]
 where
 \begin{align*}
  \bV_i &:= \sum_{\ell=1}^d
             \int_0^1
              \langle \ee^{(1-u)\tbB} \be_i, \be_\ell \rangle
              \ee^{u\tbB} \bC_\ell \ee^{u\tbB^\top}
              \dd u , \qquad i \in \{1, \ldots, d\} , \\
  \bV_0 &:= \int_0^1
             \ee^{u\tbB}
             \left( \int_{\cU_d} \bz \bz^\top \nu(\dd \bz) \right)
             \ee^{u\tbB^\top}
             \dd u
            + \sum_{\ell=1}^d
               \int_0^1
                \left( \int_0^{1-u}
                        \langle \ee^{v\tbB} \tBbeta, \be_ \ell \rangle
                        \, \dd v \right)
                \ee^{u\tbB} \bC_\ell \ee^{u\tbB^\top}
                \dd u .
 \end{align*}
\end{Pro}

\ Note that \ $\bV_0 = \var(\bX_1 \mid \bX_0 = \bzero)$.

\begin{Pro}\label{moment_formula_3}
Let \ $(\bX_t)_{t\in\RR_+}$ \ be a multi-type CBI process with parameters
 \ $(d, \bc, \Bbeta, \bB, \nu, \bmu)$ \ such that \ $\EE(\|\bX_0\|^q) < \infty$
 \ and the moment conditions \eqref{moment_condition_m} hold with some
 \ $q \in \NN$.
\ Then for all \ $j \in \{1, \ldots, q\}$ \ and
 \ $i_1, \ldots i_j \in \{1, \ldots, d\}$, \ there exists
 a polynomial \ $P_{j,i_1,\ldots,i_j} : \RR^d \to \RR$ \ having degree at most
 \ $\lfloor j/2 \rfloor$, \ such that
 \begin{align}\label{polinomP}
  \EE\left( M_{k,i_1} \cdots M_{k,i_j} \mid \cF_{k-1} \right)
  = P_{j,i_1,\ldots,i_j}(\bX_{k-1}) , \qquad k \in \NN ,
 \end{align}
 where \ $\bM_k =: (M_{k,1}, \ldots, M_{k,d})^\top$.
\ The coefficients of the polynomial \ $P_{j,i_1,\ldots,i_j}$ \ depends on
 \ $d$, $\bc$, $\Bbeta$, $\bB$, $\nu$, $\mu_1$, \ldots, $\mu_d$.
\end{Pro}

\noindent
\textbf{Proof.}
We have
 \[
   \EE\left( M_{k,i_1} \cdots M_{k,i_j} \mid \cF_{k-1} \right)
   = \EE\left[ (X_{k,i_1} - \EE(X_{k,i_1} \mid \bX_{k-1})) \cdots
               (X_{k,i_j} - \EE(X_{k,i_j} \mid \bX_{k-1}))
               \mid \bX_{k-1} \right]
 \]
 and
 \begin{align*}
  &\EE\left[ (X_{k,i_1} - \EE(X_{k,i_1} \mid \bX_{k-1})) \cdots
             (X_{k,i_j} - \EE(X_{k,i_j} \mid \bX_{k-1}))
              \mid \bX_{k-1} = \bx \right] \\
  &\qquad\qquad
   =\EE\left[ (X_{1,i_1} - \EE(X_{1,i_1} \mid \bX_0 = \bx)) \cdots
              (X_{1,i_j} - \EE(X_{1,i_j} \mid \bX_0 = \bx))
              \mid \bX_0 = \bx \right]
 \end{align*}
 for all \ $\bx \in \RR_+^d$, \ since \ $(\bX_t)_{t\in\RR_+}$ \ is a
 time-homogeneous Markov process.
Replacing \ $\bw$ \ by \ $\ee^{t\tbB^\top} \!\! \bw$ \ in the formula for
 \ $\EE\bigl[(\bw^\top \ee^{-t\tbB} (\bY_t - \EE(\bY_t))^k\bigr]$ \ in the proof
 of Barczy et al.\ \cite[Theorem 4.5]{BarLiPap3}, and then using the law of
 total probability, one obtains
 \begin{equation}\label{help9}
  \begin{aligned}
   &\EE\left[\langle \bw, (\bX_t - \EE(\bX_t)) \rangle^j\right] \\
   &= j (j - 1)
      \sum_{i=1}^d
       c_i
       \int_0^t
        (\bw^\top \ee^{(t-s)\tbB} \be_i)^2
        \EE\bigl[(\bw^\top \ee^{(t-s)\tbB} (\bX_s - \EE(\bX_s)))^{j-2}
                 X_{s,i}\bigr]
        \dd s \\
   &\quad
    + \sum_{\ell=0}^{j-2}
         \binom{j}{\ell}
         \sum_{i=1}^d
          \int_0^t \int_{\cU_d}
           (\bw^\top \ee^{(t-s)\tbB} \bz)^{j-\ell}
           \EE\Bigl[\bigl(\bw^\top \ee^{(t-s)\tbB}(\bX_s - \EE(\bX_s))\bigr)^\ell
                    X_{s,i}\Bigr]
           \dd s \, \mu_i(\dd\bz) \\
   &\quad
    + \sum_{\ell=0}^{j-2}
         \binom{j}{\ell}
         \int_0^t \int_{\cU_d}
          (\bw^\top \ee^{(t-s)\tbB} \bz)^{j-\ell}
          \EE\Bigl[\bigl(\bw^\top \ee^{(t-s)\tbB}
                         (\bX_s - \EE(\bX_s))\bigr)^\ell\Bigr]
          \dd s \, \nu(\dd\bz)
  \end{aligned}
 \end{equation}
 for all \ $t \in \RR_+$, \ $j \in \{1, \ldots, q\}$ \ and \ $\bw \in \RR^d$,
 \ and hence, for each \ $t \in \RR_+$, \ $j \in \{1, \ldots, q\}$ \ and
 \ $\bw \in \RR^d$, \ there exists a polynomial \ $P_{t,j,\bw} : \RR^d \to \RR$
 \ having degree at most \ $\lfloor j/2 \rfloor$, \ such that
 \[
   \EE\left[ \langle \bw, (\bX_t - \EE(\bX_t)) \rangle^j \right]
   = \EE\bigl[P_{t,j,\bw}(\bX_0)\bigr] ,
 \]
 where the coefficients of the polynomial \ $P_{t,j,\bw}$ \ depends on
 \ $d$, $\bc$, $\Bbeta$, $\bB$, $\nu$, $\mu_1$, \ldots, $\mu_d$.

For all \ $a_1 , \ldots, a_j \in \RR$, \ we have
 \[
   a_1 \cdots a_j
   = \frac{1}{j! 2^j}
     \sum_{\ell_1=0}^1 \ldots \sum_{\ell_j=0}^1
      (-1)^{\ell_1+\cdots+\ell_j}
      \left[ (-1)^{\ell_1} a_1 + \cdots + (-1)^{\ell_j} a_j \right]^j ,
 \]
 see, e.g., the proof of Corollary 4.4 in \cite{BarLiPap3}.
Hence
 \begin{align*}
  &\EE\left[(X_{t,i_1} - \EE(X_{t,i_1})) \cdots
            (X_{t,i_j} - \EE(X_{t,i_j}))\right] \\
  &= \frac{1}{j! 2^j}
     \sum_{\ell_1=0}^1 \ldots \sum_{\ell_j=0}^1
      (-1)^{\ell_1+\cdots+\ell_j}
      \EE\left[\langle (-1)^{\ell_1} \be_{i_1} + \cdots + (-1)^{\ell_j} \be_{i_j},
                       \bX_t - \EE(\bX_t)\rangle^j\right] \\
  &= \frac{1}{j! 2^j}
     \sum_{\ell_1=0}^1 \ldots \sum_{\ell_j=0}^1
      (-1)^{\ell_1+\cdots+\ell_j}
      \EE\bigl[P_{t,j,(-1)^{\ell_1} \be_{i_1} + \cdots + (-1)^{\ell_j} \be_{i_j}}(\bX_0)
         \bigr]
  =: \EE[P_{t,j,i_1,\ldots,i_j}(\bX_0)] ,
 \end{align*}
 which clearly implies the statement with
 \ $P_{j,i_1,\ldots,i_j} := P_{1,j,i_1,\ldots,i_j}$.
\proofend

\begin{Cor}\label{EEX_EEU_EEV}
Let \ $(\bX_t)_{t\in\RR_+}$ \ be a 2-type CBI process with parameters
 \ $(2, \bc, \Bbeta, \bB, \nu, \bmu)$ \ such that \ $\bX_0 = \bzero$,
 \ $\Bbeta \ne \bzero$ \ or \ $\nu \ne 0$, \ and \eqref{tbB2} holds with some
 \ $\gamma \in \RR$ \ and \ $\kappa \in \RR_{++}$ \ such that
 \ $s = \gamma + \kappa = 0$ \ (hence it is irreducible and critical).
Suppose that the moment conditions \eqref{moment_condition_m} hold with some
 \ $q \in \NN$.
\ Then
 \begin{gather*}
  \EE(\|\bX_k\|^i) = \OO(k^i) , \qquad
  \EE(\|\bM_k\|^{2j}) = \OO(k^j) , \qquad
  \EE(U^i_k ) = \OO(k^i) , \qquad
  \EE(V^{2j}_k ) = \OO(k^j)
 \end{gather*}
 for \ $i, j \in \ZZ_+$ \ with \ $i \leq q$ \ and \ $2 j \leq q$.

If, in addition, \ $\langle \tbC \bv, \bv \rangle = 0$, \ then
 \[
   \EE(|\langle \bv, \bM_k \rangle|^i) = \OO(1) , \qquad
   \EE(V_k^{2j}) = \OO(1)
 \]
 for \ $i, j \in \ZZ_+$ \ with \ $i \leq q$ \ and \ $2 j \leq q$.

If, in addition, \ $\langle \tbC \bu, \bu \rangle = 0$, \ then
 \[
   \EE(|\langle \bu, \bM_k \rangle|^i) = \OO(1)
 \]
 for \ $i \in \ZZ_+$ \ with \ $i \leq q$.
\end{Cor}

\noindent
\textbf{Proof.}
The first and second statements follow from Lemmas
 \ref{moment_estimations_X_critical} and \ref{moment_estimations_1_2},
 respectively.
Lemma \ref{moment_estimations_X_critical} implies
 \ $\EE(U^i_k ) = \EE(\langle \bu, \bX_k \rangle^i )
    \leq \EE(\|\bu\|^i \|\bX_k\|^i ) = \OO(k^i)$
 \ for \ $i \in \ZZ_+$ \ with \ $i \leq q$.

By \eqref{rec_V_sol} and the triangular inequality for the \ $L_{2j}$-norm,
 \begin{equation}\label{triang}
   \bigl(\EE(V_k^{2j})\bigr)^{1/(2j)}
   \leq \sum_{\ell=1}^k \delta^{k-\ell}
         \biggl[ \tdelta |\langle \bv, \tBbeta \rangle|
                 + \bigl(\EE(\langle \bv,
                                     \bM_\ell \rangle^{2j})\bigr)^{1/(2j)}
         \biggr] .
 \end{equation}
We have
 \ $\EE(\langle \bv, \bM_\ell \rangle^{2j})
    \leq \|\bv\|^{2j} \EE(\|\bM_\ell\|^{2j})
    = \OO(\ell^j) = \OO(k^j)$,
 \ $\ell \in \{1, \ldots, k\}$, \ and
 \ $\sum_{\ell=1}^k \delta^{k-\ell} < \sum_{i=0}^\infty \delta^i = (1- \delta)^{-1}$,
 \ thus \ $\EE(V^{2j}_k) = \OO(k^j)$, \ $k \in \NN$, \ as desired.

By Lemma \ref{SDE_transform_sol_1t},
 \ $\langle \tbC \bv, \bv \rangle = 0$ \ implies
 \ $\langle \bv, \bM_k \rangle \ase \langle \bv, \Beta_k \rangle$,
 \ $k \in \NN$, \ where \ $\Beta_k$, \ $k \in \NN$, \ are independent and
 identically distributed, thus
 \[
   \EE(|\langle \bv, \bM_k \rangle|^i)
   = \EE(|\langle \bv, \Beta_1 \rangle|^i)
   = \EE(|\langle \bv, \bM_1 \rangle|^i)
   \leq \|\bv\|^i \EE(\|\bM_1\|^i)
   = \OO(1)
 \]
 for \ $i \in \ZZ_+$ \ with \ $i \leq q$.
\ By \eqref{triang} and
 \ $\sum_{\ell=1}^k \delta^{k-\ell} < (1- \delta)^{-1}$, \ we obtain
 \ $\EE(V^{2j}_k) = \OO(1)$ \ for \ $j \in \ZZ_+$ \ with \ $2 j \leq q$.

In a similar way, \ $\langle \tbC \bu, \bu \rangle = 0$ \ yields
 \ $\EE(|\langle \bu, \bM_k \rangle|^i)
    = \EE(|\langle \bu, \Beta_1 \rangle|^i) = \OO(1)$
 \ for \ $i \in \ZZ_+$ \ with \ $i \leq q$.
\proofend

\begin{Cor}\label{LEM_UV_UNIFORM}
Let \ $(\bX_t)_{t\in\RR_+}$ \ be a 2-type CBI process with parameters
 \ $(2, \bc, \Bbeta, \bB, \nu, \bmu)$ \ such that \ $\bX_0 = \bzero$
 \ $\Bbeta \ne \bzero$ \ or \ $\nu \ne 0$, \ and \eqref{tbB2} holds with some
 \ $\gamma \in \RR$ \ and \ $\kappa \in \RR_{++}$ \ such that
 \ $s = \gamma + \kappa = 0$ \ (hence it is irreducible and critical).
Suppose that the moment conditions \eqref{moment_condition_m} hold with some
 \ $\ell \in \NN$.
\ Then
 \begin{itemize}
  \item[\textup{(i)}]
   for all \ $i,j\in\ZZ_+$ \ with \ $\max\{i,j\} \leq \lfloor \ell/2 \rfloor$,
    \ and for all \ $\theta > i + \frac{j}{2} + 1$, \ we have
    \begin{align}\label{seged_UV_UNIFORM1}
     n^{-\theta}
     \sum_{k=1}^n\vert U_k^i V_k^j\vert
     \stoch 0
     \qquad \text{as \ $n\to\infty$,}
    \end{align}
  \item[\textup{(ii)}]
   for all \ $i,j\in\ZZ_+$ \ with \ $\max\{i,j\} \leq \ell$, \ for all \ $T>0$,
    \ and for all \ $\theta > i + \frac{j}{2} + \frac{i+j}{\ell}$, \ we have
    \begin{align}\label{seged_UV_UNIFORM2}
     n^{-\theta} \sup_{t\in[0,T]} \vert U_\nt^i V_\nt^j \vert \stoch 0
     \qquad \text{as \ $n\to\infty$,}
    \end{align}
  \item[\textup{(iii)}]
   for all \ $i,j\in\ZZ_+$ \ with \ $\max\{i,j\} \leq \lfloor \ell/4 \rfloor$,
    \ for all \ $T>0$, \ and \ for all
    \ $\theta > i + \frac{j}{2} + \frac{1}{2}$, \ we have
    \begin{align}\label{seged_UV_UNIFORM4}
     n^{-\theta} \sup_{t\in[0,T]}
     \left|\sum_{k=1}^\nt [U_k^i V_k^j - \EE(U_k^i V_k^j \mid \cF_{k-1})] \right|
     \stoch 0
     \qquad \text{as \ $n\to\infty$.}
    \end{align}
 \end{itemize}

If, in addition, \ $\langle \tbC \bv, \bv \rangle = 0$, \ then
 \begin{itemize}
  \item[\textup{(iv)}]
   for all \ $i,j\in\ZZ_+$ \ with \ $\max\{i,j\} \leq \lfloor \ell/2 \rfloor$,
    \ and for all \ $\theta > i + 1$, \ we have
    \begin{align}\label{seged_UV_UNIFORM1_mod}
     n^{-\theta}
     \sum_{k=1}^n\vert U_k^i V_k^j\vert
     \stoch 0
     \qquad \text{as \ $n\to\infty$,}
    \end{align}
  \item[\textup{(v)}]
   for all \ $i,j\in\ZZ_+$ \ with \ $\max\{i,j\} \leq \ell$, \ for all \ $T>0$,
    \ and for all \ $\theta > i + \frac{i+j}{\ell}$, \ we have
    \begin{align}\label{seged_UV_UNIFORM2_mod}
     n^{-\theta} \sup_{t\in[0,T]} \vert U_\nt^i V_\nt^j \vert \stoch 0
     \qquad \text{as \ $n\to\infty$,}
    \end{align}
  \item[\textup{(vi)}]
   for all \ $i,j\in\ZZ_+$ \ with \ $\max\{i,j\} \leq \lfloor \ell/4 \rfloor$,
    \ for all \ $T>0$, \ and \ for all
    \ $\theta > i + \frac{1}{2}$, \ we have
    \begin{align}\label{seged_UV_UNIFORM4_mod}
     n^{-\theta} \sup_{t\in[0,T]}
     \left|\sum_{k=1}^\nt [U_k^i V_k^j - \EE(U_k^i V_k^j \mid \cF_{k-1})] \right|
     \stoch 0
     \qquad \text{as \ $n\to\infty$,}
    \end{align}
  \item[\textup{(vii)}]
   for all \ $j\in\ZZ_+$ \ with \ $j \leq \lfloor \ell/2 \rfloor$, \ for all
    \ $T > 0$, \ and \ for all \ $\theta > \frac{1}{2}$, \ we have
    \begin{align}\label{seged_UV_UNIFORM4_mod_V}
     n^{-\theta} \sup_{t\in[0,T]}
     \left|\sum_{k=1}^\nt [V_k^j - \EE(V_k^j \mid \cF_{k-1})] \right|
     \stoch 0
     \qquad \text{as \ $n\to\infty$.}
    \end{align}
 \end{itemize}
\end{Cor}

\noindent
\textbf{Proof.}
The first three statements can be derived exactly as in
 Barczy et al.~\cite[Corollary 9.2 of arXiv version]{BarIspPap2}.

If \ $\langle \tbC \bv, \bv \rangle = 0$, \ then, by
 Cauchy--Schwarz's inequality and Corollary \ref{EEX_EEU_EEV}, we have
 \begin{align*}
  \EE\left(\sum_{k=1}^n |U_k^i V_k^j|\right)
  \leq \sum_{k=1}^n \sqrt{\EE(U_k^{2i}) \EE(V_k^{2j})}
  = \sum_{k=1}^n \sqrt{\OO(k^{2i}) \OO(1)}
  = \sum_{k=1}^n \OO(k^{i})
  = \OO(n^{1+i}) .
 \end{align*}
Using Slutsky's lemma this implies \eqref{seged_UV_UNIFORM1_mod}.

Now we turn to prove \eqref{seged_UV_UNIFORM2_mod}.
First note that
 \begin{align}\label{seged_UV_UNIFORM3_mod}
  \sup_{t\in[0,T]} |U_\nt^i V_\nt^j|
  \leq \sup_{t\in[0,T]} |U_\nt^i| \sup_{t\in[0,T]} |V_\nt^j| ,
 \end{align}
 and for all \ $\vare > 0$ \ and \ $\delta > 0$, \ we have, by Markov's
 inequality,
 \begin{align*}
  \PP\biggl( n^{-\vare} \sup_{t\in[0,T]} |U_\nt^i| > \delta \biggr)
  &= \PP\biggl(n^{-\ell\vare/i} \sup_{t\in[0,T]} |U_\nt^\ell| > \delta^{\ell/i}\biggr)
   \leq \sum_{k=1}^\nT \PP(U_k^\ell > \delta^{\ell/i} n^{\ell\vare/i}) \\
  &\leq \sum_{k=1}^\nT \frac{\EE(U_k^\ell)}{\delta^{\ell/i} n^{\ell\vare/i}}
   =\sum_{k=1}^\nT \frac{\OO(k^\ell)}{\delta^{\ell/i} n^{\ell\vare/i}}
   = \OO(n^{\ell+1-\ell\vare/i}),
 \end{align*}
 for all \ $i \in \{1, 2, \ldots, \ell\}$, \ and
 \begin{align*}
  \PP\biggl( n^{-\vare} \sup_{t\in[0,T]} |V_\nt^j| >\delta \biggr)
  &= \PP\biggl(n^{-\ell\vare/j} \sup_{t\in[0,T]} |V_\nt^\ell| > \delta^{\ell/j}\biggr)
   \leq \sum_{k=1}^\nT \PP(|V_k^\ell| > \delta^{\ell/j} n^{\ell\vare/j}) \\
  &\leq \sum_{k=1}^\nT \frac{\EE(|V_k^\ell|)}{\delta^{\ell/j} n^{\ell\vare/j}}
   \leq \sum_{k=1}^\nT \frac{\sqrt{\EE(V_k^{2\ell})}}{\delta^{\ell/j} n^{\ell\vare/j}}
   = \sum_{k=1}^\nT \frac{\OO(1)}{\delta^{\ell/j} n^{\ell\vare/j}}
   = \OO(n^{1-\ell\vare/j})
 \end{align*}
 for all \ $j \in \{1, 2, \ldots, \ell\}$.
\ Hence, if \ $\ell + 1 - \ell\vare/i < 0$, \ i.e.,
 \ $\vare > \frac{\ell+1}{\ell}i$, \ then
 \[
   n^{-\vare} \sup_{t\in[0,T]} |U_\nt^i| \stoch 0
   \qquad \text{as \ $n \to \infty$,}
 \]
 and if \ $1 - \ell\vare/j < 0$, \ i.e., \ $\vare > j/\ell$, \ then
 \[
   n^{-\vare} \sup_{t\in[0,T]} |V_\nt^j| \stoch 0
   \qquad \text{as \ $n \to \infty$.}
 \]
By \eqref{seged_UV_UNIFORM3_mod}, we get \eqref{seged_UV_UNIFORM2_mod}.

Next we show \eqref{seged_UV_UNIFORM4_mod}.
Applying Doob's maximal inequality (see, e.g., Revuz and Yor
 \cite[Chapter II, Theorem 1.7]{RevYor}) for the martingale
 \[
   \sum_{k=1}^n \big[U_k^i V_k^j - \EE(U_k^i V_k^j \mid \cF_{k-1}) \big] ,
   \qquad n \in \NN ,
  \]
 (with the filtration \ $(\cF_k)_{k\in\NN}$), we obtain
 \begin{align*}
  &\EE\left(\sup_{t\in[0,T]}
            \Biggl( \sum_{k=1}^\nt
                     \bigl[U_k^i V_k^j
                           - \EE(U_k^i V_k^j \mid \cF_{k-1}) \bigr]
            \Biggr)^2 \right)
   \leq 4 \EE\left( \Biggl( \sum_{k=1}^\nT
                            \bigl[U_k^i V_k^j
                                  - \EE(U_k^i V_k^j \mid \cF_{k-1}) \bigr]
                   \Biggr)^2 \right)\\
  &= 4 \sum_{k=1}^\nT
        \EE\Bigl(\bigl[U_k^i V_k^j - \EE(U_k^i V_k^j \mid \cF_{k-1})\bigr]^2\Bigr)
   = 4 \sum_{k=1}^\nT
        \biggl\{\EE(U_k^{2i} V_k^{2j})
                - \EE\Bigl(\bigl[\EE(U_k^i V_k^j \mid \cF_{k-1})\bigr]^2\Bigr)
        \biggr\} \\
  &\leq 4 \sum_{k=1}^\nT \EE(U_k^{2i} V_k^{2j})
   = \sum_{k=1}^\nT \OO(k^{2i}) = \OO(n^{2i+1}) ,
 \end{align*}
 since
 \ $\EE(U_k^{2i} V_k^{2j}) \leq \sqrt{\EE(U_k^{4i})\EE(V_k^{4j})} = \OO(k^{2i})$
 \ by Corollary \ref{EEX_EEU_EEV}.
Thus we obtained \eqref{seged_UV_UNIFORM4_mod}.

Finally, we prove \eqref{seged_UV_UNIFORM4_mod_V}.
Applying again Doob's maximal inequality, we have
 \begin{align*}
  \EE\left(\sup_{t\in[0,T]}
            \Biggl( \sum_{k=1}^\nt \big[V_k^j - \EE(V_k^j \mid \cF_{k-1}) \big]
            \Biggr)^2 \right)
  &\leq 4 \EE\left( \Biggl( \sum_{k=1}^\nT
                             \big[V_k^j - \EE(V_k^j \mid \cF_{k-1}) \big]
                    \Biggr)^2 \right)\\
  &\leq 4 \sum_{k=1}^\nT \EE(V_k^{2j})
   = \sum_{k=1}^\nT \OO(1) = \OO(n) ,
 \end{align*}
 since \ $\EE(V_k^{2j}) = \OO(1)$ \ by Corollary \ref{EEX_EEU_EEV}.
Hence we conclude \eqref{seged_UV_UNIFORM4_mod_V}.
\proofend

\section{CLS estimators}
\label{section_estimators}

For the existence of CLS estimators we need the following approximations.

\begin{Lem}\label{main_VV}
Suppose that the assumptions of Theorem \ref{main} hold.
For each \ $T > 0$, \ we have
 \[
   n^{-2}
   \sup_{t\in[0,T]}
    \left| \sum_{k=1}^\nt V_k^2
           - \frac{\langle \tbC \, \bv, \bv \rangle}
                  {1 - \delta^2}
             \sum_{k=1}^\nt U_{k-1} \right|
   \stoch 0 , \qquad \text{as \ $n \to \infty$.}
 \]
\end{Lem}

\noindent
\textbf{Proof.}
In order to prove the statement, we derive a decomposition of
 \ $\sum_{k=1}^\nt V_k^2$.
\ Using recursion \eqref{rec_V}, Proposition \ref{moment_formula_2} and
 \eqref{XUV}, we obtain
 \begin{align*}
  &\EE(V_k^2 \mid \cF_{k-1})
   = \delta^2 V_{k-1}^2
     + 2 \delta\tdelta \langle \bv, \tBbeta \rangle V_{k-1}
     + (\tdelta)^2 \langle \bv, \tBbeta \rangle^2
     + \bv^\top \EE(\bM_k \bM_k^\top \mid \cF_{k-1}) \bv \\
  &= \delta^2 V_{k-1}^2
      + \frac{1}{2} \bv^\top (\bV_1 + \bV_2) \bv \, U_{k-1}
      + \text{constant}
      + \text{constant $\times$ $V_{k-1}$.}
 \end{align*}
Using \eqref{tbC=},
 \begin{align*}
  \sum_{k=1}^\nt V_k^2
  = \sum_{k=1}^\nt \big[ V_k^2 - \EE(V_k^2 \mid \cF_{k-1}) \big]
     + \delta^2 \sum_{k=1}^\nt V_{k-1}^2
     + \bv^\top \tbC \bv \sum_{k=1}^\nt U_{k-1}
     + \OO(n)
     + \text{const.\ $\times$ $\sum_{k=1}^\nt V_{k-1}$.}
 \end{align*}
Consequently,
 \begin{align}\label{sum_Vk2}
  \begin{aligned}
   \sum_{k=1}^\nt V_k^2
   &= \frac{1}{1 - \delta^2}
      \sum_{k=1}^\nt \big[ V_k^2 - \EE(V_k^2 \mid \cF_{k-1}) \big]
      + \frac{1}{1 - \delta^2}
        \langle \tbC \, \bv, \bv \rangle \sum_{k=1}^\nt U_{k-1} \\
   &\quad
      - \frac{\delta^2}{1 - \delta^2} V_\nt^2 + \OO(n)
      + \text{constant $\times$ $\sum_{k=1}^\nt V_{k-1}$.}
  \end{aligned}
 \end{align}
Using \eqref{seged_UV_UNIFORM4} with \ $(\ell, i, j) = (8, 0, 2)$ \ we
 obtain
 \begin{align*}
  \frac{1}{n^2}
  \sup_{t\in[0,T]}
   \left|\sum_{k=1}^\nt \big[ V_k^2 - \EE(V_k^2 \mid \cF_{k-1}) \big]\right|
  \stoch 0 , \qquad \text{as \ $n \to \infty$.}
 \end{align*}
Using \eqref{seged_UV_UNIFORM2} with \ $(\ell, i, j) = (3, 0, 2)$ \ we obtain
 \ $n^{-2} \sup_{t\in[0,T]} V_\nt^2 \stoch 0$.
\ Moreover, \ $n^{-2} \sum_{k=1}^\nt V_{k-1} \stoch 0$ \ as \ $n \to \infty$
 \ follows by \eqref{seged_UV_UNIFORM1} with \ $(\ell, i, j) = (2, 0, 1)$.
\ Consequently, by \eqref{sum_Vk2}, we obtain the statement.
\proofend

\begin{Lem}\label{main_Vt}
Suppose that the assumptions of Theorem \ref{main} hold, and
 \ $\langle \tbC \bv, \bv \rangle = 0$.
\ Then for each \ $T > 0$,
 \[
   \sup_{t\in[0,T]}
    \Biggl| \frac{1}{n} \sum_{k=1}^\nt V_k
            - \frac{\tdelta\langle \bv, \tBbeta \rangle}
                 {1-\delta} t \Biggr|
   \stoch 0 \qquad
   \text{as \ $n \to \infty$.}
 \]
\end{Lem}

\noindent
\textbf{Proof.}
Using recursion \eqref{rec_V}, we obtain
 \[
   \EE(V_k \mid \cF_{k-1})
   = \delta V_{k-1} + \tdelta \langle \bv, \tBbeta \rangle ,
   \qquad k \in \NN .
 \]
Thus
 \[
   \sum_{k=1}^\nt V_k
   = \sum_{k=1}^\nt [V_k - \EE(V_k \mid \cF_{k-1})]
     + \delta \sum_{k=1}^\nt V_{k-1}
     + \nt \tdelta \langle \bv, \tBbeta \rangle .
 \]
Consequently,
 \[
   \sum_{k=1}^\nt V_k
   = \frac{1}{1-\delta} \sum_{k=1}^\nt [V_k - \EE(V_k \mid \cF_{k-1})]
     - \frac{\delta}{1-\delta} V_\nt
     + \nt \frac{\tdelta \langle \bv, \tBbeta \rangle}{1-\delta} .
 \]
Using \eqref{seged_UV_UNIFORM4_mod_V} with \ $(\ell, j) = (2, 1)$ \ we obtain
 \begin{align*}
  n^{-1}
  \sup_{t\in[0,T]}
   \Biggl|\sum_{k=1}^\nt \big[ V_k - \EE(V_k \mid \cF_{k-1}) \big]\Biggr|
  \stoch 0 \qquad \text{as \ $n \to \infty$.}
 \end{align*}
Using \eqref{seged_UV_UNIFORM2_mod} with \ $(\ell, i, j) = (2, 0, 1)$ \ we
 obtain \ $n^{-1} \sup_{t\in[0,T]} |V_\nt| \stoch 0$ \ as \ $n \to \infty$, \ and
 hence we conclude the statement.
\proofend

\begin{Lem}\label{main_VVt}
Suppose that the assumptions of Theorem \ref{main} hold, and
 \ $\langle \tbC \bv, \bv \rangle = 0$.
\ Then for each \ $T > 0$,
 \[
   \sup_{t\in[0,T]} \Biggl| \frac{1}{n} \sum_{k=1}^\nt V_k^2 - Mt \Biggr|
   \stoch 0 \qquad
   \text{as \ $n \to \infty$,}
 \]
 where \ $M$ \ is defined in \eqref{M}.
Moreover, \ $M = 0$ \ if and only if
 \ $(\tbeta_1 - \tbeta_2)^2 + \int_{\cU_2} (z_1 - z_2)^2 \, \nu(\dd\bz) = 0$,
 \ which is equivalent to \ $X_{k,1} \ase X_{k,2}$ \ for all \ $k \in \NN$.
\end{Lem}

\noindent
\textbf{Proof.}
In order to prove the convergence in the statement, we derive a decomposition
 of \ $\sum_{k=1}^\nt V_k^2$.
\ Using recursion \eqref{rec_V}, we obtain
 \begin{align*}
  \EE(V_k^2 \mid \cF_{k-1})
  &=\EE\Bigl[\bigl(\delta V_{k-1}
                   + \langle \bv, \bM_k + \tdelta \, \tBbeta \rangle \bigr)^2
             \,\big|\, \cF_{k-1} \Bigr] \\
  &= \delta^2 V_{k-1}^2
     + 2 \delta \tdelta \langle \bv, \tBbeta \rangle V_{k-1}
     + (\tdelta)^2 \langle \bv, \tBbeta \rangle^2
     + \langle \bV_0 \bv, \bv \rangle ,
 \end{align*}
 since, using \ $\langle \tbC \bv, \bv \rangle = 0$ \ and
 \ $\tbC = (\bV_1 + \bV_2)/2$, \ we conclude
 \ $\langle \bV_i \bv, \bv \rangle  = 0$, \ $i \in \{1, 2\}$, \ thus
 by \eqref{VMk},
 \[
   \EE(\langle \bv, \bM_k \rangle^2 \mid \cF_{k-1})
   = \bv^\top \EE(\bM_k \bM_k^\top \mid \cF_{k-1}) \bv
   = \langle \bV_0 \bv, \bv \rangle .
 \]
Thus
 \begin{align*}
  \sum_{k=1}^\nt V_k^2
  &= \sum_{k=1}^\nt \big[ V_k^2 - \EE(V_k^2 \mid \cF_{k-1}) \big]
     + \sum_{k=1}^\nt \EE(V_k^2 \mid \cF_{k-1}) \\
  &= \sum_{k=1}^\nt \big[ V_k^2 - \EE(V_k^2 \mid \cF_{k-1}) \big]
     + \delta^2 \sum_{k=1}^\nt V_{k-1}^2 \\
  &\quad
     + 2 \delta \tdelta \langle \bv, \tBbeta \rangle \sum_{k=1}^\nt V_{k-1}
     + (\tdelta)^2 \langle \bv, \tBbeta \rangle^2 \nt
     + \langle \bV_0 \bv, \bv \rangle \nt .
 \end{align*}
Consequently,
 \begin{align*}
   \sum_{k=1}^\nt V_k^2
   &= \frac{1}{1 - \delta^2}
      \sum_{k=1}^\nt \big[ V_k^2 - \EE(V_k^2 \mid \cF_{k-1}) \big]
      - \frac{\delta^2}{1 - \delta^2} V_\nt^2 \\
   &\quad
      + \frac{2 \delta \tdelta \langle \bv, \tBbeta \rangle}
             {1 - \delta^2}
        \sum_{k=1}^\nt V_{k-1}
      + \frac{(\tdelta)^2 \langle \bv, \tBbeta \rangle^2
              + \langle \bV_0 \bv, \bv \rangle}
             {1 - \delta^2}
        \nt .
 \end{align*}
Using \eqref{seged_UV_UNIFORM4_mod_V} with \ $(\ell, j) = (4, 2)$ \ we obtain
 \begin{align*}
  n^{-1}
  \sup_{t\in[0,T]}
   \Biggl|\sum_{k=1}^\nt \big[ V_k^2 - \EE(V_k^2 \mid \cF_{k-1}) \big]\Biggr|
  \stoch 0 \qquad \text{as \ $n \to \infty$.}
 \end{align*}
Using \eqref{seged_UV_UNIFORM2_mod} with \ $(\ell, i, j) = (3, 0, 2)$ \ we
 obtain \ $n^{-1} \sup_{t\in[0,T]} V_\nt^2 \stoch 0$ \ as \ $n \to \infty$, \ and
 hence \ $n^{-1} \sup_{t\in[0,T]} |V_\nt| \stoch 0$ \ as \ $n \to \infty$.
\ Consequently, by Lemma \ref{main_Vt}, we obtain
 \begin{align*}
  &\sup_{t\in[0,T]}
    \Biggl|\frac{1}{n} \sum_{k=1}^\nt V_k^2
           - \frac{2 \delta (\tdelta)^2 \langle \bv, \tBbeta \rangle^2}
                  {(1-\delta^2)(1-\delta)} t
           - \frac{(\tdelta)^2 \langle \bv, \tBbeta \rangle^2
                   + \langle \bV_0 \bv, \bv \rangle}
                  {1-\delta^2} t\Biggr|\\
  &=\sup_{t\in[0,T]}
     \Biggl|\frac{1}{n} \sum_{k=1}^\nt V_k^2
            - \frac{(\tdelta)^2 \langle \bv, \tBbeta \rangle^2}
                   {(1-\delta)^2} t
            - \frac{\langle \bV_0 \bv, \bv \rangle}
                   {1-\delta^2} t\Biggr|
   =\sup_{t\in[0,T]} \Biggl|\frac{1}{n} \sum_{k=1}^\nt V_k^2 - M t\Biggr|
    \stoch 0
 \end{align*}
 as \ $n \to \infty$.

Clearly, \ $M = 0$ \ if and only if
 \ $\int_{\cU_2} (z_1 - z_2)^2 \, \nu(\dd \bz) = 0$ \ and
 \ $\tbeta_1 = \tbeta_2$, \ which is equivalent to
 \ $\langle \bV_0 \bv, \bv \rangle = 0$ \ and
 \ $\langle \bv, \tBbeta \rangle = 0$.
\ By \eqref{tbC=} and Proposition \ref{moment_formula_2}, under
 \ $\langle \tbC \bv, \bv \rangle = 0$, \ $M = 0$ is equivalent
 to
 \[
   \EE\bigl(\langle \bv, \bM_k + \tdelta \, \tBbeta \rangle^2\bigr)
   = \sum_{j=1}^2 \be_j^\top \bX_{k-1} \langle \bV_j \bv, \bv \rangle
     + \langle \bV_0 \bv, \bv \rangle
     + (\tdelta)^2 \langle \bv, \tBbeta \rangle^2 = 0
 \]
 for all \ $k \in \NN$, \ which is equivalent to
 \ $\langle \bv, \bM_k + \tdelta \, \tBbeta \rangle \ase 0$ \ for all
 \ $k \in \NN$.
\ By \eqref{rec_V_sol}, this is equivalent to \ $V_k \ase 0$ \ for
 all \ $k \in \NN$, \ which is equivalent to
 \ $X_{k,1} - X_{k,2} = \langle \bv, \bX_k \rangle = V_k \ase 0$ \ for all
 \ $k \in \NN$.
\proofend

For all \ $n \in \NN$ \ and \ $\bx_1, \ldots, \bx_n \in \RR^2$, \ let us put
 \ $\bx^{(n)} := (\bx_1, \ldots, \bx_n)$, \ and in what follows we use the
 convention \ $\bx_0 := \bzero$.
\ For all \ $n \in \NN$, \ we define the
 function
 \ $Q_n : \RR^{2n} \times \RR^4 \to \RR$ \ by
 \[
   Q_n(\bx^{(n)} ; \varrho', \delta', \oBbeta')
   := \sum_{k=1}^n
       \left\| \bx_k
            - \frac{1}{2}
              \begin{bmatrix}
               \varrho' + \delta' & \varrho' - \delta' \\
               \varrho' - \delta' & \varrho' + \delta'
              \end{bmatrix}
              \bx_{k-1}
            - \oBbeta' \right\|^2
 \]
 for all \ $(\varrho', \delta', \oBbeta') \in \RR^4$ \ and \ $\bx^{(n)} \in \RR^{2n}$.
\ By definition, for all \ $n \in \NN$, \ a CLS estimator of the parameters
 \ $(\varrho, \delta, \oBbeta)$ \ is a measurable function
 \ $(\hvarrho, \hdelta, \hoBbeta) : \RR^{2n} \to \RR^4$ \ such that
 \[
   Q_n(\bx^{(n)} ; \hvarrho_n(\bx^{(n)}), \hdelta_n(\bx^{(n)}), \hoBbeta_n(\bx^{(n)}))
   = \inf_{(\varrho', \delta', \oBbeta') \in \RR^4}
      Q_n(\bx^{(n)} ; \varrho', \delta', \oBbeta')
   \qquad \forall\; \bx^{(n)} \in \RR^{2n} .
 \]
Next we give the solutions of this extremum problem.

\begin{Lem}\label{CLSE}
For each \ $n \in \NN$, \ any CLS estimator of the parameters
 \ $(\varrho, \delta, \oBbeta)$ \ is a measurable function
 \ $(\hvarrho, \hdelta, \hoBbeta) : \RR^{2n} \to \RR^4$ \ for which
 \begin{gather}
  \hvarrho_n(\bx^{(n)})
  := \frac{n \sum_{k=1}^n u_k u_{k-1}
           - \sum_{k=1}^n u_k \sum_{k=1}^n u_{k-1}}
          {n \sum_{k=1}^n u_{k-1}^2 - \bigl(\sum_{k=1}^n u_{k-1}\bigr)^2} ,
  \label{CLSErx} \\
  \hdelta_n(\bx^{(n)})
  := \frac{n \sum_{k=1}^n v_k v_{k-1}
           - \sum_{k=1}^n v_k \sum_{k=1}^n v_{k-1}}
          {n \sum_{k=1}^n v_{k-1}^2 - \bigl(\sum_{k=1}^n v_{k-1}\bigr)^2} ,
  \label{CLSEdx} \\
  \hoBbeta_n(\bx^{(n)})
  := \frac{1}{n} \sum_{k=1}^n \bx_k
     - \frac{1}{2n}
       \sum_{k=1}^n
        \begin{bmatrix}
         u_{k-1} & v_{k-1} \\
         u_{k-1} & - v_{k-1}
        \end{bmatrix}
        \begin{bmatrix}
         \hvarrho_n(\bx^{(n)}) \\
         \hdelta_n(\bx^{(n)})
        \end{bmatrix}
  \label{CLSEbx}
 \end{gather}
 if \ $n \sum_{k=1}^n u_{k-1}^2 - \bigl(\sum_{k=1}^n u_{k-1}\bigr)^2 > 0$ \ and
 \ $n \sum_{k=1}^n v_{k-1}^2 - \bigl(\sum_{k=1}^n v_{k-1}\bigr)^2 > 0$, \ where
 \ $u_k := \langle \bu, \bx_k \rangle$ \ and
 \ $v_k := \langle \bv, \bx_k \rangle$, \ $k \in \ZZ_+$.
\end{Lem}

\noindent
\textbf{Proof.}
For any fixed \ $\bx^{(n)} \in \RR^{2n}$ \ with
 \ $n \sum_{k=1}^n u_{k-1}^2 - \bigl(\sum_{k=1}^n u_{k-1}\bigr)^2 > 0$ \ and
 \ $n \sum_{k=1}^n v_{k-1}^2 - \bigl(\sum_{k=1}^n v_{k-1}\bigr)^2 > 0$, \ the
 quadratic function
 \ $\RR^4 \ni (\varrho', \delta', \oBbeta')
          \mapsto Q_n(\bx^{(n)} ; \alpha', \beta', \mu')$
 \ can be written in the form
 \[
   Q_n(\bx^{(n)} ; \varrho', \delta', \oBbeta')
   = \sum_{k=1}^n
      \left\| \bx_k
              - \frac{1}{2}
              \begin{bmatrix}
               u_{k-1} & v_{k-1} \\
               u_{k-1} & - v_{k-1}
              \end{bmatrix}
              \begin{bmatrix}
               \varrho' \\
               \delta'
              \end{bmatrix}
              - \oBbeta' \right\|^2
   = \sum_{k=1}^n \| \bx_k - \bF_{k-1}^\top(\bx^{(n)}) \balpha \|^2 ,
 \]
 with
 \[
   \bF_{k-1}(\bx^{(n)}) := \begin{bmatrix}
             \frac{1}{2}
            \begin{bmatrix}
             u_{k-1} & v_{k-1} \\
             u_{k-1} & - v_{k-1}
            \end{bmatrix}^\top \\
             \bI_2
            \end{bmatrix} , \qquad
   \balpha := \begin{bmatrix}
               \varrho' \\
               \delta' \\
               \oBbeta'
              \end{bmatrix} .
 \]
Consequently,
 \begin{align*}
  &Q_n(\bx^{(n)} ; \varrho', \delta', \oBbeta')
   = \sum_{k=1}^n
      (\bx_k - \bF_{k-1}^\top(\bx^{(n)}) \balpha)^\top
      (\bx_k - \bF_{k-1}^\top(\bx^{(n)}) \balpha) \\
  &\qquad
   = \sum_{k=1}^n
      \left(\balpha^\top \bF_{k-1}(\bx^{(n)}) \bF_{k-1}^\top(\bx^{(n)}) \balpha
       - \balpha^\top \bF_{k-1}(\bx^{(n)}) \bx_k
       - \bx_k^\top \bF_{k-1}^\top(\bx^{(n)}) \balpha + \bx_k^\top \bx_k\right) \\
  &\qquad
   = \balpha^\top \bG_n(\bx^{(n)}) \balpha - \balpha^\top \bh_n(\bx^{(n)})
     - \bh_n(\bx^{(n)})^\top \balpha + \sum_{k=1}^n \bx_k^\top \bx_k ,
 \end{align*}
 with
 \[
   \bG_n(\bx^{(n)})
   := \sum_{k=1}^n \bF_{k-1}(\bx^{(n)}) \bF_{k-1}^\top(\bx^{(n)}) , \qquad
   \bh_n(\bx^{(n)})
   := \sum_{k=1}^n \bF_{k-1}(\bx^{(n)}) \bx_k , \qquad n \in \NN .
 \]
The matrix \ $\bG_n(\bx^{(n)})$ \ is strictly positive definite, since
 \[
   \bG_n(\bx^{(n)}) = \frac{1}{2}
           \sum_{k=1}^n
            \begin{bmatrix}
             u_{k-1}^2 & 0 & u_{k-1} & u_{k-1} \\
             0 & v_{k-1}^2 & v_{k-1} & - v_{k-1} \\
             u_{k-1} & v_{k-1} & 2 & 0 \\
             u_{k-1} & - v_{k-1} & 0 & 2
            \end{bmatrix} ,
 \]
 hence
 \begin{align*}
  \bgamma^\top \bG_n(\bx^{(n)}) \bgamma
  &= \frac{\gamma_1^2}{2} \sum_{k=1}^n u_{k-1}^2
     + \frac{\gamma_2^2}{2} \sum_{k=1}^n v_{k-1}^2
     + n \gamma_3^2 + n \gamma_4^2
     + \gamma_1 \gamma_3 \sum_{k=1}^n u_{k-1} \\
  &\quad
     + \gamma_1 \gamma_4 \sum_{k=1}^n u_{k-1}
     + \gamma_2 \gamma_3 \sum_{k=1}^n v_{k-1}
     - \gamma_2 \gamma_4 \sum_{k=1}^n v_{k-1} \\
  &= n \biggl( \gamma_3 + \frac{\gamma_1}{2n} \sum_{k=1}^n u_{k-1}
             + \frac{\gamma_2}{2n} \sum_{k=1}^n v_{k-1} \biggr)^2
    + n \biggl( \gamma_4 + \frac{\gamma_1}{2n} \sum_{k=1}^n u_{k-1}
               - \frac{\gamma_2}{2n} \sum_{k=1}^n v_{k-1} \biggr)^2 \\
  &\quad
    + \frac{\gamma_1^2}{2}
      \Biggl[ \sum_{k=1}^n u_{k-1}^2
               - \biggl( \sum_{k=1}^n u_{k-1} \biggr)^2 \Biggr]
    + \frac{\gamma_2^2}{2}
      \Biggl[ \sum_{k=1}^n v_{k-1}^2
              - \biggl( \sum_{k=1}^n v_{k-1} \biggr)^2 \Biggr]
  > 0
 \end{align*}
 for all
 \ $\bgamma = (\gamma_1, \gamma_2, \gamma_3, \gamma_4)^\top
    \in \RR^4 \setminus \{\bzero\}$.
\ Thus
 \begin{align*}
  Q_n(\bx^{(n)} ; \varrho', \delta', \oBbeta')
  &= (\balpha - \bG_n^{-1}(\bx^{(n)}) \bh_n(\bx^{(n)}))^\top \bG_n(\bx^{(n)})
     (\balpha - \bG_n^{-1}(\bx^{(n)}) \bh_n(\bx^{(n)})) \\
  &\quad
     - \bh_n(\bx^{(n)})^\top \bG_n^{-1}(\bx^{(n)}) \bh_n(\bx^{(n)})
     + \sum_{k=1}^n \bx_k^\top \bx_k \\
  &\geq - \bh_n(\bx^{(n)})^\top \bG_n^{-1}(\bx^{(n)}) \bh_n(\bx^{(n)})
        + \sum_{k=1}^n \bx_k^\top \bx_k ,
 \end{align*}
 and equality holds if and only if
 \ $\balpha = \bG_n^{-1}(\bx^{(n)}) \bh_n(\bx^{(n)})$.
\ Consequently,
 \begin{equation}\label{CLSErdb}
   \begin{bmatrix}
    \hvarrho_n(\bx^{(n)}) \\
    \hdelta_n(\bx^{(n)}) \\
    \hoBbeta_n(\bx^{(n)})
   \end{bmatrix}
   = \bG_n^{-1}(\bx^{(n)}) \bh_n(\bx^{(n)}) , \qquad \text{hence} \qquad
   \bG_n(\bx^{(n)})
   \begin{bmatrix}
    \hvarrho_n(\bx^{(n)}) \\
    \hdelta_n(\bx^{(n)}) \\
    \hoBbeta_n(\bx^{(n)})
   \end{bmatrix}
   = \bh_n(\bx^{(n)}) .
 \end{equation}
The last two coordinates of the second equation in \eqref{CLSErdb} have the form
 \[
   \frac{1}{2}
   \sum_{k=1}^n
    \begin{bmatrix}
     u_{k-1} & v_{k-1} \\
     u_{k-1} & - v_{k-1}
    \end{bmatrix}
   \begin{bmatrix}
    \hvarrho_n(\bx^{(n)}) \\
    \hdelta_n(\bx^{(n)})
   \end{bmatrix}
   + n \hoBbeta_n(\bx^{(n)})
   = \sum_{k=1}^n \bx_k ,
 \]
 hence we obtain \eqref{CLSEbx}.
The first two coordinates of the second equation in \eqref{CLSErdb} have the form
 \[
   \frac{1}{2}
   \sum_{k=1}^n
    \begin{bmatrix}
     u_{k-1}^2 & 0 \\
     0 & v_{k-1}^2
    \end{bmatrix}
   \begin{bmatrix}
    \hvarrho_n(\bx^{(n)}) \\
    \hdelta_n(\bx^{(n)})
   \end{bmatrix}
   + \frac{1}{2}
     \sum_{k=1}^n
      \begin{bmatrix}
       u_{k-1} & u_{k-1} \\
       v_{k-1} & - v_{k-1}
      \end{bmatrix}
      \hoBbeta_n(\bx^{(n)}) \\
   = \frac{1}{2}
     \sum_{k=1}^n
      \begin{bmatrix}
       u_{k-1} & u_{k-1} \\
       v_{k-1} & - v_{k-1}
      \end{bmatrix}
     \bx_k .
 \]
Using \eqref{CLSEbx}, we conclude
 \begin{align*}
  &\Biggl(n \sum_{k=1}^n
      \begin{bmatrix}
       u_{k-1}^2 & 0 \\
       0 & v_{k-1}^2
      \end{bmatrix}
   - \frac{1}{2}
     \sum_{k=1}^n
      \begin{bmatrix}
       u_{k-1} & u_{k-1} \\
       v_{k-1} & - v_{k-1}
      \end{bmatrix}
     \sum_{k=1}^n
      \begin{bmatrix}
       u_{k-1} & v_{k-1} \\
       u_{k-1} & - v_{k-1}
      \end{bmatrix}\Biggr)
     \begin{bmatrix}
      \hvarrho_n(\bx^{(n)}) \\
      \hdelta_n(\bx^{(n)})
     \end{bmatrix} \\
  &= n \sum_{k=1}^n
        \begin{bmatrix}
         u_{k-1} & u_{k-1} \\
         v_{k-1} & - v_{k-1}
        \end{bmatrix}
     \bx_k
     - \sum_{k=1}^n
        \begin{bmatrix}
         u_{k-1} & u_{k-1} \\
         v_{k-1} & - v_{k-1}
        \end{bmatrix}
       \sum_{k=1}^n \bx_k ,
 \end{align*}
 which can be written in the form
 \begin{align*}
  &\begin{bmatrix}
    n \sum_{k=1}^n u_{k-1}^2 - \bigl(\sum_{k=1}^n u_{k-1}\bigr)^2 & 0 \\
    0 & n \sum_{k=1}^n v_{k-1}^2 - \bigl(\sum_{k=1}^n v_{k-1}\bigr)^2
   \end{bmatrix}
   \begin{bmatrix}
    \hvarrho_n(\bx^{(n)}) \\
    \hdelta_n(\bx^{(n)})
   \end{bmatrix} \\
  &= \begin{bmatrix}
      n \sum_{k=1}^n u_k u_{k-1} - \sum_{k=1}^n u_k \sum_{k=1}^n u_{k-1} \\
      n \sum_{k=1}^n v_k v_{k-1} - \sum_{k=1}^n v_k \sum_{k=1}^n v_{k-1}
     \end{bmatrix} ,
 \end{align*}
 implying \eqref{CLSErx} and \eqref{CLSEdx}.
\proofend

\begin{Lem}\label{LEMMA_CLSE_exist_discrete}
Suppose that the assumptions of Theorem \ref{main} hold.
Then \ $\PP(H_n) \to 1$ \ as \ $n \to \infty$, \ and hence, the probability of
 the existence of a unique CLS estimator \ $\hvarrho_n$ \ converges to 1 as
 \ $n \to \infty$, \ and this CLS estimator has the form given in \eqref{CLSEr}
 on the event \ $H_n$.

If, in addition,
 \ $\|\bc\|^2 + \sum_{i=1}^2 \int_{\cU_2} (z_1 - z_2)^2 \, \mu_i(\dd\bz)
    + \int_{\cU_2} (z_1 - z_2)^2 \, \nu(\dd\bz) > 0$, \ then
 \ $\PP(\tH_n) \to 1$ \ as \ $n \to \infty$, \ and hence, the probability of the
 existence of a unique CLS estimator \ $\hdelta_n$ \ converges to 1 as
 \ $n \to \infty$, \ and this CLS estimator has the form given in \eqref{CLSEd}
 on the event \ $\tH_n$.
\ Consequently, the probability of the existence of the estimator \ $\hoBbeta_n$
 \ converges to 1 as \ $n \to \infty$, \ and this estimator has the form given
 in \eqref{CLSEb} on the event \ $H_n \cap \tH_n$.
\end{Lem}

\noindent
\textbf{Proof.}
First, note that for all \ $n \in \NN$,
 \begin{align*}
  \Omega \setminus H_n
  &= \biggl\{\omega \in \Omega
             : \sum_{k=1}^n \langle \bu, \bX_{k-1}(\omega) \rangle^2
               - \frac{1}{n}
                 \biggl(\sum_{i=1}^n
                         \langle \bu, \bX_{i-1}(\omega) \rangle\biggr)^2
               = 0\biggr\} \\
  &= \biggl\{\omega \in \Omega
             : \sum_{k=1}^n \biggl(\langle \bu, \bX_{k-1}(\omega)\rangle
               - \frac{1}{n}
                 \sum_{i=1}^n
                  \langle \bu, \bX_{i-1}(\omega) \rangle\biggr)^2
               = 0\biggr\} \\
  &= \biggl\{\omega \in \Omega
             : \langle \bu, \bX_{k-1}(\omega) \rangle
               = \frac{1}{n}
                 \sum_{i=1}^n \langle \bu, \bX_{i-1}(\omega) \rangle ,
               \, k \in \{1, \ldots, n\}\biggr\} \\
  &= \biggl\{\omega \in \Omega
             : 0 = \langle \bu, \bX_0(\omega) \rangle
               = \langle \bu, \bX_1(\omega) \rangle = \cdots
               = \langle \bu, \bX_{n-1}(\omega) \rangle\biggr\} \\
  &= \biggl\{\omega \in \Omega
             : \frac{1}{n^2}
               \sum_{i=1}^n \langle \bu, \bX_{i-1}(\omega) \rangle
               = 0\biggr\} ,
 \end{align*}
 where we used that \ $\bX_0 = \bzero$, \ $\bX_k \in \RR_+^2$,
 \ $k \in \ZZ_+$, \ and \ $\bu \in \RR_{++}$.

By continuous mapping theorem, we obtain
 \begin{align}\label{seged2}
  \frac{1}{n^2} \sum_{k=1}^n \langle \bu, \bX_{k-1} \rangle
  \distr \int_0^1 \cZ_t \, \dd t \qquad
  \text{as \ $n \to \infty$,}
 \end{align}
 where \ $\cZ$ \ is given by the SDE \eqref{SDE_Y},
 see, e.g., the method of the proof of Proposition 3.1 in
 Barczy et al.\ \cite{BarIspPap1}.

We have already proved \ $\PP\bigl( \int_0^1 \cZ_t \, \dd t > 0 \bigr) = 1$.
\ Thus the distribution function of \ $\int_0^1 \cZ_t \, \dd t$ \ is
 continuous at 0, and hence, by \eqref{seged2},
 \[
   \PP(H_n)
   = \PP\biggl(\frac{1}{n^2}
               \sum_{k=1}^n \langle \bu, \bX_{k-1} \rangle > 0\biggr)
   \to \PP\left( \int_0^1 \cZ_t \, \dd t > 0 \right) = 1
 \]
 as \ $n \to \infty$.

In a similar way, for all \ $n \in \NN$,
 \begin{align*}
  \Omega \setminus \tH_n
  &= \{\omega \in \Omega
       : 0 = \langle \bv, \bX_0(\omega) \rangle
         = \langle \bv, \bX_1(\omega) \rangle = \cdots
         = \langle \bv, \bX_{n-1}(\omega) \rangle\} \\
  &= \biggl\{\omega \in \Omega
             : \frac{1}{n^2}
               \sum_{i=1}^n \langle \bv, \bX_{i-1}(\omega) \rangle^2
               = 0\biggr\} .
 \end{align*}

Now suppose that
  \ $\|\bc\|^2 + \sum_{i=1}^2 \int_{\cU_2} (z_1 - z_2)^2 \, \mu_i(\dd\bz) > 0$
 \ holds.
Then \ $\langle \tbC \bv, \bv \rangle > 0$, \ see
 \eqref{tbCbv_LEFT,bv_LEFT=0}.
Applying Lemma \ref{main_VV} and \eqref{seged2}, we obtain
 \[
   \frac{1}{n^2} \sum_{k=1}^n \langle \bv, \bX_{k-1}\rangle^2
   \distr
   \frac{\langle \tbC \, \bv, \bv \rangle}{1-\delta^2}
   \int_0^1 \cZ_t \, \dd t
   \qquad \text{as \ $n \to \infty$,}
 \]
 thus
 \[
   \PP(\tH_n)
   = \PP\biggl(\frac{1}{n^2}
               \sum_{k=1}^n \langle \bv, \bX_{k-1} \rangle^2 > 0\biggr)
   \to \PP\biggl(\frac{\langle \tbC \, \bv, \bv \rangle}{1-\delta^2}
                 \int_0^1 \cZ_t \, \dd t > 0\biggr)
   = 1 ,
 \]
 hence we obtain the statement under the assumption
 \ $\|\bc\|^2 + \sum_{i=1}^2 \int_{\cU_2} (z_1 - z_2)^2 \, \mu_i(\dd\bz) > 0$.

Next we suppose that
 \ $\|\bc\|^2 + \sum_{i=1}^2 \int_{\cU_2} (z_1 - z_2)^2 \, \mu_i(\dd\bz) = 0$
 \ and \ $\int_{\cU_2} (z_1 - z_2)^2 \, \nu(\dd\bz) > 0$.
\ Then \ $\langle \tbC \bv, \bv \rangle = 0$ \ and
 \ $\langle \bV_0 \bv, \bv \rangle > 0$ \ (see
 \eqref{tbCbv_LEFT,bv_LEFT=0} and \eqref{bV_0bv_LEFT,bv_LEFT}), \ hence
 \ $M > 0$.
\ Applying Lemma \ref{main_VVt}, we obtain
 \[
   \frac{1}{n^2} \sum_{k=1}^n \langle \bv, \bX_{k-1}\rangle^2
   \stoch M
   \qquad \text{as \ $n \to \infty$,}
 \]
 thus
 \[
   \PP(\tH_n)
   = \PP\biggl(\frac{1}{n^2}
               \sum_{k=1}^n \langle \bv, \bX_{k-1} \rangle^2 > 0\biggr)
   \to 1 ,
 \]
 hence we conclude the statement under the assumptions
 \ $\|\bc\|^2 + \sum_{i=1}^2 \int_{\cU_2} (z_1 - z_2)^2 \, \mu_i(\dd\bz) = 0$
 \ and \ $\int_{\cU_2} (z_1 - z_2)^2 \, \nu(\dd\bz) > 0$.
\proofend

\section{A version of the continuous mapping theorem}
\label{CMT}

The following version of continuous mapping theorem can be found for example
 in Kallenberg \cite[Theorem 3.27]{K}.

\begin{Lem}\label{Lem_Kallenberg}
Let \ $(S, d_S)$ \ and \ $(T, d_T)$ \ be metric spaces and
 \ $(\xi_n)_{n \in \NN}$, \ $\xi$ \ be random elements with values in \ $S$
 \ such that \ $\xi_n \distr \xi$ \ as \ $n \to \infty$.
\ Let \ $f : S \to T$ \ and \ $f_n : S \to T$, \ $n \in \NN$, \ be measurable
 mappings and \ $C \in \cB(S)$ \ such that \ $\PP(\xi \in C) = 1$ \ and
 \ $\lim_{n \to \infty} d_T(f_n(s_n), f(s)) = 0$ \ if
 \ $\lim_{n \to \infty} d_S(s_n,s) = 0$ \ and \ $s \in C$.
\ Then \ $f_n(\xi_n) \distr f(\xi)$ \ as \ $n \to \infty$.
\end{Lem}

For the case \ $S = \DD(\RR_+, \RR^d)$ \ and \ $T = \RR^q$
 \ (or \ $T = \DD(\RR_+,\RR^q)$), \ where \ $d$, \ $q \in \NN$, \ we formulate
 a consequence of Lemma \ref{Lem_Kallenberg}.

For functions \ $f$ \ and \ $f_n$, \ $n \in \NN$, \ in \ $\DD(\RR_+, \RR^d)$,
 \ we write \ $f_n \lu f$ \ if \ $(f_n)_{n \in \NN}$ \ converges to \ $f$
 \ locally uniformly, that is, if \ $\sup_{t \in [0,T]} \|f_n(t) - f(t)\| \to 0$
 \ as \ $n \to \infty$ \ for all \ $T > 0$.
\ For measurable mappings
 \ $\Phi : \DD(\RR_+, \RR^d) \to \RR^q$
 \ (or \ $\Phi : \DD(\RR_+, \RR^d) \to \DD(\RR_+,\RR^q)$) \ and
 \ $\Phi_n : \DD(\RR_+, \RR^d) \to \RR^q$
 \ (or \ $\Phi_n : \DD(\RR_+, \RR^d) \to \DD(\RR_+,\RR^q)$), \ $n \in \NN$,
 \ we will denote by \ $C_{\Phi, (\Phi_n)_{n \in \NN}}$ \ the set of all functions
 \ $f \in \CC(\RR_+, \RR^d)$ \ such that
 \ $\Phi_n(f_n) \to \Phi(f)$ \ (or \ $\Phi_n(f_n) \to \lu \Phi(f)$) \ whenever
 \ $f_n \lu f$ \ with \ $f_n \in \DD(\RR_+, \RR^d)$, \ $n \in \NN$.

We will use the following version of the continuous mapping theorem several
 times, see, e.g., Barczy et al. \cite[Lemma 4.2]{BarIspPap1} and Isp\'any and
 Pap \cite[Lemma 3.1]{IspPap}.

\begin{Lem}\label{Conv2Funct}
Let \ $d, q \in \NN$, \ and \ $(\bcU_t)_{t\in\RR_+}$ \ and
 \ $(\bcU^{(n)}_t)_{t\in\RR_+}$, \ $n \in \NN$, \ be \ $\RR^d$-valued stochastic
 processes with c\`adl\`ag paths such that \ $\bcU^{(n)} \distr \bcU$.
\ Let \ $\Phi : \DD(\RR_+, \RR^d) \to \RR^q$
 \ (or \ $\Phi : \DD(\RR_+, \RR^d) \to \DD(\RR_+,\RR^q)$) \ and
 \ $\Phi_n : \DD(\RR_+, \RR^d) \to \RR^q$
 \ (or \ $\Phi_n : \DD(\RR_+, \RR^d) \to \DD(\RR_+,\RR^q)$), \ $n \in \NN$,
 \ be measurable mappings such that there exists
 \ $C \subset C_{\Phi,(\Phi_n)_{n\in\NN}}$ \ with \ $C \in \cD_\infty(\RR_+, \RR^d)$
 \ and \ $\PP(\bcU \in C) = 1$.
\ Then \ $\Phi_n(\bcU^{(n)}) \distr \Phi(\bcU)$.
\end{Lem}

In order to apply Lemma \ref{Conv2Funct}, we will use the following statement
 several times, see Barczy et al.~\cite[Lemma B.3]{BarIspPap2}.

\begin{Lem}\label{Marci}
Let \ $d, p, q \in \NN$, $h:\RR^d\to\RR^q$ \ be a continuous function and
 \ $K : [0,1] \times \RR^{2d} \to \RR^p$ \ be a function such that for all
 \ $R > 0$ \ there exists \ $C_R > 0$ such that
 \begin{equation}\label{Lipschitz}
  \| K(s, x) - K(t, y) \| \leq C_R \left( | t - s | + \| x - y \| \right)
 \end{equation}
 for all \ $s, t \in [0, 1]$ \ and \ $x, y \in \RR^{2d}$ \ with
 \ $\| x \| \leq R$ \ and \ $\| y \| \leq R$.
\ Moreover, let us define the mappings
 \ $\Phi, \Phi_n : \DD(\RR_+, \RR^d) \to \RR^{q+p}$, \ $n \in \NN$, \ by
 \begin{align*}
  \Phi_n(f)
  &:= \left( h(f(1)),
             \frac{1}{n}
             \sum_{k=1}^n
              K\left( \frac{k}{n}, f\left( \frac{k}{n} \right),
                      f\left( \frac{k-1}{n} \right) \right) \right) , \\
  \Phi(f)
  & := \left(  h(f(1)), \int_0^1 K( u, f(u), f(u) ) \, \dd u \right)
 \end{align*}
 for all \ $f \in \DD(\RR_+, \RR^d)$.
\ Then the mappings \ $\Phi$ \ and \ $\Phi_n$, \ $n \in \NN$, \ are measurable,
 and
 \ $C_{\Phi,(\Phi_n)_{n \in \NN}} = \CC(\RR_+, \RR^d) \in \cD_\infty(\RR_+, \RR^d)$.
\end{Lem}

\section{Convergence of random step processes}
\label{section_conv_step_processes}

We recall a result about convergence of random step processes towards a
 diffusion process, see Isp\'any and Pap \cite{IspPap}.
This result is used for the proof of convergence \eqref{conv_Z}.

\begin{Thm}\label{Conv2DiffThm}
Let \ $\bgamma : \RR_+ \times \RR^d \to \RR^{d \times r}$ \ be a continuous
 function.
Assume that uniqueness in the sense of probability law holds for the SDE
 \begin{equation}\label{SDE}
  \dd \, \bcU_t
  = \gamma (t, \bcU_t) \, \dd \bcW_t ,
  \qquad t \in \RR_+,
 \end{equation}
 with initial value \ $\bcU_0 = \bu_0$ \ for all \ $\bu_0 \in \RR^d$, \ where
 \ $(\bcW_t)_{t \in \RR_+}$ \ is an $r$-dimensional standard Wiener process.
Let \ $(\bcU_t)_{t \in \RR_+}$ \ be a solution of \eqref{SDE} with initial value
 \ $\bcU_0 = \bzero \in \RR^d$.

For each \ $n \in \NN$, \ let \ $(\bU^{(n)}_k)_{k \in \NN}$ \ be a sequence of
 $d$-dimensional martingale differences with respect to a filtration
 \ $(\cF^{(n)}_k)_{k \in \ZZ_+}$, \ that is,
 \ $\EE(\bU^{(n)}_k \mid \cF^{(n)}_{k-1}) = 0$, \ $n \in \NN$, \ $k \in \NN$.
\ Let
 \[
   \bcU^{(n)}_t := \sum_{k=1}^{\nt} \bU^{(n)}_k \, ,
   \qquad t \in \RR_+, \quad n \in \NN .
 \]
Suppose that \ $\EE \big( \|\bU^{(n)}_k\|^2 \big) < \infty$ \ for all
 \ $n, k \in \NN$.
\ Suppose that for each \ $T > 0$,
 \begin{enumerate}
  \item [\textup{(i)}]
        $\sup\limits_{t\in[0,T]}
         \left\| \sum\limits_{k=1}^{\nt}
                  \var\bigl( \bU^{(n)}_k \mid \cF^{(n)}_{k-1} \bigr)
                 - \int_0^t
                    \bgamma(s,\bcU^{(n)}_s) \bgamma(s,\bcU^{(n)}_s)^\top
                    \dd s \right\|
         \stoch 0$,\\
  \item [\textup{(ii)}]
        $\sum\limits_{k=1}^{\lfloor nT \rfloor}
          \EE \big( \|\bU^{(n)}_k\|^2 \bbone_{\{\|\bU^{(n)}_k\| > \theta\}}
                    \bmid \cF^{(n)}_{k-1} \big)
         \stoch 0$
        \ for all \ $\theta>0$,
 \end{enumerate}
 where \ $\stoch$ \ denotes convergence in probability.
Then \ $\bcU^{(n)} \distr \bcU$ \ as \ $n\to\infty$.
\end{Thm}

Note that in (i) of Theorem \ref{Conv2DiffThm}, \ $\|\cdot\|$ \ denotes
 a matrix norm, while in (ii) it denotes a vector norm.

\end{document}